\documentclass[reqno,12pt]{amsart}
\usepackage{amssymb}
\usepackage{euscript}

\newcommand{\ds}{\displaystyle}
\newcommand{\fr}[2]{\frac{#1}{#2}}
\newcommand{\dfr}[2]{\dfrac{#1}{#2}}%ams
\newcommand{\tfr}[2]{\tfrac{#1}{#2}}%ams
\newcommand{\cd}{\cdot}
\newcommand{\cds}{\cdots}
\newcommand{\dsum}{\displaystyle \sum}

\newcommand{\vsv}{\vspace{5mm}}
\newcommand{\vsb}{\vspace{2mm}}

\newcommand{\B}{{\mathcal B}}

\newcommand{\I}{{\mathcal I}}
\newcommand{\CP}{\mathcal{P}}

\renewcommand{\l}{\left}
\renewcommand{\r}{\right}

\newcommand{\q}{\quad}
\newcommand{\qq}{\qquad}

\theoremstyle{plain}
\newtheorem{theorem}{Theorem}[section]

\newtheorem{lemma}[theorem]{Lemma}
\newtheorem{proposition}[theorem]{Proposition}
\theoremstyle{remark}
\newtheorem{remark}[theorem]{Remark}

\theoremstyle{definition}

\numberwithin{equation}{section}
\newcommand{\hf}{\frac{1}{2}}
\newcommand{\al}{\alpha}
\newcommand{\be}{\beta}
\newcommand{\bbe}{\boldsymbol{\beta}}

\newcommand{\om}{\omega}
\newcommand{\tom}{\tilde{\omega}}

\newcommand{\ep}{\epsilon}
\newcommand{\tensor}{\otimes}
\newcommand{\la}{\langle}

\newcommand{\ra}{\rangle}

\newcommand{\Vir}{\mathrm{Vir}}
\newcommand{\Aut}{\mathrm{Aut}}

\newcommand{\Z}{\mathbb{Z}}

\newcommand{\C}{\mathbb{C}}
\newcommand{\R}{\mathbb{R}}

\newcommand{\Q}{\mathbb{Q}}

\newcommand{\ba}{\mathbf{a}}

\newcommand{\vir}{\mathrm{Vir}}
\newcommand{\aut}{\mathrm{Aut}}

\newcommand{\tr}{\mathrm{tr}}

\newcommand{\ch}{\mathrm{ch}}

\newcommand{\Span}{\mathrm{span}}

\renewcommand{\L}{\mathcal{L}}

\begin{document}

\baselineskip 6mm

\title[McKay's observation]
{McKay's observation and vertex operator algebras generated by two
conformal vectors of central charge $1/2$}

\author[C.H. Lam]{Ching Hung Lam $^\dagger$}
\address[C.H.  Lam ]{Department of Mathematics, National Cheng Kung University,
Tainan, Taiwan 701}%
\email{chlam@mail.ncku.edu.tw}

\author[H. Yamada]{Hiromichi Yamada $^\ddagger$}
\address[H. Yamada]{
Department of Mathematics, Hitotsubashi University, Kunitachi,
Tokyo 186-8601,  Japan} \email{yamada@math.hit-u.ac.jp}

\author[H. Yamauchi]{  Hiroshi Yamauchi $^\S$}
\address[H. Yamauchi]{Graduate School of Mathematical
Sciences, The University of Tokyo, Tokyo, 153-8914, Japan}

\email{yamauchi@ms.u-tokyo.ac.jp}

\thanks{$^\dagger$ Partially supported by NSC grant 93-2115-M-006-012 of Taiwan,
R.O.C.\\
\hspace*{3mm} $^\ddagger$ Partially supported by JSPS Grant-in-Aid for
Scientific Research No. 15540015
\\
\hspace*{3mm} $^\S$ Supported by JSPS Research Fellowships for Young Scientists. }

\subjclass{17B68, 17B69, 20D08}

\begin{abstract}
This paper is a continuation of \cite{lyy} at which several coset
subalgebras of the lattice VOA $V_{\sqrt{2}E_8}$ were constructed
and the relationship between such algebras with the famous McKay
observation on the extended $E_8$ diagram and  the Monster simple
group were discussed. In this article, we shall provide the
technical details. We completely determine the structure of the
coset subalgebras constructed and show that they are all generated
by two conformal vectors of central charge $1/2$. We also study
the representation theory of these coset subalgebras and show that
the product of two Miyamoto involutions is in the desired
conjugacy class of the Monster simple group if a coset subalgebra
$U$ is actually contained in the Moonshine VOA $V^\natural$. The
existence of $U$ inside the Moonshine VOA $V^\natural$ for the
cases of $1A,2A,2B$ and $4A$ is also established. Moreover, the
cases for $3A$, $5A$ and $3C$ are discussed.
\end{abstract}
\maketitle

\section{Introduction}

This paper is a continuation of the authors' work \cite{lyy} at
which several coset subalgebras of the lattice VOA
$V_{\sqrt{2}E_8}$ were constructed and the relationship between
such algebras with the famous McKay observation \cite{c1,McK} on
the extended $E_8$ diagram

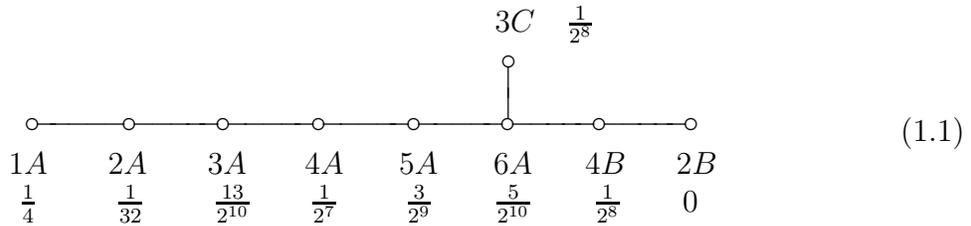
\begin{equation}\label{mckay}
\begin{array}{l}
  \hspace{184pt} 3C\quad   \frac{1}{2^8}\\
  \hspace{186.0pt}\circ \vspace{-6pt}\\
 \hspace{187.4pt}| \vspace{-6pt}\\
 \hspace{187.4pt}| \vspace{-6pt}\\
  \hspace{6pt} \circ\hspace{-5pt}-\hspace{-7pt}-\hspace{-5pt}-\hspace{-5pt}-
  \hspace{-5pt}-\hspace{-5pt}\circ\hspace{-5pt}-\hspace{-5pt}-
  \hspace{-5pt}-\hspace{-6pt}-\hspace{-7pt}-\hspace{-5pt}\circ
  \hspace{-5pt}-\hspace{-5.5pt}-\hspace{-5pt}-\hspace{-5pt}-
  \hspace{-7pt}-\hspace{-5pt}\circ\hspace{-5pt}-\hspace{-5.5pt}-
  \hspace{-5pt}-\hspace{-5pt}-\hspace{-7pt}-\hspace{-5pt}\circ
  \hspace{-5pt}-\hspace{-6pt}-\hspace{-5pt}-\hspace{-5pt}-
  \hspace{-7pt}-\hspace{-5pt}\circ\hspace{-5pt}-\hspace{-5pt}-
  \hspace{-6pt}-\hspace{-6pt}-\hspace{-7pt}-\hspace{-5pt}\circ
  \hspace{-5pt}-\hspace{-5pt}-
  \hspace{-6pt}-\hspace{-6pt}-\hspace{-7pt}-\hspace{-5pt}\circ
  \vspace{-6.2pt}\\
  \vspace{-6pt} \\
  1A\hspace{23pt} 2A\hspace{23  pt} 3A\hspace{22pt} 4A\hspace{21pt} 5A\hspace{21pt}
  6A\hspace{20pt} 4B\hspace{19pt} 2B\\
  \hspace{4pt} \frac{1}{4} \hspace{30pt} \frac{1}{32}\hspace{26  pt} \frac{13}{2^{10}}
  \hspace{22pt} \frac{1}{2^7} \hspace{25pt} \frac{3}{2^{9}}\hspace{23pt}
  \frac{5}{2^{10}}\hspace{23pt} \frac{1}{2^8}\hspace{23pt} 0 \\
  \
\end{array}
\end{equation}
and the Monster simple group were discussed.

In this article, we shall provide the technical details. We shall
determine the structure of the coset subalgebras and show that
they are all generated by two conformal vectors of central charge
$1/2$. We also study the representation theory of these coset
subalgebras and show that the product of two Miyamoto involutions
is in the desired conjugacy class of the Monster simple group if a
coset subalgebra $U$ is actually contained in the Moonshine VOA
$V^\natural$. The existence of $U$ inside the Moonshine VOA
$V^\natural$ for the cases of $1A,2A,2B$ and $4A$ is also
established. Moreover, the cases for $3A$, $5A$ and $3C$ are
discussed.

The organization of the article is as follows. In Section 2 we
shall review some important notation and terminology from
\cite{lyy}. We review certain conformal vectors in the lattice VOA
$V_{\sqrt{2}R}$, where $R$ is a root lattice of type $A$, $D$, or
$E$ (cf. \cite{dlmn}).  We then consider the sublattice $L$ of
$E_8$ and define the coset subalgebra $U$ and two conformal
vectors $\hat{e}$ and $\hat{f}$ of central charge $1/2$.  A
canonical automorphism $\sigma$ of order $n=|E_8/L|$ induced by
the quotient group $E_8/L$ is also discussed. In Section 3, we
study the structure of $U$ in each of the nine cases corresponding
to the McKay's diagram. We also study the representation theory of
these coset subalgebras and show that the product of two Miyamoto
involutions is in the desired conjugacy class of the Monster
simple group if a coset subalgebra $U$ is actually contained in
the Moonshine VOA $V^\natural$. The existence of $U$ inside the
Moonshine VOA $V^\natural$ for the cases of $1A,2A,2B$ and $4A$ is
also established. Moreover, the cases for $3A$, $5A$ and $3C$ are
discussed. Appendix contains the classification of conformal
vectors in $U$,  calculations of certain characters which are used
in Section 3 and the classification of irreducible modules for
$5A$ and $3C$ cases.

The authors thank Masahiko Miyamoto and Masaaki Kitazume for
stimulating discussions and Kazuhiro Yokoyama for helping them to
compute the conformal vectors for the cases of $5A$ and $6A$ by a
computer algebra system Risa/Asir. In Appendix C we study an
extension of a simple rational VOA by an irreducible module which is
not a simple current module. A similar extension is also considered
in \cite{TY}. The authors thank Kenichiro Tanabe for valuable
discussions concerning it. Part of the work was done while the third
author (H. Yamauchi) was visiting the National Center for
Theoretical Science of Taiwan in August, 2004. He thanks the center
for the hospitality during the stay.

\section{Preliminary}
In this section, we shall recall the notation and the constructions
of certain coset subalgebras of $V_{\sqrt{2}E_8}$ and their
automorphisms from \cite{lyy}. We shall mainly deal with lattice
VOAs introduced by \cite{flm}. Let $V_N$ be a lattice VOA associated
with any positive definite even lattice $N$. By \cite[Theorem
3.1]{Li}, there is a unique symmetric invariant bilinear form
$\la\,\cdot\,,\,\cdot\,\ra$ on $V_N$ such that
$\la\mathbf{1},\mathbf{1}\ra = 1$. That the form
$\la\,\cdot\,,\,\cdot\,\ra$ is invariant, i.e.,
\begin{equation*}
\la Y(u,z)v,w \ra = \la v, Y(e^{zL(1)}(-z^{-2})^{L(0)}u, z^{-1})w\ra
\quad \text{ for }u,v,w \in V_N.
\end{equation*}
For $u, v \in (V_N)_2$ with $L(1)u = 0$, the invariance of the form
implies that $\la u,v\ra = \la \mathbf{1}, u_3 v\ra$ which induces a
bilinear form on $V_2$. It is also well known that $V_N$ possesses a
positive definite invariant hermitian form $(\,\cdot\,,\,\cdot\,)$.
Indeed, let $V_{N,\R}$ be the $\R$-form of $V_N$ defined as in
\cite[Section 12.4]{flm}. Then $V_{N,\R}$ is invariant under the
automorphism $\theta$, where $\theta$ is a lift of $-1$ isometry of
the lattice $N$. Let $V_{N,\R}^{\pm}$ be the eigenspaces for
$\theta$ with eigenvalues $\pm 1$. Then $\la\,\cdot\,,\,\cdot\,\ra$
is positive definite on $V_{N,\R}^+$ and negative definite on
$V_{N,\R}^-$. Moreover, $\la V_{N,\R}^+, V_{N,\R}^-\ra = 0$. Hence
$\la\,\cdot\,,\,\cdot\,\ra$ is positive definite on the $\R$-vector
space $\widetilde{V}_{N,\R} = V_{N,\R}^+ + \sqrt{-1}V_{N,\R}^-$.
Clearly, $V_N = \C \otimes_{\R} \widetilde{V}_{N,\R}$ and so
$\widetilde{V}_{N,\R}$ is an $\R$-form of $V_N$. Define a hermitian
form $(\,\cdot\,,\,\cdot\,)$ on $V_N$ by $(\lambda u, \mu v) =
\lambda \overline{\mu}(u,v)$ for $\lambda, \mu \in \C$ and $u,v \in
\widetilde{V}_{N,\R}$. Then $(\,\cdot\,,\,\cdot\,)$ is positive
definite on $V_N$. Furthermore, it is
$\widetilde{V}_{N,\R}$-invariant, that is,
\begin{equation*}
(Y(u,z)v,w) = (v, Y(e^{zL(1)}(-z^{-2})^{L(0)}u, z^{-1})w)
\end{equation*}
for $u \in \widetilde{V}_{N,\R}$ and $v,w \in V_N$, where
$L(0)=\omega_1$ and $L(1)=\omega_2$ with $\omega$ being the
Virasoro element of $V_N$. These two forms
$\la\,\cdot\,,\,\cdot\,\ra$ and $(\,\cdot\,,\,\cdot\,)$ for the
case $N = \sqrt{2}E_8$ will be used in Section 2 and 3.

\subsection{Conformal vectors}

We shall now review the construction of certain conformal vectors in
the lattice VOA $V_{\sqrt{2}R}$ from \cite{dlmn}, where $R$ is a
root lattice of type $A_n$, $D_n$, or $E_n$.

Let $\Phi$ be the root system of $R$ and $\Phi^+$ and $\Phi^-$ the
set of all positive roots and negative roots, respectively. Then
$\Phi=\Phi^+\cup \Phi^-= \Phi^+\cup ( -\Phi^+)$. The Virasoro
element $\omega$ of $V_{\sqrt{2}R}$ is given by
\[
\omega =\om (\Phi)=\frac{1}{2h}\sum_{\al\in \Phi^+} \al(-1)^2\cdot
1,
\]
where $h$ is the Coxeter number of $\Phi$. Now define
\begin{equation}\label{cv1}
\begin{split}
s=s(\Phi)&=\frac{1}{2(h+2)}\sum_{\al\in \Phi^+}\left(
\al(-1)^2\cdot
1 -2(e^{\sqrt{2}\al}+ e^{-\sqrt{2}\al})\right),\\
\tilde{\om}=\tilde{\om}(\Phi)&= \om - s.
\end{split}
\end{equation}
It is shown in \cite{dlmn} that $\tilde{\om}$ and $s$ are mutually
orthogonal conformal vectors and the central charge of
$\tilde{\om}$ is $2n/(n+3)$ if $R$ is of type $ A_n$, $1$ if $R$
is of type $ D_n$ and $6/7, 7/10$ and $1/2$ if $R$ is of type
$E_6, E_7$ and $E_8$, respectively.
In this article, we denote by $\vir(x)$ the Virasoro sub VOA generated by
a conformal vector $x$ of $V$.

Now let $\mathrm{Aut\,}(\Phi)$ be the automorphism group of $\Phi$.
For any element $g\in \mathrm{Aut\,}(\Phi)$, $g$ induces an
automorphism on the lattice $R$ and hence also defines an
automorphism of the VOA $V_{\sqrt{2}R}$ by
\[
g(u\otimes e^{\sqrt{2}\al})= gu\otimes e^{\sqrt{2} g\al} \qquad
\text{ for } \quad u\otimes e^{\sqrt{2}\al}\in M(1)\otimes
e^{\sqrt{2}\al}\subset V_{\sqrt{2}R}.
\]
Note that both $s$ and $\tilde{\om}$ are fixed by
$\mathrm{Aut\,}(\Phi)$ and thus also fixed by the Weyl group
$W(\Phi)$ of $\Phi$.
%In fact, $\om,\tom$ and $s$ are the only conformal
%vectors in $V_{\sqrt{2}R}$ which are fixed by the Weyl group.

Let $R^* = \{ \al \in \Q \otimes_{\Z} R\,|\,\la\al,R\ra \subset
\Z\}$ be the dual lattice of $R$. Then we have the following
Proposition.

\begin{proposition}[cf. Section 2 of \cite{lyy}]\label{hw}
Let $\gamma+R$ be a coset of $R$ in $R^*$ and $k=\min\{ \la
\al,\al\ra | \al\in \gamma+ R\}$. Define
\[
v=\sum_{{\al\in \gamma+R}\atop \la \al,\al\ra=k} e^{\sqrt{2}\al}
\in V_{\sqrt{2}(\gamma+R)}.
\]
Then $v$ is a highest weight vector of highest weight $(0,k)$ in
$V_{\sqrt{2}(\gamma+ R)}$ with respect to  $\Vir(s)\otimes
\Vir(\tilde{\om})$, that is, $s_j v = \tilde{\om}_j v=0$ for all
$j \ge 2$, $s_1 v=0$, and $\tilde{\om}_1 v=kv$. In other words,
with respect to $s$, there is always a highest weight vector of
weight $0$ in $V_{\sqrt{2}(\gamma+R)}$.
\end{proposition}

\subsection{ Extended $E_8$ diagram and coset subalgebras of $V_{\sqrt{2}E_8}$}\label{S3}
Next, we shall review the construction some coset VOAs $U$ using
the extended $E_8$ diagram. In each case, $U$ contains some
conformal vectors of central charge $1/2$ and the inner products
among these conformal vectors are the same as the numbers given in
the McKay diagram \ref{mckay}. First, we shall consider certain
sublattices of the root lattice $E_8$ by using the extended $E_8$
diagram
\begin{equation}\label{ext}
\begin{array}{l}
  \hspace{183pt} \al_8\\
  \hspace{186.0pt}\circ \vspace{-6pt}\\
 \hspace{187.4pt}| \vspace{-6pt}\\
 \hspace{187.4pt}| \vspace{-6pt}\\
  \hspace{6pt} \circ\hspace{-5pt}-\hspace{-7pt}-\hspace{-5pt}-\hspace{-5pt}-
  \hspace{-5pt}-\hspace{-5pt}\circ\hspace{-5pt}-\hspace{-5pt}-
  \hspace{-5pt}-\hspace{-6pt}-\hspace{-7pt}-\hspace{-5pt}\circ
  \hspace{-5pt}-\hspace{-5.5pt}-\hspace{-5pt}-\hspace{-5pt}-
  \hspace{-7pt}-\hspace{-5pt}\circ\hspace{-5pt}-\hspace{-5.5pt}-
  \hspace{-5pt}-\hspace{-5pt}-\hspace{-7pt}-\hspace{-5pt}\circ
  \hspace{-5pt}-\hspace{-6pt}-\hspace{-5pt}-\hspace{-5pt}-
  \hspace{-7pt}-\hspace{-5pt}\circ\hspace{-5pt}-\hspace{-5pt}-
  \hspace{-6pt}-\hspace{-6pt}-\hspace{-7pt}-\hspace{-5pt}\circ
  \hspace{-5pt}-\hspace{-5pt}-
  \hspace{-6pt}-\hspace{-6pt}-\hspace{-7pt}-\hspace{-5pt}\circ
  \vspace{-6.2pt}\\
  \vspace{-6pt} \\
  \hspace{6pt}\al_0\qq \alpha_1\qq \alpha_2\qq \alpha_3\qq
  \alpha_4\qq \alpha_5\qq \alpha_6 \qq \alpha_7
\end{array}
\end{equation}
where $\al_1,\al_2,\dots,\al_8$ are the simple roots of $E_8$ and
\begin{equation}\label{eq:A_0}
\al_0+2\al_1+3\al_2+4\al_3+ 5\al_4+6\al_5+4\al_6+2\al_7+3\al_8 =
0.
\end{equation}
Then $\la \al_i,\al_i\ra = 2$, $0 \le i \le 8$. Moreover, for
$i\ne j$, $\la \al_i, \al_j \ra = -1$ if the nodes $\al_i$ and
$\al_j$ are connected by an edge and $\la \al_i, \al_j \ra = 0$
otherwise. Note that $-\alpha_0$ is the highest root.

For any $i=0,1,\dots, 8$, let $L(i)$ be the sublattice generated
by $\al_j, 0 \le j \le 8, j \ne i$. Then $L(i)$ is a rank $8$
sublattice of $E_8$. Note $L(i)$ is the lattice associated with
the Dynkin diagram obtained by removing the corresponding node
$\al_i$ in the extended $E_8$ diagram \eqref{ext} and the index $|
E_8/L(i) |$ is equal to $n_i$, where $n_i$ is the coefficient of
$\al_i$ in the left hand side of \eqref{eq:A_0}. Actually,
\begin{gather}\label{eq:DEC}
\begin{aligned}
    L(0)&\cong E_8,\\
    L(3)&\cong A_3\oplus D_5,\\
    L(6)&\cong A_7\oplus A_1,
\end{aligned}
\qquad
\begin{aligned}
    L(1)&\cong A_1\oplus E_7,\\
    L(4)&\cong A_4\oplus A_4,\\
    L(7)&\cong D_8,
\end{aligned}
\qquad
\begin{aligned}
    L(2)&\cong A_2\oplus E_6, \\
    L(5)&\cong A_5\oplus A_2\oplus A_1,\\
    L(8)&\cong A_8.
\end{aligned}
\end{gather}
\medskip

Now let us explain the details of our construction. First, we fix
$i\in \{0,1,\dots,8\}$ and denote $L(i)$ by $L$. In each case,
$|E_8/L|=n_i$ and $\al_i+L$ is a generator of the quotient group
$E_8/L$. Hence we have
\begin{equation}\label{eq:E8L}
E_8=L\cup (\al_i+L)\cup (2\al_i+L)\cup \cdots \cup (
(n_i-1)\al_i+L).
\end{equation}
Let $\lambda=\sqrt{2}\al_i$. Then
\[
\sqrt{2}E_8=\sqrt{2}L\cup (\lambda+\sqrt{2}L)\cup
(2\lambda+\sqrt{2}L) \cup \cdots \cup ( (n_i-1)\lambda+\sqrt{2}L)
\]
and the lattice VOA $V_{\sqrt{2}E_8}$ can  be decomposed as
\[
V_{\sqrt{2}E_8}=V_{\sqrt{2}L}\oplus V_{\lambda+\sqrt{2}L}\oplus
\cdots\oplus  V_{(n_i-1)\lambda+\sqrt{2}L},
\]
where $V_{j\lambda+\sqrt{2}L}$, $j=0,1,\dots,n_i-1$, are
irreducible modules of $V_{\sqrt{2}L}$ (cf. \cite{d1}).

\begin{remark}\label{sigma}
The abelian group $E_8/L$ actually induces an automorphism
$\sigma$ of $ V_{\sqrt{2}E_8}$ such that
\begin{equation}\label{au}
\sigma(u)=\xi^j u\qquad \text{ for any } \quad u\in
V_{j\lambda+\sqrt{2}L},
\end{equation}
where $\xi=e^{2\pi \sqrt{-1}/ n_{i}}$ is a primitive $n_i$-th root
of unity. More precisely, let
\begin{equation}\label{adef}
\ba = \begin{cases}
\al_1 & \mbox{if } i=0,\\
-\frac{1}{i+1}(\al_0 + 2\al_1 + \cdots + i\al_{i-1}) &
\mbox{if } 1\le i\le 5,\\
-\frac{1}{8}(\al_0+2\al_1 + \cdots + 6\al_5+7\al_8) &
\mbox{if } i=6,\\
\frac{1}{2}(\al_6+\al_8) &
\mbox{if } i=7,\\
-\frac{1}{9}(\al_0+2\al_1 + \cdots + 8\al_7) & \mbox{if } i=8.
\end{cases}
\end{equation}
Then $\la\ba,\al_j\ra \in \Z$ for $0 \le j\le 8$ with $j\ne i$ and
$\la\ba,\al_i\ra \equiv -1/n_i \pmod \Z$. The automorphism
$\sigma: V_{\sqrt{2}E_8} \to V_{\sqrt{2}E_8}$ is in fact defined
by
\begin{equation}\label{sigmadef}
\sigma = e^{-\pi \sqrt{-1} \bbe(0)} =\exp(-\sqrt{-2}\pi\ba(0))
\qquad \text{ with } \quad \bbe=\sqrt{2} \ba.
\end{equation}
For $u\in M(1)\otimes e^{\al} \subset V_{\sqrt{2}E_8}$, we have
$\sigma(u)=e^{-\pi \sqrt{-1} \la\bbe, \al\ra} u$. Note that $\ba
+R$ is a generator of the group $R^*/R$ for the cases $i\ne 0,7$,
where $R$ is an indecomposable component of the lattice $L$ of
type $A$.

For any lattice VOA $V_N$ associated with a positive definite even
lattice $N$, there is a natural involution $\theta$ induced by the
isometry $\al \to -\al$ for $\al \in N$. If $N = \sqrt{2}E_8$,
which is doubly even, we may define $\theta: V_{\sqrt{2}E_8} \to
V_{\sqrt{2}E_8}$ by
\[
\alpha(-n)\rightarrow -\alpha(-n)\qquad \mbox{ and }\qquad
e^{\alpha}\rightarrow e^{-\alpha}
\]
for any $\al \in \sqrt{2}E_8$ (cf. \cite{flm}). Then $\theta
\sigma \theta= \sigma^{-1}$ and the group generated by $\theta$
and $\sigma$ is isomorphic to a dihedral group of order $2n_i$.
\end{remark}

Next we shall recall the definition of certain coset subalgebras
from \cite{lyy}. Let $R_1,\dots, R_l$ be the indecomposable
components of the lattice $L$ and $\Phi_1,\dots, \Phi_l$ the
corresponding root systems of $R_1,\dots, R_l$ (cf.
\eqref{eq:DEC}). Then $L=R_1\oplus \cdots\oplus R_l$ and
\begin{equation*}
V_{\sqrt{2}L}\cong V_{\sqrt{2}R_1}\otimes \cdots \otimes
V_{\sqrt{2}R_l}.
\end{equation*}

By \eqref{cv1} in Section 2, one obtains $2l$ mutually orthogonal
conformal vectors
\begin{equation}\label{eq:CV}
s^k=s(\Phi_k),\quad \tilde{\om}^k=\tilde{\om}(\Phi_k),\quad 1 \le
k \le l
\end{equation}
such that the Virasoro element $\om$ of $V_{\sqrt{2}L}$, which is
also the Virasoro element of $V_{\sqrt{2}E_8}$, can be written as
a sum of these conformal vectors
\begin{equation*}
\om=s^1+\cdots+ s^l +\tilde{\om}^1+\cdots+ \tilde{\om}^l.
\end{equation*}

Now we define $U$ to be a coset (or commutant) subalgebra
\begin{equation}\label{eq:UDEF}
U=\{v\in V_{\sqrt{2}E_8}\,| \, (s^k)_1 v=0 \text{ for all } k=1,
\dots, l\}.
\end{equation}
Note that $U$ is a VOA with the Virasoro element
$\om'=\tilde{\om}^1+\cdots+ \tilde{\om}^l$ and the automorphism
$\sigma$ defined by (\ref{au}) induces an automorphism of order
$n_i$ on $U$. By abuse of notation, we shall denote it by $\sigma$
also.

\begin{remark}
In \cite{ly3}, it is shown that $\{ v\in V_{\sqrt{2}A_n}\,|\,
s(A_n)_1 v=0\}$ is isomorphic to a parafermion algebra
$W_{n+1}(2n/(n+3))$ of central charge $2n/(n+3)$. Thus, if $L$ has
some indecomposable component of type $A_n$, then $U$ will contain
some subalgebra isomorphic to a parafermion algebra. It is well
known \cite{zf1} that the parafermion algebra $W_{n+1}(2n/(n+3))$
processes a certain $\Z_{n+1}$ symmetry among its irreducible
modules. The automorphism $\sigma$ is actually related to such a
symmetry. More details about the relation between coset subalgebra
$U$ and the parafermion algebra $W_{n+1}(2n/(n+3))$ can be found
in Appendix $B$.
\end{remark}

Next we shall recall the definition of two conformal vectors of
central charge $1/2$ from \cite{lyy}. Note that
\begin{equation}\label{eq:EF}
\begin{aligned}
\hat{e} &= \frac{1}{16}\om +\frac{1}{32}\sum_{\al\in\Phi^+(E_8)}
(e^{\sqrt{2}\al}+e^{-\sqrt{2}\al})
\end{aligned}
\end{equation}
is a conformal vector of central charge $1/2$. Let $\sigma$ be the
automorphism defined in Remark \ref{sigma}. Then we have

\begin{theorem}[cf. \cite{lyy}]
Let $\hat{e}$ be defined as above and $\hat{f} = \sigma \hat{e}$.
Then $\hat{e}, \hat{f}\in U$ and we have
\begin{equation}
\la \hat{e}, \hat{f}\ra = {\displaystyle
\begin{cases}
1/4 & \text{ if } i=0,\\
1/32 & \text{ if } i=1,\\
13/2^{10} & \text{ if } i=2,\\
1/2^{7} & \text{ if } i=3,\\
3/2^{9} & \text{ if } i=4,\\
5/2^{10} & \text{ if } i=5,\\
1/2^{8} & \text{ if } i=6,\\
0               & \text{ if } i=7,\\
1/2^{8} & \text{ if } i=8.
\end{cases}}
\end{equation}
In other words, the values of $\la \hat{e}, \hat{f}\ra$ are
exactly the values given in McKay's diagram \eqref{mckay}.
\end{theorem}

\section{The coset subalgebra $U$ and Miyamoto's
$\tau$-involutions}\label{s5}

This section is the main part of this article. We shall study the
structure of the coset subalgebra $U$ defined by \eqref{eq:UDEF}
for each of the nine cases. Except for the case of $4A$, we shall
show that the subalgebra $U$ always contain a set of mutually
orthogonal conformal vectors such that their sum is the Virasoro
element of $U$ and the central charges of these conformal vectors
are all coming from the unitary series
\begin{equation}\label{eq:VIRC}
c=c_m=1-\frac{6}{(m+2)(m+3)}, \qquad m=1,2,3, \dots.
\end{equation}

Such a conformal vector generates a simple Virasoro VOA isomorphic to
$L(c_m,0)$ inside $U$. The irreducible modules of $L(c_m,0)$ are of the form
$L(c_m,h^m_{r,s})$, where
\begin{equation}\label{eq:VIRH}
h^m_{r,s} = \frac{(r(m+3)-s(m+2))^2 - 1}{4(m+2)(m+3)},\qquad 1 \le
r \le m+1,\quad 1 \le s \le m+2.
\end{equation}
Note that $h^m_{r,s} = h^m_{m+2-r,m+3-s}$ and that
$L(c_m,h^m_{r,s})$, $1 \le s \le r \le m+1$ are all the
inequivalent irreducible $L(c_m,0)$-modules.

For the $4A$ case, we shall show that $U$ is isomorphic to the
fixed point subalgebra $V_\mathcal{N}^+$ of $\theta$ for some rank
two lattice $\mathcal{N}$.

Furthermore, we shall discuss the relation between Miyamoto's
$\tau$-involutions and the structure of the coset subalgebra $U$.
We shall show that for any VOA $V$ which contains a subalgebra
isomorphic to $U$, the product of the Miyamoto involutions
$\tau_{\hat{e}}$ and $\tau_{\hat{f}}$ naturally defines an
automorphism of order $n_i$ or $n_i/2$ on $V$. If $U$ is actually
contained in the Moonshine VOA $V^\natural$, then we shall show
that  $\tau_{\hat{e}}\tau_{\hat{f}}$ is of the desired conjugacy
class of the Monster simple group mentioned in the McKay diagram.
The existence of $U$ inside the Moonshine VOA $V^\natural$ will
also be established for the cases $1A,2A,2B,3A$, and $4A$.

As in Section \ref{S3}, $L=L(i), i=0,1,\ldots,8$ denotes the
lattice associated with the Dynkin diagram obtained by removing
the $i$-th node $\al_i$ in the extended $E_8$ diagram. The coset
decomposition of $E_8$ by $L$ is given in \eqref{eq:E8L}. For
$j=1, \dots, n_i-1$, we define
\begin{equation}\label{eq:HWV}
X^j = \ds\dsum_{\alpha \in j\al_i + L \atop \la \alpha, \alpha \ra
=2} e^{\sqrt{2}\alpha}.
\end{equation}
Clearly $X^j$ is of weight $2$. By Proposition \ref{hw}, it is
easy to see that $X^j\in U$ for all $j=1, \dots, n_i-1$.

Recall the positive definite invariant hermitian form
$(\,\cdot\,,\,\cdot\,)$ on lattice VOAs mentioned in Introduction.
We shall consider the form $(\,\cdot\,,\,\cdot\,)$ for the lattice
VOA $V_{\sqrt{2}E_8}$. Let
\begin{equation*}
Y^{j,+} = \frac{1}{2}(X^j + X^{n_i -j}), \quad
Y^{j,-} = \frac{1}{2}\sqrt{-1}(X^j - X^{n_i -j}), \quad
1 \le j \le [n_i/2]
\end{equation*}
and set $\B_{\R}=\mathrm{span}_{\R} \{ \tilde{\omega}^1, \ldots,
\tilde{\omega}^l, Y^{1,\pm},\ldots, Y^{[n_i/2], \pm} \}$. Then
$\B_{\R}$ is contained in $\widetilde{V}_{\sqrt{2}E_8,\R}$ and so
the form $(\,\cdot\,,\,\cdot\,)$ is $\B_{\R}$-invariant. It is
clear that the conformal vectors $\hat{e}$ and $\hat{f}$ defined
by \eqref{eq:EF} are contained in $\B_{\R}$. In the following
argument we use the fact that $V_{\sqrt{2}E_8}$ possesses a
positive definite hermitian form which is $\B_{\R}$-invariant.

\medskip
Now let us study the structure of $U$ and the Miyamoto involutions
$\tau_{\hat{e}}$ and $\tau_{\hat{f}}$ associated with the
conformal vectors $\hat{e}$ and $\hat{f}$ defined by \eqref{eq:EF}
in each of the nine cases. First we shall note that
\begin{equation} \label{tau-ef}
\tau_{\hat{e}}\tau_{\hat{f}} = e^{2\pi \sqrt{-1} \bbe(0)}
\end{equation}
as an automorphism of $V_{\sqrt{2}E_8}$. (cf. \cite[Section
4]{lyy}).

\subsection{$1A$ case} In this case, $L=L(0)\cong E_8$,
$n_0=|E_8/L|=1$, and $U\cong L(1/2,0)$. The conformal vector
$\tilde{\om}^1 \in V_{\sqrt{2}E_8}$ defined by \eqref{eq:CV} is
the only conformal vector in $U$. Its central charge is $1/2$.
Moreover, $\hat{e}=\hat{f}=\tilde{\om}^1$ and $U=\Vir(\hat{e})
\subset V^+_\Lambda\subset V^\natural$. Thus we have
$\tau_{\hat{e}}\tau_{\hat{f}}=1$.

\subsection{$2A$ case} In this case, $ L=L(1)\cong A_1\oplus E_7$,
$n_1=|E_8/L|=2$, and the conformal vectors $\tom^1 \in
V_{\sqrt{2}A_1}$ and $\tom^2 \in V_{\sqrt{2}E_7}$ defined by
\eqref{eq:CV} are of central charge $1/2$ and $7/10$,
respectively.

\begin{proposition}
The vector $X^1$ defined by \eqref{eq:HWV} is a highest weight
vector of highest weight $(1/2, 3/2)$ with respect to
$\Vir(\tom^1)\otimes \Vir(\tom^2)$ . Thus as a module of
$\Vir(\tom^1)\otimes \Vir(\tom^2)$,
\[
U\cong L(\frac{1}2,0)\otimes L(\frac{7}{10},0) \oplus
L(\frac{1}2,\frac{1}2)\otimes L(\frac{7}{10},\frac{3}{2}).
\]
\end{proposition}

This VOA has been well studied in \cite{kly,ly}. In fact, $U$ has
exactly three conformal vectors of central charge 1/2, namely
$\hat{e}$, $\hat{f}$ and $w=\tom^1$. The automorphism group
$\Aut\, U$ of $U$ is a symmetric group $S_3$ of degree $3$. Note
that $U$ is generated by $\hat{e}$ and $\hat{f}$ and they are both
fixed by $\theta$. Thus, $U \subset V^+_\Lambda\subset
V^\natural$.

The VOA $U$ is rational and it has exactly eight inequivalent
irreducible modules $M^j, W^j, j \in \{0,a,b,c\}\cong \Z_2\times
\Z_2$. As $\Vir(\tom^1)\otimes \Vir(\tom^2)$-modules, they are of
the following form,
\begin{align*}
M^{0} & \cong [0,0] \oplus [\frac{1}{2},\frac{3}{2}], & W^{0} &
\cong [0,\frac{3}{5}] \oplus [\frac{1}{2},
\frac{1}{10}],\\
M^{a} & \cong M^{b} \cong [\frac{1}{16},\frac{7}{16}], & W^{a} &
\cong W^{b} \cong [\frac{1}{16},\frac{3}{80}],\\
M^{c} & \cong [\frac{1}{2},0] \oplus [0,\frac{3}{2}], & W^{c} &
\cong [\frac{1}{2},\frac{3}{5}] \oplus [0,\frac{1}{10}],
\end{align*}
where $[h_1,h_2]$ denotes $L(1/2,h_1) \otimes L(7/10,h_2)$.

It is known that the fusion rules among irreducible $U$-modules
have a symmetry of $\Z_2\times \Z_2$. For any VOA $V$ containing a
subalgebra isomorphic to $U$, there are three automorphisms of
order $2$ or $1$ associated with $U$ (cf. \cite{ly}). They are
given by
\[
\tau^a=
\begin{cases}
1 &\text{ on } U^0,U^a,\\
-1& \text{ on } U^b,U^c,
\end{cases}
\qquad \tau^b=
\begin{cases}
1 &\text{ on } U^0,U^b,\\
-1& \text{ on } U^c,U^a,
\end{cases}
\qquad \tau^c=
\begin{cases}
1 &\text{ on } U^0,U^c,\\
-1& \text{ on } U^a,U^b,
\end{cases}
\]
where $U^j$ is the sum of all irreducible $U$-submodules of $V$
which are isomorphic to either $M^j$ or $W^j$ for $j=0,a,b,c$.
Actually,
\[
\tau^a=\tau_w,\quad \tau^b=\tau_{\hat{e}},\quad
\tau^c=\tau_{\hat{f}},
\]
and we have
$\tau_{\hat{e}}\tau_{\hat{f}}=\tau^b\tau^c=\tau^a=\tau_w$.

If $V$ is the Moonshine VOA $V^\natural$, then we have the
following theorem.

\begin{theorem}
As automorphisms of $V^{\natural}$,
$\tau_{\hat{e}}\tau_{\hat{f}}=\tau_w$ and thus
$\tau_{\hat{e}}\tau_{\hat{f}}$ is of class $2A$.
\end{theorem}

\subsection{$3A$ case}
In this case, $L=L(2)\cong A_2\oplus E_6$,
$n_2=|E_8/L|=3$, and the conformal vectors $\tom^1 \in
V_{\sqrt{2}A_2}$ and $\tom^2 \in V_{\sqrt{2}E_6}$ defined by
\eqref{eq:CV} are of central charge $4/5$ and $6/7$, respectively.
>From \cite{kmy,ly2}, we know that
\begin{align*}
U^1 &=\{u\in V_{\sqrt{2}A_2}| (s^1)_1u=0\}\cong L(\frac{4}{5},
0)\oplus L(\frac{4}{5}, 3),\\
U^2 &=\{u\in V_{\sqrt{2}E_6}| (s^2)_1u=0\}\cong L(\frac{6}{7},
0)\oplus L(\frac{6}{7}, 5).
\end{align*}
Hence,
\[
U\supset U^1\otimes U^2\cong \Big(L(\frac{4}{5}, 0)\oplus
L(\frac{4}{5}, 3)\Big)\otimes \Big(L(\frac{6}{7}, 0)\oplus
L(\frac{6}{7}, 5)\Big).
\]

\begin{proposition}
Both of the vectors $X^1$ and $X^2$ defined by \eqref{eq:HWV} are
highest weight vectors of highest weight $(2/3,4/3)$ with respect
to $\Vir(\tom^1)\otimes \Vir(\tom^2)$.
\end{proposition}

\begin{proof}
First, note that $X^1$ and $X^2 $ are highest weight vectors of
weight $2$ with respect to the Virasoro element $\om'
=\tom^1+\tom^2$ of $U$ and that
\begin{equation*}
\begin{split}
\tom^1 &= \frac{1}{15} \sum_{\al\in\Phi^+(A_2)} \al(-1)^2 \cdot 1
+ \frac{1}5 \sum_{\al\in\Phi^+(A_2)} (e^{\sqrt{2}\al} +
e^{-\sqrt{2}\al})\\
&=\frac{2}{5} \om' + \frac{1}5 \sum_{\al\in\Phi^+(A_2)}
(e^{\sqrt{2}\al} + e^{-\sqrt{2}\al}).
\end{split}
\end{equation*}
Clearly, $(\tom^1)_3X^1= (\tom^1)_2X^1 =0$ and
\begin{equation*}
\begin{split}
(\tom^1)_1 X^1 &= \Big(\frac{2}{5} \om' + \frac{1}5
\sum_{\al\in\Phi^+(A_2)}
(e^{\sqrt{2}\al} + e^{-\sqrt{2}\al})\Big)_1 X^1 \\
& = (\frac{2}{5}\times \frac{2}3+ \frac{1}5\times 2) X^1 =
\frac{2}3 X^1.
\end{split}
\end{equation*}
Hence, $X^1$ is a highest weight vector of highest weight
$(2/3,4/3)$ with respect to $\Vir(\tom^1)\otimes \Vir(\tom^2)$.
Similarly, $X^2$ is also a highest weight vector of highest weight
$(2/3,4/3)$ with respect to $\Vir(\tom^1)\otimes \Vir(\tom^2)$.
\end{proof}

Since  $U^1\otimes U^2$ and $L(4/5,2/3)^{\pm}\otimes L(6/7,
4/3)^{\pm}$ are the only irreducible modules of $U^1\otimes U^2$
which have integral weights (cf. \cite{kmy,lly,ly2}), by comparing
dimensions of the homogeneous subspaces of small weights, we have
the following proposition.
\begin{proposition}
As a module of $\Vir(\tom^1)\otimes \Vir(\tom^2)$,
\begin{equation*}
\begin{split}
U \cong & \Big(L(\frac{4}{5}, 0)\oplus L(\frac{4}{5},
3)\Big)\otimes
\Big(L(\frac{6}{7}, 0) \oplus L(\frac{6}{7}, 5)\Big)\\
& \oplus L(\frac{4}{5}, \frac{2}{3})\otimes L(\frac{6}{7},
\frac{4}{3})\oplus L(\frac{4}{5}, \frac{2}{3})\otimes
L(\frac{6}{7}, \frac{4}{3}).
\end{split}
\end{equation*}
\end{proposition}

By using Appendix A, we know that there are four distinct pairs of
mutually orthogonal conformal vectors of central charge $4/5$ and
$6/7$, namely, $\{\tom^1, \tom^2\}$, $\{x^j, y^j\}$, $j=0,1,2$,
where
\begin{equation*}
x^j = \frac{1}{16}\tom^1 + \frac{7}{8}\tom^2 - \frac{1}{48}\xi^j
X^1 -\frac{1}{48}\xi^{2j} X^2,\quad y^j = \frac{15}{16}\tom^1 +
\frac{1}{8}\tom^2 + \frac{1}{48}\xi^j X^1 +
\frac{1}{48}\xi^{2j}X^2
\end{equation*}
and $\xi=e^{2\pi \sqrt{-1}/3}$ is a primitive cubic root of unity.

\begin{lemma}\label{hwv-3a}
Let $u=135\tom^1-126\tom^2-13(X^1+X^2)$ and $v=X^1-X^2$. Then $u$
and $v$ are highest weight vectors of highest weight $(2/3,4/3)$
and $(13/8, 3/8)$ with respect to $\Vir(x^0)\otimes \Vir(y^0)$,
respectively.
\end{lemma}

\begin{proof}
We have
\[
\begin{split}
(x^0)_1 X^1 &= (\frac{1}{16}\tom^1+\frac{7}8\tom^2
-\frac{1}{48}X^1
  -\frac{1}{48}X^2)_1 X^1\\
  & = (\frac{1}{16} \times \frac{2}3 + \frac{7}8 \times \frac{4}3) -
  \frac{20}{48} X^2 - \frac{1}{48} ( 135\tom^1+252\tom^2)\\
  &= \frac{29}{24} X^1 - \frac{5}{12} X^2 - \frac{1}{16} (
  45\tom^1+84\tom^2).
\end{split}
\]
Similarly,
\[
\begin{split}
(x^0)_1 X^2&= \frac{29}{24} X^2 - \frac{5}{12} X^1 - \frac{1}{16}
(
  45\tom^1+84\tom^2).
\end{split}
\]
Thus $(x^0)_1v = \frac{13}{8} v$ and hence $v$ is a highest weight
vector of highest weight $(13/8, 3/8)$. Furthermore,

\begin{equation*}
\begin{split}
(x^0)_1 u = &\  135 \Big( \frac{1}{8} \tom^1 -\frac{1}{48} \times
\frac{2}3 (X^1+X^2)\Big) - 126 \Big( \frac{7}{4} \tom^2
-\frac{1}{48} \times \frac{4}3 (X^1+X^2)\Big)\\
& \ - 13 \Big( \frac{19}{24} (X^1+X^2)-\frac{1}8(
45\tom^1+84\tom^2)\Big)\\
=&\ \frac{2}{3} \left( 135\tom^1 -126\tom^2-13(X^1+X^2)\right)=
\frac{2}3 u.
\end{split}
\end{equation*}
Thus $u$ is a highest weight vector of highest weight $(2/3,
4/3)$.
\end{proof}

\begin{proposition}
As a module of $\Vir(x^0)\otimes \Vir(y^0)$,
\[
\begin{split}
U\cong & L(\frac{4}{5}, 0)\otimes L(\frac{6}{7}, 0) \oplus
L(\frac{4}{5}, 3)\otimes L(\frac{6}{7}, 5) \oplus L(\frac{4}{5},
\frac{2}{3})\otimes L(\frac{6}{7}, \frac{4}{3})\\
& \oplus L(\frac{4}{5}, \frac{13}{8})\otimes L(\frac{6}{7},
\frac{3}{8})\oplus L(\frac{4}{5}, \frac{1}{8})\otimes
L(\frac{6}{7}, \frac{23}{8}).
\end{split}
\]
\end{proposition}

\begin{proof}
The fixed point subalgebra $U^+$ of $\theta$ contains $\Vir(x^0)
\otimes \Vir(y^0)$. Moreover, $u \in U^+$ and so $U^+$ contains a
submodule of the form $L(4/5, 2/3)\otimes L(6/7, 4/3)$ by Lemma
\ref{hwv-3a}. Hence comparing the first several terms of the
characters, we know that
\[
U^+\cong L(\frac{4}{5}, 0)\otimes L(\frac{6}{7}, 0) \oplus
L(\frac{4}{5}, 3)\otimes L(\frac{6}{7}, 5) \oplus L(\frac{4}{5},
\frac{2}{3})\otimes L(\frac{6}{7}, \frac{4}{3})
\]
as a module of $\Vir(x^0)\otimes \Vir(y^0)$. Note that $U^+$ has
the same form as a module of $\Vir(\tom^1)\otimes \Vir(\tom^2)$.
Lemma \ref{hwv-3a} also implies that $U^-$ contains a submodule of
the form $L(4/5, 13/8)\otimes L(6/7, 3/8)$ since $v\in U^-$. Thus
\[
U^-\cong L(\frac{4}{5}, \frac{13}{8})\otimes L(\frac{6}{7},
\frac{3}{8})\oplus L(\frac{4}{5}, \frac{1}{8})\otimes
L(\frac{6}{7}, \frac{23}{8})
\]
and we have the desired result.
\end{proof}

The VOA $U$ has been constructed and studied by Sakuma and
Yamauchi \cite{sy} (see also Miyamoto \cite {m6}). It is known
that the automorphism group $\Aut\,U$ of $U$ is isomorphic to the
symmetric group $S_3$ and $U$ is generated by two conformal
vectors of central charge 1/2, namely, $\hat{e}$ and $\hat{f}$.

There are exactly six irreducible modules of $U$ (cf. \cite{sy}),
namely,
\begin{align*}
& M^0\otimes A^0 \oplus M^1\otimes A^1 \oplus M^2\otimes A^2,&
& M^0\otimes B^0 \oplus M^1\otimes B^1 \oplus M^2\otimes B^2,\\
& M^0\otimes C^0 \oplus M^1\otimes C^1 \oplus M^2\otimes C^2,&
& W^0\otimes A^0 \oplus W^1\otimes A^1 \oplus W^2\otimes A^2,\\
& W^0\otimes B^0 \oplus W^1\otimes B^1 \oplus W^2\otimes B^2,& &
W^0\otimes C^0 \oplus W^1\otimes C^1 \oplus W^2\otimes C^2,
\end{align*}
where
\begin{align*}
& M^{0}\cong L( \frac{4}{5},0) \oplus L( \frac{4}{5} ,3),& &
M^{1}\cong L( \frac{4}{5},\frac{2}{3})^+,& & M^{2}\cong L(
\frac{4}{5},\frac{2}{3})^-,\\
& W^{0}\cong L( \frac{4}{5},\frac{2}{5}) \oplus L(
\frac{4}{5},\frac{7}{5}),& & W^{1}\cong L(
\frac{4}{5},\frac{1}{15})^+, & & W^{2}\cong L(
\frac{4}{5},\frac{1}{15})^-
\end{align*}
are the irreducible modules of $L(4/5,0)\oplus L(4/5, 3)$ and
\begin{align*}
&A^0 \cong L(\frac{6}7,0)\oplus L(\frac{6}7,5),&& A^1 \cong
L(\frac{6}7,\frac{4}3)^+,
&& A^2\cong L(\frac{6}7,\frac{4}3)^-,\\
&B^0 \cong L(\frac{6}7,\frac{1}7)\oplus L(\frac{6}7,\frac{22}7),&&
B^1 \cong L(\frac{6}7,\frac{10}{21})^+,
&& B^2\cong L(\frac{6}7,\frac{10}{21})^-,\\
&C^0 \cong L(\frac{6}7,\frac{5}7)\oplus L(\frac{6}7,\frac{12}7),&&
C^1 \cong L(\frac{6}7,\frac{1}{21})^+, && C^2\cong
L(\frac{6}7,\frac{1}{21})^-
\end{align*}
are the irreducible modules of $L(6/7,0)\oplus L(6/7,5)$.

Both of $L(4/5,0)\oplus L(4/5,3)$ and $L(6/7,0)\oplus L(6/7,5)$
are rational VOAs and the fusion rules among their irreducible
modules have been determined in \cite{m5} and \cite{ly2}. There
are two $\Z_3$-symmetries given as follows.
\begin{equation*}
\tau_1=
\begin{cases}\ds
1 & \text{on } M^0, W^0,\\
e^{2\pi \sqrt{-1}/3} &\text{on } M^1, W^1,\\
e^{4\pi \sqrt{-1}/3}& \text{on } M^2, W^2,
\end{cases}
\qquad
\tau_2=
\begin{cases}\ds
1 & \text{on } A^0, B^0, C^0,\\
e^{2\pi \sqrt{-1}/3} &\text{on } A^1, B^1, C^1,\\
e^{4\pi \sqrt{-1}/3}& \text{on } A^2,B^2, C^2.
\end{cases}
\end{equation*}

If $V$ is a VOA which contains $U$ as a subalgebra, then both
$\tau_1$ and $\tau_2$ induce automorphisms of $V$. Moreover, as
automorphisms of $V$, $\tau_1=\tau_2$ and we have that
$\tau_{\hat{e}}\tau_{\hat{f}}= \tau_1$ is of order $3$.

\begin{remark}
Recall that $\hat{e}$ and $\hat{f}$ are fixed by the Weyl group
$W(\Phi) = W(A_2) \times W(E_6)$ of the root system $\Phi = A_2
\oplus E_6$ of $L$. Since $U$ is generated by $\hat{e}$ and
$\hat{f}$, $W(\Phi)$ leaves every element of $U$ invariant. Let
$g$ be an element of order $3$ in $W(\Phi)$ which acts
fixed-point-freely on $\sqrt{2}E_8$. Then $g$ induces a
fixed-point-free action on the Leech lattice $\Lambda$ also (cf.
\cite{kly3}). Since every element of $U$ is fixed by $g$, $U$ is
actually contained in the Moonshine VOA $V^\natural$ by the
$\Z_3$-orbifold construction of $V^{\natural}$ given by Dong and
Mason \cite{dm-z3}. Thus, $\tau_{\hat{e}}\tau_{\hat{f}}= \tau_1$
is of class $3A$ (cf. \cite{kly3,m5}).

In \cite{m6}, Miyamoto showed that for any two conformal vectors
$e$ and $f$ of central charge $1/2$ with $\la e,f\ra=13/2^{10}$ in
the Moonshine VOA $V^\natural$, the vertex subalgebra $W$
generated by $e$ and $f$ must contain a subalgebra of the form
$L(4/5,0)\otimes L(6/7,0)$. In fact, Sakuma and
Yamauchi\,\cite{sy} showed that the algebra $W$ must be isomorphic
to $U$. This  gives a proof that $U$ is contained in $V^\natural$.
\end{remark}

\subsection{$4A$ case} In this case, $L=L(3)= A_3\oplus D_5$,
$n_3=|E_8/L|=4$, and the conformal vectors $\tom^1\in
V_{\sqrt{2}A_3}$ and $\tom^2\in V_{\sqrt{2}D_5}$ defined by
\eqref{eq:CV} are both of central charge $1$. Let
$\ep_1,\dots,\ep_8 \in \R^8$ be such that  $\la \ep_i, \ep_j\ra
=2\delta_{ij}$ for any $i,j\in\{1,\ldots,8\}$. Then
\[
\sqrt{2}E_8=\left\{ \sum_{i=1}^8 a_i\ep_i\left | \,{ {\text{ all }
a_i\in \Z \text{ or all } a_i\in \frac{1}2 +\Z}} \atop {\text{
with } \sum_{i=1}^8 a_i \equiv 0 \mod 2 }\right. \right\}.
\]
Let
\[
\mathcal{K}=\left\{ \sum_{i=1}^8 a_i\ep_i\left | \,{ {\text{ all }
a_i\in \Z \text{ or all } a_i\in \frac{1}2 +\Z}}\right. \right\}.
\]
Then $\sqrt{2}E_8\subset \mathcal{K}$ and $|\mathcal{K}/
\sqrt{2}E_8|=2$.

Now we consider two automorphisms $\theta$ and $\psi$ of
$V_\mathcal{K}$ defined by
\begin{gather*}
\theta (\al(-1)\cdot 1)= -\al(-1)\cdot 1\quad \text{ and } \quad
\theta(e^\al) =e^{-\al},\\
\psi(u\otimes e^\beta)= (-1)^{\la \ep_1+\cdots+\ep_8,\beta\ra/2}
u\otimes e^\beta
\end{gather*}
for any $\al, \beta \in \mathcal{K}$ and $u\in M(1)$. Note that
$V_{\sqrt{2}E_8}= (V_{\mathcal{K}})^\psi$. It is well known (cf.
\cite[Chapter 10]{flm}) that $\theta$ and $\psi$ are conjugate in
$\Aut\,V_\mathcal{K}$. Thus, we have
\[
V_{\sqrt{2}E_8}= (V_{\mathcal{K}})^\psi\cong V_{\mathcal{K}}^+.
\]

By using the same argument as in Dong et al.\ \cite{dly,dly2}, one
can show that the VOA $V_{\sqrt{2}A_3}$ contains a subalgebra
isomorphic to $V_{\sqrt{2}A_2}^+\otimes V_{\Z \gamma_3}^+$ and the
VOA $V_{\sqrt{2}D_5}$ contains a subalgebra isomorphic to
$V_{\sqrt{2}A_4}^+\otimes V_{\Z\gamma_5}^+$ where
$\sqrt{2}A_2\cong
\mathrm{span}_{\mathbb{Z}}\{-\ep_1+\ep_2,-\ep_2+\ep_3\}$,
$\gamma_3=\ep_1+\ep_2+\ep_3$,  $ \sqrt{2}A_4\cong
\mathrm{span}_{\mathbb{Z}}\{-\ep_4+\ep_5,-\ep_5+\ep_6,-\ep_6+\ep_7,
-\ep_7+\ep_8\}$ and $\gamma_5=\ep_4+\ep_5+\ep_6+\ep_7+\ep_8$.
Moreover, $U$ contains a subalgebra isomorphic to  $V_{\Z
\gamma_3}^+\otimes V_{\Z\gamma_5}^+$. Define
\[
\mathcal{N}=\{ x\in \mathcal{K}\,|\, \la x, -\ep_j+\ep_{j+1}  \ra
=0 \text{ for } j = 1,2,4,5,6,7 \}.
\]
Then
\[
V_\mathcal{N}= \{ v\in V_{\mathcal{K}}\,|\, u_nv=0 \text{ for any
} u\in V_{\sqrt{2}A_2} \text{ or } V_{\sqrt{2}A_4} \text{ and }
n\geq 0\}.
\]
Since $U=(V_\mathcal{N})^\psi$, it is now easy to see that $U\cong
V^{+}_{\mathcal{N}}$. Note that
\[
\mathcal{N}=\left\{ \sum_{i=1}^8 a_i\ep_i\left | \,{ {\text{ all }
a_i\in \Z \text{ or all } a_i\in \frac{1}2 +\Z} \text{ with }} \atop {
a_1+a_2+a_3= a_4+a_5+a_6+a_7+a_8=0}\right. \right\}.
\]
It is of rank two and generated by the elements
\begin{equation*}
\xi_1 = \frac{1}{2}(\gamma_3+\gamma_5),\qquad
\xi_2 = \frac{1}{2}(-\gamma_3+\gamma_5).
\end{equation*}

\begin{proposition}
The VOA $V_{\mathcal{N}}^+$ is generated by its weight $2$
subspace $(V_{\mathcal{N}}^+)_{_2}$.
\end{proposition}

\begin{proof}
Recall that $V_{\mathcal{N}}^+$ contains a subalgebra
isomorphic to $V_{\Z \gamma_3}^+\otimes V_{\Z\gamma_5}^+$. Moreover,
\begin{equation*}
\begin{split}
V_{\mathcal{N}}^+ =& \left(V_{\Z \gamma_3}^+\otimes
V_{\Z\gamma_5}^+\right) \oplus \left(V_{\Z \gamma_3}^-\otimes
V_{\Z\gamma_5}^-\right)\\
& \oplus \left(V_{\gamma_3/2+\Z \gamma_3}^+\otimes
V_{\gamma_5/2+\Z\gamma_5}^+\right)\oplus \left(V_{\gamma_3/2+\Z
\gamma_3}^-\otimes V_{\gamma_5/2+\Z\gamma_5}^-\right)
\end{split}
\end{equation*}

Note that $(V_{\mathcal{N}}^+)_{_2}$ is of dimension $5$ and it
has a basis consisting of the following five elements.
\[
\be_1(-1)^2\cdot 1,\q \be_2(-1)^2 \cdot 1,\q
\be_1(-1)\be_2(-1)\cdot 1,\q e^{\xi_1}+e^{-\xi_1},\q
e^{\xi_2}+e^{-\xi_2},
\]
where $\be_1=\frac{1}{\sqrt{6}}\gamma_3$ and
$\be_2=\frac{1}{\sqrt{10}}\gamma_5$. Let $W$ be the subalgebra
generated by $(V_{\mathcal{N}}^+)_{_2}$. We want to show that
$W=V_{\mathcal{N}}^+$. Since the minimal weights of $V_{\Z
\gamma_3}^-\otimes V_{\Z\gamma_5}^-$, $V_{\gamma_3/2+\Z
\gamma_3}^+\otimes V_{\gamma_5/2+\Z\gamma_5}^+$, and
$V_{\gamma_3/2+\Z \gamma_3}^-\otimes V_{\gamma_5/2+\Z\gamma_5}^-$
are all equal to $2$ and since they are irreducible $V_{\Z
\gamma_3}^+\otimes V_{\Z\gamma_5}^+$-modules, it suffices to show
that $V_{\Z \gamma_3}^+\otimes V_{\Z\gamma_5}^+\subset W$.

For any rank one even lattice $\Z\al$ with $\la
\al,\al\ra=2k$, it was shown in \cite{dg} that the VOA
$V_{\Z\al}^+$ is generated by three elements
\begin{equation}\label{eq:3G}
\begin{split}
\om=&\frac{1}2 \be(-1)^2\cdot 1, \\
J=&\be(-1)^4\cdot 1 -2\be(-3)\be(-1)\cdot
1+\frac{3}2\be(-2)^2\cdot
1,\\
E=&e^\al+e^{-\al},
\end{split}
\end{equation}
where $\be=\frac{1}{\sqrt{2k}}\al$.

Let $J_1$ and $J_2$ be the elements obtained by replacing $\beta$
with $\beta_1$ and $\beta_2$ in the element $J$ of \eqref{eq:3G},
respectively. Likewise, let $E_1=e^{\gamma_3}+e^{-\gamma_3}$ and
$E_2=e^{\gamma_5}+e^{-\gamma_5}$. We only need to show that $W$
contains $J_1,J_2,E_1$, and $E_2$.

By direct computation, we have
\begin{equation*}
\begin{split}
&(\be_1(-1)\be_2(-1)\cdot 1)_{-1} (\be_1(-1)\be_2(-1)\cdot
1)\\
& \qquad =\be_1(-3)\be_1(-1)\cdot 1 +\be_2(-3)\be_2(-1)\cdot 1 +
\be_1(-1)^2\be_2(-1)^2\cdot 1
\end{split}
\end{equation*}
and
\[
(\be_1(-1)^2\cdot 1)_{-1} (\be_2(-1)^2\cdot 1)
=\be_1(-1)^2\be_2(-1)^2\cdot 1.
\]
Thus, $\be_1(-3)\be_1(-1)\cdot 1 +\be_2(-3)\be_2(-1)\cdot 1 \in W$
and
\[
\be_1(-3)\be_1(-1)\cdot 1= \frac{1}{8}(\be_1(-1)^2\cdot 1)_1 (
\be_1(-3)\be_1(-1)\cdot 1 +\be_2(-3)\be_2(-1)\cdot 1)
\]
is also contained in $W$. Moreover,
\begin{gather*}
(\be_1(-1)^2\cdot 1)_{-1} (\be_1(-1)^2\cdot 1) =\be_1(-1)^4 \cdot
1 + 4\be_1(-3)\be_1(-1)\cdot 1,\\
((\be_1(-1)^2\cdot 1)_{0})^2 (\be_1(-1)^2\cdot 1) = 16
\be_1(-3)\be_1(-1)\cdot 1 + 8 \be_1(-2)^2 \cdot 1.
\end{gather*}
Hence, we have $\be_1(-1)^4 \cdot 1, \be_1(-2)^2 \cdot 1\in W$ and
thus  $J_1\in W$. Similarly, $J_2\in W$.

Clearly, $W$ contains
\[
E_1=e^{\gamma_3}+e^{-\gamma_3}= ( e^{\xi_1}+e^{-\xi_1})_0
( e^{\xi_2}+e^{-\xi_2}).
\]
Furthermore, $W$ contains the following three elements.
\begin{align*}
(\gamma_3(-1)^2\cdot 1)_{-1}E_1 &= \gamma_3(-1)^2E_1 +
12\gamma_3(-2)(e^{\gamma_3}-e^{-\gamma_3}),\\
(\gamma_5(-1)^2\cdot 1)_{-1}E_1 &=
\gamma_5(-1)^2 E_1,\\
(\gamma_3(-1)\gamma_5(-1)\cdot 1)_{-1}E_1 &=
\gamma_3(-1)\gamma_5(-1)E_1 +
6\gamma_5(-2)(e^{\gamma_3}-e^{-\gamma_3}).
\end{align*}
Then $W$ also contains
\begin{equation*}
\begin{split}
&(\gamma_3(-1)\gamma_5(-1)\cdot 1)_{1}\big(\gamma_3(-1)\gamma_5(-1)
E_1 +
6\gamma_5(-2)(e^{\gamma_3}-e^{-\gamma_3})\big)\\
& \qquad = 10\left(\gamma_3(-1)^2E_1 +12\gamma_3(-2)
(e^{\gamma_3}-e^{-\gamma_3})\right)+6\gamma_5(-1)^2E_1
+6\gamma_5(-2) (e^{\gamma_3}-e^{-\gamma_3}).
\end{split}
\end{equation*}
Hence, we have $\gamma_5(-2) (e^{\gamma_3}-e^{-\gamma_3})\in W$
and so
\[
\gamma_3(-2)
(e^{\gamma_3}-e^{-\gamma_3})=\frac{1}{20}(\gamma_3(-1)\gamma_5(-1)\cdot
1)_{1}(\gamma_5(-2)(e^{\gamma_3}-e^{-\gamma_3}))
\]
is also contained in $W$. Therefore, $\gamma_3(-1)^2E_1$,
$\gamma_5(-1)^2E_1$, and $\gamma_5(-1)\gamma_3(-1)E_1$ are
contained in $W$.

Finally, $W$ contains
\begin{equation*}
\begin{split}
&( e^{\xi_1}+e^{-\xi_1})_{-2} (
e^{\xi_2}+e^{-\xi_2})\\
& \qquad = E_2 + \frac{1}2\xi_1(-1)^2 E_1 +
\frac{1}2
\xi_1(-2)(e^{\gamma_3}-e^{-\gamma_3})\\
& \qquad = E_2 +
\frac{1}8(\gamma_3+\gamma_5)(-1)^2 E_1 + \frac{1}4
(\gamma_3+\gamma_5)(-2)(e^{\gamma_3}-e^{-\gamma_3})
\end{split}
\]
and thus $E_2 \in W$.
\end{proof}

\begin{theorem}
$(1)$\ The Griess algebra $U_2$ of $U$ is generated by $\hat{e}$
and $\hat{f}$.

$(2)$\ The coset subalgebra  $U$ is generated by $\hat{e}$ and
$\hat{f}$.
\end{theorem}

\begin{proof}
We only need to show the first assertion.  Note that
$U_2=\mathrm{span}_\C\{ \tom^1,\tom^2,X^1,X^2,X^3\}$ is of
dimension $5$ and that
\begin{align*}
\hat{e} &= \tom^1+\tom^2 +X^1+X^2+X^3,\\
\hat{f} &= \tom^1+\tom^2 +\sqrt{-1}X^1-X^2-\sqrt{-1}X^3.
\end{align*}

Let $\mathcal{G}$ be the Griess subalgebra generated by $\hat{e}$
and $\hat{f}$. Then $\hat{e}_1 \hat{f}$, $\hat{e}_1(\hat{e}_1
\hat{f})$, and $\hat{f}_1(\hat{e}_1 \hat{f})$ are also in
$\mathcal{G}$. By direct computation, it is easy to see that
$\hat{e}$, $\hat{f}$, $\hat{e}_1 \hat{f}$, $\hat{e}_1(\hat{e}_1
\hat{f})$, and $\hat{f}_1(\hat{e}_1 \hat{f})$ are linearly
independent. Thus $\mathcal{ G}=U_2$.
\end{proof}

\begin{theorem}
The automorphism group $\mathrm{Aut}\, U $ of $U$ is a dihedral
group of order $8$.
\end{theorem}

\begin{proof}
There are exactly four conformal vectors of
central $1/2$ in $U_2$, namely,
\begin{equation*}
e_j=\sigma^j \hat{e} = \frac{3}{16}\tom^1 + \frac{5}{16}\tom^2 +
\frac{1}{32} \big((\sqrt{-1})^j X^1 + (-1)^j X^2 + (-\sqrt{-1})^j
X^3\big), \quad 0 \le j \le 3.
\end{equation*}

Since $U$ is generated by $\hat{e}$ and $\hat{f}$, we can consider
$ \mathrm{Aut}\, U$ as a subgroup of the permutation group on the
set $\{{e}_0, {e}_1,{e}_2,{e}_3\}$. Now let $g\in \mathrm{Aut}\,
U$. Then $g$ also preserves the inner product and thus
\[
\la g {e}_i, g {e}_j\ra =\la {e}_i, {e}_j\ra=
\begin{cases}
\ds
1/2^7 & \text{ if } i+j \text{ is odd},\\
0 & \text{ if } i+j \text{ is even}.
\end{cases}
\]
Hence, $g$ will either keep $\{e_0, e_2\}$ and $\{e_1,e_3\}$
invariant or $g$ will map $\{e_0, e_2\}$ to  $\{e_1,e_3\}$. Thus
$|\mathrm{Aut}\,U| \le 8$. On the other hand, $\theta$ and
$\sigma$ defined in Remark \ref{sigma} generate a subgroup
isomorphic to a dihedral group of order $8$ inside $\mathrm{Aut}\,
U$. Hence the assertion holds.
\end{proof}

Note that $\la \xi_1,\xi_1\ra =\la \xi_2,\xi_2\ra=4$, $ \la
\xi_1,\xi_2\ra =-1$ and the Leech lattice $\Lambda$ does contain
some sublattice isomorphic to $\mathcal{N}$. Therefore, $U \subset
V_\Lambda^+ \subset V^\natural$.

In \cite{abe}, the fusion rules for $V_{\Z\gamma_3}^+$ and
$V_{\Z\gamma_5}^+$ are determined. It is known that there are
$\Z_4$-symmetries among the irreducible modules of
$V_{\Z\gamma_3}^+$ and also among the irreducible modules of
$V_{\Z\gamma_5}^+$. By direct computation, it is easy to verify
that the automorphism $\tau_{\hat{e}}\tau_{\hat{f}}$ agrees with
the $\Z_4$-symmetries of $V_{\Z\gamma_3}^+$ and $V_{\Z\gamma_5}^+$
and thus $\tau_{\hat{e}}\tau_{\hat{f}}$ is of class $4A$ (cf.
\cite{mat}). An explicit construction of $4A$-elements as
automorphisms of $V^{\natural}$ has already been obtained by
Shimakura \cite{Shi}.

\subsection{$5A$ case.}\label{sub5-5A}

In this case, $L=L(4)\cong A_4\oplus A_4$, $n_4=|E_8/ L|=5$, and
the conformal vectors $\tom^1 \in V_{\sqrt{2}A_4}$ and $\tom^2 \in
V_{\sqrt{2}A_4}$ defined by \eqref{eq:CV} are both of central
charge $8/7$.

\begin{lemma}[Lemma A.7]
The coset subalgebra $U$ contains a set of three mutually orthogonal conformal
vectors of central charge $1/2, 25/28$, and $25/28$, respectively, namely,
\begin{align*}
& u=\hat{e}=\frac{7}{32} (\tom^1+\tom^2) + \frac{1}{32}( X^1+X^2+X^3+X^4),\\
& v= \frac{15}{64}\tom^1+\frac{35}{64}\tom^2 -
\frac{3}{64}(X^1+X^4) +\frac{1}{64}(X^2+X^3),\\
& w= \frac{35}{64}\tom^1+\frac{15}{64}\tom^2 +
\frac{1}{64}(X^1+X^4) -\frac{3}{64}(X^2+X^3).
\end{align*}
\end{lemma}

By the above lemma, $U$ contains $\Vir(u)\otimes \Vir(v) \otimes
\Vir(w) \cong L(1/2,0)\otimes L(25/28,0)\otimes
L(25/28,0)$. All irreducible modules of $\Vir(u)\otimes
\Vir(v) \otimes \Vir(w)$ are known (cf.
\cite{ff,W}). They are of the form
$L(1/2,h_1)\otimes L(25/28,h_2)\otimes
L(25/28,h_3)$.
Among them, the highest weights $(h_1,h_2,h_3)$ of the irreducible
modules which have
integral weights are as follows.
\begin{align*}
&(0,0,0), \  (\frac{1}{16},\frac{5}{32},\frac{57}{32}), \
       (\frac{1}{16}, \frac{57}{32}, \frac{5}{32}),\ (\frac{1}2, \frac{3}4, \frac{3}4) ,\
       (0, \frac{3}4, \frac{13}4) ,\
   (0, \frac{13}4, \frac{3}4),\\
   & (\frac{1}{16}, \frac{57}{32}, \frac{165}{32}),\
   (\frac{1}{16},\frac{165}{32},\frac{57}{32}), \
   (\frac{1}2, \frac{13}4, \frac{13}4), \
   (\frac{1}2, 0, \frac{15}2) ,\ (\frac{1}2, \frac{15}2, 0) ,\ (0,\frac{15}2,\frac{15}2).
\end{align*}

The following lemma can be proved by direct computation.

\begin{lemma}
  Let
  \begin{align*}
    a &= (\tom^1-\tom^2)- \frac{1}{35} ( X^1+X^4) + \frac{1}{35}( X^2+X^3),
    \\
    b^1 &= 3(X^1-X^4)+ 2(X^2-X^3),
    \\
    b^2 &= 2(X^1-X^4) -3(X^2-X^3).
  \end{align*}
  Then $a$, $b^1$ and $b^2$ are highest weight vectors of highest
  weight $(\fr{1}{2},\fr{3}{4},\fr{3}{4})$,
  $(\fr{1}{16},\fr{5}{32},\fr{57}{32})$ and
  $(\fr{1}{16},\fr{57}{32},\fr{5}{32})$ with respect to
  $\Vir(u)\otimes \Vir(v)\otimes \Vir(w)$, respectively.
\end{lemma}

\medskip
We denote the irreducible module $L(1/2,h_1)\otimes L(25/28,h_2)\otimes
L(25/28,h_3)$ by $[h_1,h_2,h_3]$ for simplicity of notation.
By using the theory of characters (cf. Appendix B), we actually
have the following decomposition of $U$ into a direct sum of
$[h_1,h_2,h_3]$'s.

\begin{theorem}[Theorem \ref{B-5A}]
As a module of $\Vir(u)\otimes \Vir(v)\otimes \Vir(w)$,
\begin{align*}
U\cong & [0,0,0] \oplus [\frac{1}{16},\frac{5}{32},\frac{57}{32}]
 \oplus [\frac{1}{16},\frac{57}{32},\frac{5}{32}]
 \oplus [\frac{1}{2},\frac{3}{4},\frac{3}{4}]\\
& \oplus [0,\frac{3}{4},\frac{13}{4}]
 \oplus [0,\frac{13}{4},\frac{3}{4}]
 \oplus [\frac{1}{16},\frac{57}{32},\frac{165}{32}]
 \oplus [\frac{1}{16},\frac{165}{32},\frac{57}{32}]\\
 &\oplus [\frac{1}{2},\frac{13}{4},\frac{13}{4}]
\oplus [\frac{1}{2},0,\frac{15}{2}]
\oplus [\frac{1}{2},\frac{15}{2},0] \oplus
[0,\frac{15}{2},\frac{15}{2}].
\end{align*}
\end{theorem}

Next, we shall discuss the generators of $U$.

\begin{theorem}\label{5a-gen}
The coset subalgebra  $U$ is generated by its weight $2$ subspace
$U_2$.
\end{theorem}

We shall divide the proof into several steps. By direct
computation, we can verify the following lemma.

\begin{lemma}
  Let
  \begin{align*}
    y^1 &= a_{-1}a-\frac{192}{343}u_{-3}\cdot 1 -\frac{128}{343}u_{-1}u
      -\frac{192}{455} v_{-3}\cdot 1 -\frac{128}{325}v_{-1}v
      -\frac{192}{455} w_{-3}\cdot 1 -\frac{128}{325}w_{-1}w,
   \\
   y^2 &= (b^1)_{-1}b^2 +105u_{-1}a -\frac{735}{2} v_{-1}a +\frac{5355}{8}v_0v_0a
    -\frac{735}{2} w_{-1}a +\frac{1715}{8} w_0 w_0 a
  \\
   & \quad +\frac{665}{2}u_0v_0a -\frac{245}{2} u_0w_0a
     -\frac{4655}{12}v_0w_0a.
  \end{align*}
  Then $y^1$ and $y^2$ are non-zero singular vectors for
  $\Vir(u)\otimes \Vir(v)\otimes \Vir(w)$, that is, $u_n y^j=v_n y^j
  = w_n y^j = 0$, $j=1,2$, for any $n \geq 2$.
\end{lemma}

Now let $W$ be the subalgebra of $U$ generated by $U_2$.

\begin{lemma}\label{lem:3.16}
There are highest weight vectors of highest weight $(0,\fr{3}{4},\fr{13}{4})$
and $(0,\fr{13}{4},\fr{3}{4})$ with respect to $\Vir(u)\otimes \Vir(v)\otimes
\Vir(w)$ in $W$.
\end{lemma}

\begin{proof}
By fusion rules, we know that
\begin{gather*}
a_{-1} a\in [0,0,0]_4\oplus [0,\fr{3}{4},\fr{13}{4}]_4\oplus
[0,\fr{13}{4},\fr{3}{4}]_4,\\
(b^1)_{-1}b^2\in [\fr{1}{2},\fr{3}{4},
\fr{3}{4}]_4\oplus [0,\fr{3}{4},\fr{13}{4}]_4\oplus
[0,\fr{13}{4},\fr{3}{4}]_4.
\end{gather*}

Since both $y^1$ and $y^2$ are singular vectors of weight $4$, we
have
\[
y^1,y^2 \in [0,\fr{3}{4},\fr{13}{4}]\oplus
[0,\fr{13}{4},\fr{3}{4}].
\]
Moreover, by direct computation, one can  show that $ \la y^1,
y^2\ra=\la a_{-1}a, (b^1)_{-1}b^2\ra = 0$. Note that $a,b^1$, and
$b^2$ are highest weight vectors and $\la a, b^1\ra =\la a, b^2\ra
=\la b^1,b^2\ra=0$. Hence $y^1$ and $y^2$ are linearly
independent. Thus the assertion holds.
\end{proof}

\begin{proof}[\textbf{Proof of Theorem \ref{5a-gen}}]
\

First, we note that the coset subalgebra $U$ is simple. Then the
subalgebra $U^{\tau_u}$ consisting of the fixed points of the
Miyamoto involution $\tau_u$ associated with $u$ in $U$ is also
simple. As a module of  $\Vir(u)\otimes \Vir(v)\otimes \Vir(w)$,
\begin{align*}
U^{\tau_u} \cong& [0,0,0]\oplus [0,\fr{3}{4},
\fr{13}{4}]\oplus [0,\fr{13}{4},\fr{3}{4}]\oplus
[0,\fr{15}{2},\fr{15}{2}]\oplus
[\fr{1}{2},\fr{3}{4},\fr{3}{4}]\\
&\oplus  [\fr{1}{2},\fr{13}{4},\fr{13}{4}] \oplus
[\fr{1}{2},0,\fr{15}{2}] \oplus [\fr{1}{2},\fr{15}{2},0].
\end{align*}
Moreover, we can define an automorphism $\sigma_u$ on $U^{\tau_u}$
(cf. Miyamoto \cite{m1}) by
\begin{equation*}
\sigma_u=
\begin{cases}
1 &\text{ on } L(\frac{1}{2},0),\\
-1 &\text{ on } L(\frac{1}{2},\frac{1}{2}).
\end{cases}
\end{equation*}
Then the subalgebra $(U^{\tau_u})^{\sigma_u}$ consisting of the
fixed points of $\sigma_u$ in $U^{\tau_u}$ is again simple.
As a module of  $\Vir(u)\otimes \Vir(v)\otimes \Vir(w)$,
\[
(U^{\tau_u})^{ \sigma_u} \cong [0,0,0]\oplus [0,\fr{3}{4},
\fr{13}{4}]\oplus [0,\fr{13}{4},\fr{3}{4}]\oplus
[0,\fr{15}{2},\fr{15}{2}].
\]

By the previous lemma, we know that $W$ contains $[0,0,0]$,
$[0,3/4,13/4]$, $[0,13/4,3/4]$, $[1/2,3/4,3/4]$,
$[1/16,5/32,57/32]$ and $[1/16,57/32,5/32]$. Hence,
$(W^{\tau_u})^{\sigma_u}$ must contain an irreducible module
isomorphic to $[0,15/2,15/2]$; otherwise $(W^{ \tau_u
})^{\sigma_u}\cong [0,0,0]\oplus [0,3/4,13/4]\oplus [0,13/4,3/4]$
and the orthogonal complement of $(W^{\tau_u })^{\sigma_u}$ in
$(U^{\tau_u })^{\sigma_u}$ with respect to a positive definite
invariant hermitian form is isomorphic to $[0,15/2,15/2]$.
However, the orthogonal complement is a module for $(W^{
\tau_u})^{\sigma_u}$ because of the invariance of the form, which
is impossible by the fusion rules. Thus $(W^{ \tau_u})^{
\sigma_u}=(U^{\tau_u})^{ \sigma_u}$. Then we also have $W^{
\tau_u}=U^{\tau_u}$ since $U^{\tau_u}$ is a direct sum of two
irreducible $(W^{ \tau_u})^{ \sigma_u}$-modules and $W^{
\tau_u}\neq (W^{ \tau_u})^{ \sigma_u}$. Hence, $W$ contains all
the simple current $\vir (u)\tensor \vir (v)\tensor \vir
(w)$-modules. Now it is easy to see that $W=U$ by the fusion
rules.
\end{proof}

\begin{theorem}
$(1)$\ The Griess algebra $U_2$ of $U$ is generated by $\hat{e}$
and $\hat{f}$.

$(2)$\  The coset subalgebra $U$ is generated by $\hat{e}$ and
$\hat{f}$.
\end{theorem}

\begin{proof}
By Theorem \ref{5a-gen}, it suffices to show the first assertion.
Let $\mathcal{G}$ be the Griess subalgebra generated by $\hat{e}$
and $\hat{f}$. Then by direct computation, we can verify that
$\hat{e}$,
  $\hat{f}$,
  $\hat{e}_1\hat{f}$,
  $\hat{e}_1(\hat{e}_1 \hat{f})$,
  $\hat{f}_1(\hat{e}_1(\hat{e}_1 \hat{f}))$,
  $\hat{e}_1(\hat{f}_1(\hat{e}_1(\hat{e}_1 \hat{f})))$
are linearly independent. Thus $\mathcal{ G}=U_2$, since $\dim U_2
= 6$.
\end{proof}

\begin{theorem}
The automorphism group $\mathrm{Aut}\, U $ of $U$ is a dihedral
group of order $10$.
\end{theorem}

\begin{proof}
Recall that $\sigma$ and $\theta$ defined in Remark \ref{sigma}
generate a subgroup isomorphic to a dihedral group of order $10$
in $\mathrm{Aut}\,U$. By Lemma \ref{A-5A}, there are exactly five
conformal vectors of central charge $1/2$ in $U$, namely, $e_j =
\sigma^j \hat{e}, 0 \le j \le 4$. Since $U$ is generated by
$\hat{e}$ and $\hat{f}=\sigma \hat{e}$, $\mathrm{Aut}\,U$ can be
considered as a subgroup of the permutation group on the five
elements set $\{e_0, e_1,e_2,e_3, e_4\}$. In fact, $U$ is
generated by any $2$ distinct elements in the set $\{e_0,
e_1,e_2,e_3, e_4\}$. Thus $\mathrm{Aut}\,U$ can not contain any
$2$-cycle nor $3$-cycle since such automorphisms must fix at least
two elements in the set and thus fix the whole $U$. Hence there is
no element of order $6$ either. Now let $(i_1\,i_2\,i_3\,i_4)$ be a
$4$-cycle. Then
\[
(i_0\,i_1\,i_2\,i_3\,i_4)(i_1\,i_2\,i_3\,i_4)= (i_0\,i_1\,i_3)(i_2\,i_4).
\]
Since $(i_0\,i_1\,i_2\,i_3\,i_4)\in \mathrm{Aut}\,U$ but
$(i_0\,i_1\,i_3)(i_2\, i_4) \notin \mathrm{Aut}\,U$, there is no
$4$-cycle in $\mathrm{Aut}\,U$. Therefore, the only possible
elements are $5$-cycles and the products of two disjoint
$2$-cycles. Thus the assertion holds.
\end{proof}

\begin{theorem}
There are exactly nine irreducible modules $U(i,j),i,j=1,3,5$ for
$U$. As $\Vir(u) \otimes \Vir(v)\otimes \Vir(w)$-modules, they are
of the following form.
\begin{equation*}
\begin{split}
U(i,j) & \cong [0, h_{i,1}, h_{j,1}]\oplus [0, h_{i,3},
h_{j,5}]\oplus
[0, h_{i,5}, h_{j,3}]\oplus [0, h_{i,7}, h_{j,7}]\\
& \quad \oplus [\frac{1}2, h_{i,1}, h_{j,7}]\oplus [\frac{1}2,
h_{i,3}, h_{j,3}]\oplus
[\frac{1}2, h_{i,5}, h_{j,5}]\oplus [\frac{1}2, h_{i,7}, h_{j,1}]\\
& \quad \oplus [\frac{1}{16}, h_{i,2}, h_{j,4}]\oplus
[\frac{1}{16}, h_{i,4}, h_{j,2}]\oplus [\frac{1}{16}, h_{i,6},
h_{j,4}]\oplus [\frac{1}{16}, h_{i,4}, h_{j,6}],
\end{split}
\end{equation*}
where $h_{r,s}=h_{r,s}^5$ is defined by \eqref{eq:VIRH}.
\end{theorem}

The proof of this theorem will be given at the Appendix C. We
shall first note that $U= U(1,1)$ and the lattice VOA
$V_{\sqrt{2}E_8}$ can be decomposed as follows (cf. Appendix B.2).
\[
V_{\sqrt{2}E_8}\cong \bigoplus_{ { 1\leq k_j, \ell_j\leq j+2}
\atop{{\ k_j, \ell_j \equiv 1\,\mathrm{mod}\, 2} \atop
{j=0,1,\ldots, 4}}} L(c_1,h^1_{k_{0},k_1})\otimes
L(c_1,h^1_{\ell_{0},\ell_1})\otimes \cdots\otimes
L(c_4,h^4_{k_{3},k_4})\otimes L(c_4,h^4_{\ell_{3},\ell_4}) \otimes
U(k_4,\ell_4),
\]
where $c_m$ and $h_{r,s}^m$ are defined by \eqref{eq:VIRC} and
\eqref{eq:VIRH}, respectively

Recall that $\tau_{\hat{e}}\tau_{\hat{f}}=e^{2\pi
\sqrt{-1}\bbe(0)}$ as automorphisms of $V_{\sqrt{2}E_8}$ (cf.
Eq.\ \ref{tau-ef}). It thus induces a natural action on each of
the $U(i,j)$, $i,j=1,3,5$. Hence if a VOA $V$ contains a
subalgebra isomorphic to $U$, then $\tau_{\hat{e}}\tau_{\hat{f}}$
will define an automorphism of order $5$ on $V$. The subalgebra
$U^{\tau_{\hat{e}}\tau_{\hat{f}}}$ consisting of the fixed points
of $\tau_{\hat{e}}\tau_{\hat{f}}$ in $U$ is of the form
\begin{equation*}
U^{\tau_{\hat{e}}\tau_{\hat{f}}}\cong W_5(8/7)\otimes W_5(8/7),
\end{equation*}
where $W_5(8/7)$ is a parafermion algebra of central charge $8/7$
(cf. \cite{ly3}). It is well known that $W_5(8/7)$ possesses a
$\Z_5$ symmetry (cf. \cite{dl,zf1}). The automorphism
$\tau_{\hat{e}}\tau_{\hat{f}}$ in fact agrees with this symmetry.

\begin{remark}
Recall that $\hat{e}$ and $\hat{f}$ are fixed by the Weyl group
$W(\Phi) = W(A_4) \times W(A_4)$ of the root system
$\Phi=A_4\oplus A_4$ of $L$. Since $U$ is generated by $\hat{e}$
and $\hat{f}$, $W(\Phi)$ acts trivially on $U$. There is an
element $\psi$ of order $5$ in $W(\Phi)$ such that it induces a
fixed-point-free action on $\sqrt{2}E_8$ and on the Leech lattice
$\Lambda$. Therefore, if the conjectured $\Z_5$-orbifold
construction of the Moonshine VOA $V^\natural$ holds, then one can
prove that $U$ is contained in $V^\natural$ by using $\psi$ and
that as an automorphism of $V^\natural$,
$\tau_{\hat{e}}\tau_{\hat{f}}$ is of class $5A$ (cf. \cite{mat}).
\end{remark}

\subsection{$6A$ case.}
In this case,  $L=L(5)\cong A_2\oplus A_1\oplus A_5 $, $n_5=|E_8/L|=6$,
and the conformal vectors $\tom^1\in V_{\sqrt{2}A_2}$, $\tom^2\in V_{\sqrt{2}A_1}$,
and $\tom^3\in V_{\sqrt{2}A_5}$ defined by \eqref{eq:CV} are of
central charge $4/5$, $1/2$, and $5/4$, respectively.
Let $K = \Span_{\Z} \{\al_0,\al_1,\al_2,\al_3,\al_4,\al_8\} \cong A_1\oplus
A_5$, $J = \Span_{\Z} \{\al_6,\al_7\} \cong A_2$,
and $\tilde{\al}= \al_3 + 2\al_4 + 3\al_5 + 2\al_6 + \al_7 + \al_8$.
Then $\Z\tilde{\alpha}+K = K\cup (\tilde{\al}+K)$.
We have $\la \al_j,\tilde{\alpha}\ra = 0$ for $j=0,1,3,4,6,7$,
$\la \al_2,\tilde{\alpha}\ra = \la \al_8,\tilde{\alpha}\ra = -1$,
$\la \tilde{\alpha},\tilde{\alpha}\ra = 2$, and
$\alpha_0 + \alpha_4 + \alpha_8 + 2\alpha_1 + 2\alpha_3
+ 2\tilde{\alpha} + 3\alpha_2 = 0$.
Hence $\{\alpha_0, \alpha_1,\alpha_2,\alpha_3,\alpha_4,\tilde{\alpha},\alpha_8\}$
forms an extended $E_6$ diagram and so $K\cup (\tilde{\al}+K) \cong E_6$.
Moreover, $L=J\oplus K$ and $\la J,\Z\tilde{\alpha}+K\ra = 0$.
Therefore, we have isometric embeddings
\[
  A_2\oplus A_1\oplus A_5 \subset A_2\oplus E_6 \subset E_8.
\]
Then we obtain a conformal vector
\[
  \tom(E_6)= \frac{1}{168} \sum_{\al\in \Phi^+(E_6)} \al(-1)^2\cdot
  1 + \frac{1}{14} \sum_{\al\in \Phi^+(E_6)} ( e^{\sqrt{2}\al}+e^{-\sqrt{2}\al})
\]
of central charge $6/7$ in $V_{\sqrt{2}E_8}$.
By our embeddings, we also know that
\[
  V_{\sqrt{2}E_6}\cong \{ v\in V_{\sqrt{2}E_8}|\ (s^1)_1v=(\tom^1)_1
  v=0\}.
\]

If $\alpha \in 3\alpha_5 + L$ satisfies $\la\alpha,\alpha\ra = 2$,
then $\alpha \in \tilde{\al}+K$ since
$3\alpha_5 + L = \tilde{\al} + L$ and $\la J, \tilde{\al}+K\ra = 0$.
Thus $\Phi(E_6) = \Phi(K) \cup \{ \alpha \in 3\alpha_5+L\,|\,\la\alpha,\alpha\ra=2\}$.
Note also that the Virasoro element of $V_{\sqrt{2}E_6}$ coincides with
that of $V_{\sqrt{2}K}$. Now we can verify that
\begin{equation*}
  \tilde{\omega}(E_6) = \frac{2}{7}\tilde{\omega}^2 +
  \frac{4}{7}\tilde{\omega}^3 + \frac{1}{14}X^3 \in U.
\end{equation*}

Let $w^1=\tom^1$, $w^2=\tom(E_6)$, and $w^3=\tom^2+\tom^3-w^2$.
Then $\{w^1,w^2,w^3\}$ is a set of mutually orthogonal conformal
vectors of central charge $4/5, 6/7$, and $25/28$, respectively
and the Virasoro element
$\omega' = \tilde{\omega}^1 + \tilde{\omega}^2 + \tilde{\omega}^3$
of $U$ is a sum of $w^1,w^2,w^3$.
Note that
\begin{equation*}
  s(E_6)=\frac{1}{28} \sum_{\al\in \Phi^+(E_6)} \Big(
  \al(-1)^2\cdot 1 - 2 ( e^{\sqrt{2}\al}+e^{ -\sqrt{2}\al}) \Big)
\end{equation*}
is a linear combination of $s^2$, $s^3$, and $w^3$.
Recall that $s^1 = s(A_2)$, $s^2 = s(A_1)$, and $s^3 = s(A_5)$ are defined
by \eqref{eq:CV}.

Now let $ U'=\{ v\in V_{\sqrt{2}E_8}\, |\, (s^1)_1v=s(E_6)_1 v=0\}$.
Then
\begin{equation*}
\begin{split}
  U'
  &= \{ v\in V_{\sqrt{2}E_8}\, |\, (s^1)_1v=(s^2)_1v=(s^3)_1v =(w^3)_1 v=0\}
  \\
  &= \{ u\in U \,|\, (w^3)_1 u=0\}.
\end{split}
\end{equation*}
Hence $U'\subset U$ and by using the results in the $3A$ case, we
know that
\begin{equation*}
\begin{split}
  U' \cong
  & \Big( L(\frac{4}{5}, 0)\oplus L(\frac{4}{5}, 3)\Big)
    \otimes \Big(L(\frac{6}{7},0)\oplus L(\frac{6}{7}, 5)\Big)
  \\
  & \oplus L(\frac{4}{5},\frac{2}{3})\otimes L(\frac{6}{7}, \frac{4}{3})
    \oplus L(\frac{4}{5}, \frac{2}{3})\otimes L(\frac{6}{7}, \frac{4}{3}).
\end{split}
\end{equation*}

Note that $U'\otimes L(25/28,0)$, $U'(5/7)\otimes L(25/28,9/7)$,
and $U'(1/7)\otimes L(25/28,34/7)$ are the only irreducible
$U'\otimes L(25/28,0)$-modules which have integral weights (cf.\ \cite{sy}),
where
$$
\begin{array}{ll}
  U'(\dfrac{1}7) \cong
  & \Big(L(\dfrac{4}{5},0)\oplus L(\dfrac{4}{5},3)\Big) \otimes
    \Big(L(\dfrac{6}{7}, \dfrac{1}7)\oplus L(\dfrac{6}{7},\dfrac{22}7)\Big)
  \vsb\\
  & \oplus L(\dfrac{4}{5}, \dfrac{2}{3})\otimes L(\dfrac{6}{7}, \dfrac{10}{21})
    \oplus L(\dfrac{4}{5}, \dfrac{2}{3})\otimes L(\dfrac{6}{7}, \dfrac{10}{21}),
  \vsb\\
  U'(\dfrac{5}7) \cong
  & \Big(L(\dfrac{4}{5},0)\oplus L(\dfrac{4}{5},3)\Big) \otimes
    \Big(L(\dfrac{6}{7}, \dfrac{5}7)\oplus L(\dfrac{6}{7},\dfrac{12}7)\Big)
  \vsb\\
  & \oplus L(\dfrac{4}{5}, \dfrac{2}{3})\otimes L(\dfrac{6}{7}, \dfrac{1}{21})
    \oplus L(\dfrac{4}{5}, \dfrac{2}{3})\otimes L(\dfrac{6}{7}, \dfrac{1}{21}).
\end{array}
$$

By direct computation, it is straightforward to verify that $U$ is
generated by its weight $2$ subspace and it contains some highest
weight vectors of weight $5/7$ and $1/7$ with respect to $U'$
(cf. see Appendix B.3.1 for details).
Thus we have the following theorem.

\begin{theorem}
As a module of $U'\otimes L(25/28,0)$,
\[
U\cong U'\otimes L(\frac{25}{28},0) \oplus U'(\frac{5}7)\otimes
L(\frac{25}{28},\frac{9}7)\oplus U'(\frac{1}7)\otimes
L(\frac{25}{28},\frac{34}7).
\]
Moreover, $U$ is generated by its weight $2$ subspace $U_2$.
\end{theorem}

\begin{theorem}
$(1)$\ The Griess algebra $U_2$ of $U$ is generated by $\hat{e}$
and $\hat{f}$.

$(2)$\ The coset subalgebra $U$ is generated by $\hat{e}$ and
$\hat{f}$.
\end{theorem}

\begin{proof}
Let $\mathcal{G}$ be the Griess subalgebra generated by $\hat{e}$
and $\hat{f}$. By \eqref{6a-eq}, it is easy to verify that
  $\hat{e}$,
  $\hat{f}$,
  $\hat{e}_1 \hat{f}$,
  $\hat{e}_1 (\hat{e}_1 \hat{f})$,
  $\hat{f}_1 (\hat{e}_1 \hat{f})$,
  $\hat{f}_1 (\hat{e}_1 (\hat{e}_1 \hat{f}))$,
  $\hat{e}_1 (\hat{f}_1 (\hat{e}_1 \hat{f}))$,
  $\hat{e}_1 (\hat{f}_1 (\hat{e}_1 (\hat{e}_1 \hat{f})))$
are linearly independent. Thus $\mathcal{ G}=U_2$, since $\dim U_2
= 8$. The second assertion follows from the preceding theorem.
\end{proof}

\begin{theorem}
The automorphism group $\mathrm{Aut}\, U $ of $U$ is a dihedral
group of order $12$.
\end{theorem}

\begin{proof}
By Lemma \ref{A-6A}, there are exactly seven conformal vectors of
central charge $1/2$ in $U$, namely, $\tom^2$ and $e_j=\sigma^j
\hat{e}, 0 \le j \le 5$. Moreover, we have $\la \tom^2, e_j\ra
=1/32$ for any $0 \le j \le 5$ and
\begin{equation}\label{inn-6a}
\la e_i,e_j\ra=
\begin{cases}
    1/4 &\text{ if } i=j,\\
    1/32 &\text{ if } i-j\equiv 3 \mod 6,\\
    13/2^{10} &\text{ if } i-j\equiv \pm 2 \mod 6,\\
    5/2^{10} &\text{ if } i-j\equiv \pm 1 \mod 6.
\end{cases}
\end{equation}

Since $U$ is generated by $\hat{e}$ and $\hat{f}=\sigma\hat{e}$,
$\Aut\,U$ can be considered as a subgroup of the permutation group
on the set $\{\tom^2, e_0, e_1,\ldots,e_5\}$. Now let $g\in \Aut\,
U$. Then $g$ must preserve the inner product and so $g$ fixes
$\tom^2$. Let $g e_j = e_{\nu (j)}, 0\le j \le 5$. Then  by
\eqref{inn-6a}, $\nu(i)-\nu(j)\equiv \pm (i-j) \mod 6 $ for any
$i,j=0,1,\dots, 5$. Hence there are only $12$ possible choices for
$\nu$. Thus the assertion holds since $\sigma$ and $\theta$ in
Remark \ref{sigma} generate a subgroup of $\Aut\,U$ isomorphic to
a dihedral group of order $12$.
\end{proof}

There is a $\Z_2$ symmetry among the irreducible modules of
$L(25/28,0)$. It is given by
\begin{equation*}
\rho=
\begin{cases}
1& \text{ on }  L(25/28,h_{r,s})\text{ for odd } s,\\
-1& \text{ on } L(25/28,h_{r,s})\text{ for even } s,
\end{cases}
\end{equation*}
where $h_{r,s} = h_{r,s}^5$ is defined by \eqref{eq:VIRH}. On the
other hand, $U'$ defines an automorphism $\tau$ of order $3$ (cf.
the $3A$ case). If $V$ is a VOA which contains a subalgebra
isomorphic to $U$, then there is an automorphism of order $6$
defined by $\tau\rho$. In fact, $\tau_{\hat{e}}\tau_{\hat{f}}=
\tau\rho $ in this case.

\begin{theorem}
If the Moonshine VOA $V^\natural$ contains a subalgebra isomorphic
to $U$, then as an automorphism of $V^\natural$,
$\tau_{\hat{e}}\tau_{\hat{f}}$ is of class $6A$.
\end{theorem}

\begin{proof}
Since $\tau_{\hat{e}}\tau_{\hat{f}}= \tau\rho $, we have
$(\tau_{\hat{e}}\tau_{\hat{f}})^2= \tau$ is of class $3A$ and
$(\tau_{\hat{e}}\tau_{\hat{f}})^3= \rho$ is of class $2A$. Thus
$\tau_{\hat{e}}\tau_{\hat{f}}$ is of class $6A$.
\end{proof}

\subsection{$4B$ case.} In this case, $L=L(6)\cong A_1\oplus A_7$,
$n_6=|E_8/L|=4$, and the conformal vectors $\tom^1\in
V_{\sqrt{2}A_1}$ and $\tom^2\in V_{\sqrt{2}A_7}$ defined by
\eqref{eq:CV} are of central charge $1/2$ and $7/5$, respectively.
Let $\tilde{\al}=2\al_6$. Then $E_7 = A_7 \cup (\tilde{\al}+
A_7)$. Thus we obtain a conformal vector $\tom(E_7)$ of central
charge $7/10$. Let $w^1=\tom^1$, $w^2 =\tom(E_7)$, and $w^3=
\tom^2- w^2$. Then $\{w^1 ,w^2, w^3\}$ is a set of mutually
orthogonal conformal vectors of central charge $1/2$, $7/10$, and
$7/10$, respectively. Actually,
\begin{equation*}
w^1= \tom^1,\qquad w^2= \frac{1}2\tom^2+ \frac1{20}X^2,\qquad w^3=
\frac{1}2\tom^2 - \frac1{20}X^2.
\end{equation*}

\begin{lemma}
Let $u= X^1+X^3$ and $v=X^1-X^3$. Then $u$ and $v$ are highest
weight vectors of highest weight $(1/2, 3/2,0)$ and $(1/2,0,3/2)$
with respect to $\Vir(w^1)\otimes \Vir(w^2)\otimes \Vir(w^3)$,
respectively.
\end{lemma}

\begin{proof}
It follows from \eqref{4b-eq} that
\begin{equation*}
(w^2)_1 X^1 = \frac{3}{4}X^1 + \frac{3}4 X^3,\qquad (w^2)_1 X^3 =
\frac{3}{4}X^1 + \frac{3}4 X^3.
\end{equation*}
Hence we have $(w^2)_1 u= \frac{3}2 u$ and $(w^2)_1 v =0$.
Similarly, $(w^3)_1 u=0$ and $(w^3)_1 v= \frac{3}2 v$.
\end{proof}

For simplicity of notation, we denote $L(1/2,h_1) \otimes
L(7/10,h_2) \otimes L(7/10,h_3)$ by $[h_1,h_2,h_3]$. The following
proposition is an immediate consequence of the above lemma.

\begin{proposition}
As a module of $\Vir(w^1)\otimes \Vir(w^2)\otimes \Vir(w^3)$,
\begin{equation*}
U \cong [0,0,0]
    \oplus [\frac{1}{2},\frac{3}{2},0]
    \oplus [\frac{1}{2},0,\frac{3}{2}]
    \oplus [0,\frac{3}{2},\frac{3}{2}].
\end{equation*}
\end{proposition}

The coset subalgebra $U$ contains four more sets $\{x^j, y^j,
z^j\}, 0 \le j \le 3$, of three mutually orthogonal conformal
vectors of central charge $1/2$, $7/10$, $7/10$, respectively (cf.
Appendix A). They are given by
\begin{align*}
x^j &= \frac{1}{8}\tom^1+\frac{5}{16}\tom^2 +
\frac{1}{32}\big((\sqrt{-1})^j X^1 + (-1)^j X^2 + (-\sqrt{-1})^j X^3\big),\\
y^j &= \quad \qquad \frac{1}{2}\tom^2 \qquad \qquad \qquad
\qquad\ - (-1)^j \frac{1}{20}X^2,\\
z^j &= \frac{7}{8}\tom^1+\frac{3}{16}\tom^2
-\frac{1}{32}\big((\sqrt{-1})^j X^1 - \frac{3}{5}(-1)^j X^2 +
(-\sqrt{-1})^j X^3\big).
\end{align*}

Note that $x^0 = \hat{e}$. The following proposition can be
verified by direct computation.

\begin{proposition}
As a module of $\Vir(x^j)\otimes \Vir(y^j)\otimes \Vir(z^j)$,
\begin{equation*}
U \cong [0,0,0] \oplus [\frac{1}{2},0,\frac{3}{2}] \oplus
[\frac{1}{16},\frac{3}{2},\frac{7}{16}]
\end{equation*}
if $j=0,2$ and
\begin{equation*}
U \cong [0,0,0] \oplus [\frac{1}{2},\frac{3}{2},0] \oplus
[\frac{1}{16},\frac{7}{16},\frac{3}{2}]
\end{equation*}
if $j=1,3$.
\end{proposition}

\begin{theorem}
$(1)$\ The Griess algebra $U_2$ of $U$ is generated by $\hat{e}$
and $\hat{f}$.

$(2)$\ The coset subalgebra $U$ is generated by $\hat{e}$ and
$\hat{f}$.
\end{theorem}

\begin{proof}
Let $\mathcal{G}$ be the Griess subalgebra generated by $\hat{e}$
and $\hat{f}$. By \eqref{4b-eq}, it is easy to  verify that
$\hat{e}$, $\hat{f}$, $\hat{e}_1\hat{f}$, $\hat{e}_1(\hat{e}_1
\hat{f})$, $\hat{f}_1(\hat{e}_1\hat{f})$ are linearly independent.
Thus $\mathcal{G}=U_2$, since $\dim U_2 = 5$. By the structure of
$U$, it is east to show that $U$ is generated by $U_2$. Hence the
second assertion holds.
\end{proof}

\begin{theorem}
The automorphism group $\mathrm{Aut}\, U $ of $U$ is a dihedral
group of order $8$.
\end{theorem}

\begin{proof}
There are exactly five conformal vectors of central charge $1/2$
in $U$, namely, $\tom^1$ and $e_j=\sigma^j \hat{e}$, $0 \le j \le
3$. Moreover,  we have $\la \tom^1, e_j\ra= 1/32$ for any
$j=0,1,2,3$, and
\begin{equation*}
\la e_i, e_j \ra =
\begin{cases}
1/32 & \text{ if } i+j\equiv 0 \mod 2, \\
1/2^8 & \text{ if } i+j \equiv 1 \mod 2.
\end{cases}
\end{equation*}

Since $U$ is generated by $\hat{e}$ and $\hat{f}=\sigma\hat{e}$,
$\Aut\,U$ can be considered as a subgroup of the permutation group
on the set $\{\tom^1, e_0, e_1, e_2, e_3\}$. Furthermore, $\Aut\,
U$ must preserve the inner product $\la \,\cdot\,,\,\cdot\,\ra$ so
that $\Aut\,U$ fixes $\tom^1$ and $|\Aut\,U| \le 8$. Since
$\sigma$ and $\theta$ generate a subgroup isomorphic to a dihedral
group of order $8$ in $\mathrm{Aut}\,U$, we have the assertion.
\end{proof}

The set of all irreducible modules of $U$ can be classified easily
by using the same method as in \cite{lam,lly,Y}. They are given by
\begin{equation*}
  \begin{array}{ll}
    [0,0,0]
    \oplus [0,\fr{3}{2},\fr{3}{2}]
    \oplus [\fr{1}{2},\fr{3}{2},0]
    \oplus [\hf,0,\fr{3}{2}],
    &
    {} [0,0,\fr{3}{2}]
    \oplus [0, \fr{3}{2}, 0]
    \oplus [\hf,0,0]
    \oplus [\hf,\fr{3}{2},\fr{3}{2}],
    \vsb\\
    {} [0,0,\fr{3}{5}]
    \oplus [0,\fr{3}{2},\fr{1}{10}]
    \oplus [\hf,\fr{3}{2},\fr{3}{5}]
    \oplus [\hf,0,\fr{1}{10}],
    &
    {} [0,\fr{3}{5},0]
    \oplus [0,\fr{1}{10},\fr{3}{2}]
    \oplus [\hf,\fr{3}{5},\fr{3}{2}]
    \oplus [\hf,\fr{1}{10},0],
    \vsb\\
    {} [0,\fr{3}{5},\fr{3}{5}]
    \oplus [0,\fr{1}{10},\fr{1}{10}]
    \oplus [\hf,\fr{1}{10},\fr{3}{5}]
    \oplus [\hf,\fr{3}{5},\fr{1}{10}],
    &
    {} [0,\fr{3}{5},\fr{1}{10}]
    \oplus [0,\fr{1}{10},\fr{3}{5}]
    \oplus [\hf,\fr{1}{10},\fr{1}{10}]
    \oplus [\hf,\fr{3}{5},\fr{3}{5}],
    \vsb\\
    {} [0,0,\fr{1}{10}]
    \oplus [0,\fr{3}{2},\fr{3}{5}]
    \oplus [\hf,\fr{3}{2},\fr{1}{10}]
    \oplus [\hf,0,\fr{3}{5}],
    &
    {} [0,\fr{1}{10},0]
    \oplus [0,\fr{3}{5},\fr{3}{2}]
    \oplus [\hf,\fr{3}{5},0]
    \oplus [\hf,\fr{1}{10},\fr{3}{2}],
  \end{array}
\end{equation*}
and
\begin{equation*}
  \begin{array}{llll}
    {} [\fr{1}{16},\fr{7}{16},\fr{7}{16}]\tensor Q,
    &  [\fr{1}{16},\fr{7}{16},\fr{3}{80}]\tensor Q,
    &  [\fr{1}{16},\fr{3}{80},\fr{7}{16}]\tensor Q,
    &  [\fr{1}{16},\fr{3}{80},\fr{3}{80}]\tensor Q,
  \end{array}
\end{equation*}
where $Q$ is the unique $2$-dimensional irreducible module of the
quaternion group of order $8$.

The fixed point subalgebra $U^{\tau_{\hat{e}}\tau_{\hat{f}}}$ of
$\tau_{\hat{e}}\tau_{\hat{f}}$ in $U$ is isomorphic to $[0,0,0]
\oplus [0,3/2,3/2]$.

The fusion rules among irreducible modules of $L(7/10,0)\otimes
L(7/10,0)\oplus L(7/10,3/2)\otimes L(7/10,3/2)$ can be computed
easily. There is a $\Z_4$ symmetry given as follows.
\begin{equation*}
\tau=
\begin{cases} \ds
1 & \text{ on }\  [0,0]
    \oplus [\fr{3}{2},\fr{3}{2}],\
    [0,\fr{3}{5}]
    \oplus [\fr{3}{2},\fr{1}{10}], \\ \ds
    & \qquad [\fr{3}{5},\fr{3}{5}]\oplus
    [\fr{1}{10},\fr{1}{10}],\
    [\fr{1}{10},\fr{3}{2}]
    \oplus [\fr{3}{5},0],
    \\
    \ \\
\sqrt{-1} &\text{ on }\ [\fr{7}{16},\fr{7}{16}]^+, \
 [\fr{7}{16},\fr{3}{80}]^+, \\ \ds
    & \qquad [\fr{3}{80},\fr{7}{16}]^+,\
     [\fr{3}{80},\fr{3}{80}]^+,\\

     \ \\ \ds
 -1& \text{ on } \  [0,\fr{3}{2}]
    \oplus [\fr{3}{2},0],\
    [\fr{3}{2},\fr{3}{5}]
    \oplus [0,\fr{1}{10}], \\ \ds
    & \qquad [\fr{3}{5},\fr{1}{10}]\oplus
    [\fr{1}{10},\fr{3}{5}],\
    [\fr{1}{10},0]
    \oplus [\fr{3}{5},\fr{3}2],
    \\
\ \\
     \ds
 -\sqrt{-1} & \text{ on } \ [\fr{7}{16},\fr{7}{16}]^-, \
 [\fr{7}{16},\fr{3}{80}]^-, \\ \ds
    & \qquad [\fr{3}{80},\fr{7}{16}]^-,\
     [\fr{3}{80},\fr{3}{80}]^-,
\end{cases}
\end{equation*}
where $[h_1,h_2]$ denotes the irreducible module
$L(7/10,h_1)\otimes L(7/10,h_2)$ of $L(7/10,0)\otimes L(7/10,0)$.
Suppose $U$ is contained in a VOA $V$. Then all $\tau,
\tau_{\hat{e}}$, and $\tau_{\hat{f}}$ are well defined
automorphisms of $V$. In this case, $\tau=
\tau_{\hat{e}}\tau_{\hat{f}}$ and $\tau^2= \tau_{w^1}$. Thus we
have the following theorem.

\begin{theorem}
If the Moonshine VOA $V^\natural$ contains a subalgebra isomorphic
to $U$, then as an automorphism of $V^\natural$,
$\tau_{\hat{e}}\tau_{\hat{f}}$ is of class $4B$.
\end{theorem}

\subsection{$2B$ case.} In this case, $L=L(7) \cong D_8$,
$n_7=|E_8/L|=2$, and the conformal vector $\tilde{\omega}^1 \in
V_{\sqrt{2}D_8}$ defined by \eqref{eq:CV} is of central charge
$1$. The Virasoro element $\omega'$ of $U$ is equal to
$\tilde{\omega}^1$. Let $\ep_1,\ep_2, \dots, \ep_8 \in \R^8$ be
such that $\la \ep_i,\ep_j\ra=2\delta_{ij}$ for any $i,j$. Then
\[
\sqrt{2}D_8=\{a_1\ep_1 +\cdots +a_8\ep_8\, |\, a_i \in \Z, a_1 +
\cdots +a_8 \equiv 0 \mod 2\}.
\]

Let $\gamma=\frac{1}2 ( \ep_1+ \cdots+ \ep_8)$. Then
$\sqrt{2}E_8\cong \sqrt{2}D_8 \cup (\gamma + \sqrt{2}D_8)$. Denote
\[
U^0=\{u\in V_{\sqrt{2}D_8}\,|\, s(D_8)_1 u=0\}, \qquad U^1=\{u\in
V_{\gamma + \sqrt{2}D_8}\,|\, s(D_8)_1 u=0\},
\]
where $s(D_8)$ is defined as in \eqref{cv1}. By \cite{dly,dly2},
we know that $U^0\cong V^{+}_{\mathbb{Z}2\gamma}$ and $U^1\cong
V^{+}_{\gamma+\mathbb{Z}2\gamma}$. Note that $\Z\gamma = \Z2\gamma
\cup (\gamma+ \Z2\gamma)$ and hence we have $U=U^0 \oplus U^1
\cong V_{\Z\gamma}^+$. Since  $\la\gamma, \gamma\ra =4$, it is
well know that
\[
V_{\Z\gamma}^+\cong L(\frac{1}2,0)\otimes L(\frac{1}2,0).
\]
In fact, $\frac{1}{16} \gamma(-1)^2 \cdot 1+ \frac{1}4 (e^\gamma
+e^{-\gamma})$ and $\frac{1}{16} \gamma(-1)^2 \cdot 1 - \frac{1}4
(e^\gamma +e^{-\gamma})$ are the two mutually orthogonal conformal
vectors of central charge $1/2$ in $V_{\Z\gamma}^+$ (cf.
\cite{dmz,m1,m2}).

The following theorem is clear from the structure of $U$.

\begin{theorem}
The coset subalgebra $U$ is generated by $\hat{e}$ and $\hat{f}$.
Moreover, the automorphism group $\mathrm{Aut}\, U $ of $U$ is of
order $2$.
\end{theorem}

Both of $\hat{e}$ and $\hat{f}$ are fixed by $\theta$ and thus $U$
is contained in $V_\Lambda^+\subset V^\natural$. In this case, $
\hat{e}$ and $\hat{f}$ are mutually orthogonal and so
$\tau_{\hat{e}}$ and $\tau_{\hat{f}}$ are commutative. Hence,
$|\tau_{\hat{e}}\tau_{\hat{f}}|=2$ and as an automorphism of
$V^\natural$, $\tau_{\hat{e}}\tau_{\hat{f}}$ is of class $2B$ (see
also Shimakura \cite{Shi}).

\subsection{$3C$ case.}\label{sub5-3c}

In this case, $L=L(8)\cong A_8$, $n_8=|E_8/L|=3$, and the
conformal vector $\tilde{\omega}^1 \in V_{\sqrt{2}A_8}$ defined by
\eqref{eq:CV} is of central charge $16/11$. The Virasoro element
$\omega'$ of $U$ is equal to $\tilde{\omega}^1$. The following
lemma can be easily verified.

\begin{lemma}
Let
\[
x=\frac{11}{32} \tom^1 +\frac{1}{32}(X^1+X^2),\qquad y=
\frac{21}{32}\tom^1 -\frac{1}{32}(X^1+X^2).
\]
Then $x$ and $y$ are mutually orthogonal conformal vectors of
central charge $1/2$ and $21/22$, respectively. Moreover,
$\tom^1=x+y$.
\end{lemma}

The above lemma implies that $U$ contains $\Vir(x) \otimes \Vir(y)
\cong L(1/2,0)\otimes L(21/22,0)$.

\begin{lemma}\label{31-16}
Let $v=X^1-X^2$. Then $v$ is a highest weight vector of highest
weight $(1/16,31/16)$ with respect to $\Vir(x) \otimes \Vir(y)$.
\end{lemma}

\begin{proof}
By \eqref{3c-eq}, we have
\begin{equation*}
x_1 v = \Big(\frac{11}{32} \tom^1 +\frac{1}{32}(X^1+X^2)\Big)_1
(X^1-X^2) = \frac{1}{16} (X^1-X^2).
\end{equation*}
Hence the assertion holds.
\end{proof}

Similarly, we obtain the following lemma by direct computation.

\begin{lemma}\label{7-2}
Let
\[
  u =4 v_{-1}v -45 x_{-3}\cdot 1 -9 x_{-1}x
    -682 x_{-1}y -1089 y_{-1}y -165 y_{-3}\cdot 1.
\]
Then $u$ is a highest weight vector of highest weight $(1/2,7/2)$
with respect to $\Vir(x) \otimes \Vir(y)$ and $\la u,u\ra
=1280076$.
\end{lemma}

Note that $L(1/2,h_1) \otimes L(21/22,h_2)$ with
$(h_1,h_2)=(0,0)$, $(0,8)$, $(1/2,7/2)$, $(1/2,45/2)$,
$(1/16,31/16)$ and $(1/16,175/16)$ are the only irreducible
$L(1/2,0)\otimes L(21/22,0)$-modules which are integrally graded.
For simplicity of notation, we shall denote $L(1/2,h_1)
\otimes L(21/22,h_2)$ by $[h_1,h_2]$. The following theorem can be
obtained by using the characters (cf. Appendix B).

\begin{theorem}[Theorem \ref{B-3C}]
As a module of $\Vir(x)\otimes \Vir(y)$,
\begin{equation*}
U \cong  [0,0] \oplus [0,8] \oplus [\frac{1}{2},\frac{7}{2}]
\oplus [\frac{1}{2},\frac{45}{2}] \oplus
[\frac{1}{16},\frac{31}{16}] \oplus [\frac{1}{16},\frac{175}{16}].
\end{equation*}
\end{theorem}

\medskip

\begin{theorem}\label{thm:3.36}
$(1)$\ The Griess algebra $U_2$ of $U$ is generated by $\hat{e}$
and $\hat{f}$.

$(2)$\ The coset subalgebra $U$ is generated by $\hat{e}$ and
$\hat{f}$.
\end{theorem}

\begin{proof}
The first assertion is clear since $\dim U_2 = 3$. Let $W$ be the
subalgebra generated by $U_2$. We want to show that $W = U$. By
Lemmas \ref{31-16} and \ref{7-2}, $W$ contains submodules of the
form $[0,0]$, $[1/2,7/2]$, and $[1/16,31/16]$.

We shall show that there is a highest weight vector of highest
weight $(0,8)$ with respect to $\Vir(x) \otimes \Vir(y)$ in $W$.
Suppose $W$ does not contain a highest weight vector of highest
weight $(0,8)$. Then the fixed point subalgebra $W^{\tau_x}$ of
the Miyamoto involution $\tau_x$ in $W$ is of the form $W^{\tau_x}
\cong [0,0] \oplus [1/2,7/2]$. Let $w$ be a highest weight vector
of highest weight $(1/2,7/2)$. We normalize $w$ so that $\la
w,w\ra=1$. For example, we may take $w= u /\sqrt{\la u,u \ra}$.
Then by fusion rules, $w_{-1}w \in [0,0]$. It is well known (cf.\
\cite{FFR} \cite{FRW}) that there is an explicit construction of
the vertex operator superalgebra $L(1/2,0)\oplus L(1/2,1/2)$ by
using one free fermionic field and we can find an orthonormal
basis of $L(1/2,0)\oplus L(1/2,1/2)$ (cf. \cite{kr}). By the
assumption $W^{\tau_x} \cong [0,0] \oplus [1/2,7/2]$ is
$\Z_2$-graded and so we can use the orthonormal basis of
$L(1/2,0)\oplus L(1/2,1/2)$ to make a computation of the inner
product $\la w_{-1}w,w_{-1}w\ra$ easier. As a result, we can
obtain
\[
\la w_{-1}w, w_{-1}w\ra = \dfr{231091005602}{134258427}.
\]
Note that this value is deduced based on the assumption that
$[0,8]$ is not contained in $W$ as a $\Vir(x) \otimes
\Vir(y)$-submodule.

On the other hand, we can also compute $\la w_{-1}w, w_{-1}w\ra$
directly by using the definition of $w$ and the Jacobi identity.
In this case, we have
\begin{equation*}
\begin{split}
\la w_{-1}w,w_{-1}w\ra
  &= \la w, w_{7}w_{-1}w\ra = \la w, (w_{-1}w)_{7}w\ra
    \\
  & = \dsum_{n\geq 0} (-1)^n \binom{-1}{n} \la w, (w_{-1-n}
    w_{7+n}+w_{6-n}w_{n}) w\ra \\
  & = \la w,w_{-1}w_{7}w\ra +\dsum_{n=0}^7 \la w,w_{6-n}
    w_{n}w\ra\\
    &= 2 \la w,w\ra^2 +\dsum_{n=0}^5 \la
    w_{n}w,w_{n}w\ra.
    \end{split}
\end{equation*}

Then by computing $\la w_{5-m} w,w_{5-m} w\ra$, $m=0,1,\dots,5$,
inductively, we obtain another value
\[
\la w_{-1}w, w_{-1}w\ra = \frac{1036050}{589} \left( \neq
\dfr{231091005602}{134258427} \right).
\]
Note that this value is deduced from the structure of the Griess
algebra $U_2$ and independent of the shape of $W$ as a
$\Vir(x)\tensor \Vir (y)$-module. Thus this contradiction comes
from the assumption that $W$ does not contain a highest weight
vector of highest weight $(0,8)$. Hence we conclude that there is
a highest weight vector of highest weight $(0,8)$ in $W$.
By the computation above, we also note that $w_{-1}w$ contains
a singular vector with highest weight $(0,8)$ as a non-trivial summand.

Then, $W$ also contains a highest weight vector of highest weight
$(1/2,45/2)$; otherwise, $W^{\tau_x} \cong [0,0] \oplus [1/2,7/2]
\oplus [0,8]$ and the orthogonal complement of $W^{\tau_x}$ in
$U^{\tau_x}$ with respect to a positive definite invariant
hermitian form is isomorphic to $[1/2,45/2]$. But the orthogonal
complement must be a module for $W^{\tau_x}$, which is impossible
by the fusion rules. Thus $W^{\tau_x} = U^{\tau_x}$.

Since the fusion product of $[1/2,45/2]$ and $[1/16,31/16]$ is
$[1/16,175/16]$, $W$ contains a highest weight vector of highest
weight $(1/16,175/16)$ also, and hence $W = U$ as desired.
\end{proof}

\begin{theorem}
The automorphism group $\mathrm{Aut}\, U $ of $U$ is a symmetric
group $S_3$ of degree $3$.
\end{theorem}

\begin{proof}
By Lemma \ref{A-3C}, there are exactly three conformal vectors of
central charge $1/2$ in $U$, namely, $\sigma^j \hat{e}$,
$j=0,1,2$. Thus $\mathrm{Aut}\,U$ can be considered as a subgroup
of $S_3$. Since $\sigma$ and $\theta$ already generate a subgroup
of order $6$ in $\mathrm{Aut}\,U$, we conclude that
$\mathrm{Aut}\,U\cong S_3$.
\end{proof}

\begin{theorem}
There are exactly five irreducible $U$-modules $U(2k), 0 \le k \le
4$. In fact, $U(0)=U$ and as $\Vir(x)\otimes \Vir(y)$-modules,
\begin{align*}
U(2) & \cong [0,\frac{13}{11}] \oplus [0,\frac{35}{11}] \oplus
[\frac{1}{2},\frac{15}{22}] \oplus [\frac{1}{2},\frac{301}{22}]
\oplus [\frac{1}{16},\frac{21}{176}] \oplus [\frac{1}{16},
\frac{901}{176}],\\
U(4) & \cong [0,\frac{6}{11}]  \oplus [0,\frac{50}{11}] \oplus
[\frac{1}{2},\frac{1}{22}] \oplus [\frac{1}{2},\frac{155}{22}]
\oplus [\frac{1}{16},\frac{85}{176}] \oplus [\frac{1}{16},
\frac{261}{176}],\\
U(6) & \cong  [0,\frac{1}{11}] \oplus [0,\frac{111}{11}] \oplus
[\frac{1}{2},\frac{35}{22}] \oplus [\frac{1}{2},\frac{57}{22}]
\oplus [\frac{1}{16}, \frac{5}{176}]
\oplus [\frac{1}{16},\frac{533}{176}],\\
U(8) &\cong [0,\frac{20}{11}] \oplus [0,\frac{196}{11}] \oplus
[\frac{1}{2},\frac{7}{22}] \oplus [\frac{1}{2},\frac{117}{22}]
\oplus [\frac{1}{16},\frac{133}{176}] \oplus
[\frac{1}{16},\frac{1365}{176}].
\end{align*}
\end{theorem}

We shall again give a proof at Appendix C. Note that the lattice
VOA $V_{\sqrt{2}E_8}$ can be decomposed as follows (cf.
\cite{lam-3c,ly3}).
\begin{align}
V_{\sqrt{2}E_8} \cong
 & \bigoplus_{ { 0\leq k_j\leq j+1}
 \atop{{ k_j \equiv 0\,\mathrm{mod}\, 2} \atop{ j=0,1,\dots, 8}}}
  L(c_1,h^1_{k_{0}+1,k_1+1})\otimes \cdots  L(c_8,h^8_{k_{7}+1,k_8+1})
  \otimes U(k_8),
\end{align}
where $c_m$ and $h_{r,s}^m$ are given by \eqref{eq:VIRC} and
\eqref{eq:VIRH}. Moreover, as an automorphism of
$V_{\sqrt{2}E_8}$, $\tau_{\hat{e}} \tau_{\hat{f}}= e^{2\pi
\sqrt{-1} \bbe(0)}$ is of order $3$ (cf. Remark \ref{tau-ef}) and
it induces a natural action on each of $U(2k)$, $0\le k\le 4$.
Hence, for any VOA $V$ which contains a subalgebra isomorphic to
$U$, $\tau_{\hat{e}} \tau_{\hat{f}}$ defines an automorphism of
order $3$ on $V$.

\begin{remark}
In Miyamoto\cite{m6}, it is shown that if $e$ and $f$ are two
conformal vectors of central charge $1/2$ in the Moonshine VOA
$V^\natural$ such that $\tau_{e} \tau_f$ is of order $3$ and $\la
e,f \ra =1/2^8$, then the subalgebra $W$ generated by $e$ and $f$
must contain a subalgebra isomorphic to $L(1/2,0)\otimes
L(21/22,0)$ and the weight $2$ subspace $W_2$ of $W$ is of
dimension $3$. In fact, $W\cong U$ as $L(1/2,0) \otimes
L(21/22,0)$-modules and $\tau_{e}\tau_{f}$ is of class $3C$.
\end{remark}

\appendix
\section{ Conformal Vectors in $U$ }

In this appendix, we shall compute the conformal vectors of the
coset subalgebra $U$ defined by \eqref{eq:UDEF} in each of the
nine cases. Except for the cases of $5A$ and $6A$, the computation
was done by Maple 7. The cases for $5A$ and $6A$ were computed by
Kazuhiro Yokoyama of Kyushu University using a computer algebra
system Risa/Asir.

As in Section \ref{S3}, $L=L(i)$ denotes  the lattice associated
with the Dynkin diagram obtained by removing the $i$-th node
$\al_i$ in the extended $E_8$ diagram \eqref{ext} and $n_i =
|E_8/L|$. Let
\[
\B=\mathrm{span}_\C \{ \tom^1, , \dots,\tom^l, X^1,\dots,
X^{n_i-1}\},
\]
where $\tilde{\omega}^j$ and $X^j$ are defined by \eqref{eq:CV}
and \eqref{eq:HWV}. Then $\dim \B = l + n_i - 1$. Actually, $\B$
is equal to the Griess algebra $U_2$ of $U$. Note that
\begin{equation}\label{eq:GC}
\begin{split}
& (\tilde{\omega}^j)_1\tilde{\omega}^k = 2 \tilde{\omega}^j
\delta_{j,k},\qquad \la \tilde{\omega}^j, \tilde{\omega}^k \ra =
c\delta_{j,k}/2, \qquad
\la \tilde{\omega}^j, X^k \ra = 0,\\
& \la X^j, X^k \ra = 0 \quad \text{if} \quad j+k \not\equiv 0 \mod
n_i,
\end{split}
\end{equation}
where $c$ denotes the central charge of $\tilde{\omega}^j$.

\bigskip \textbf{$1A$ case.} In this case, $L=L(0)\cong E_8$,
$n_0=1$, and $\tilde{\om}^1 \in V_{\sqrt{2}E_8}$ is the only
conformal vector in $\B$, whose central charge is $1/2$.

\bigskip \textbf{$2A$ case.} In this case, $ L=L(1)\cong A_1\oplus
E_7$, $n_1=2$, and the conformal vectors $\tom^1 \in
V_{\sqrt{2}A_1}$ and $\tom^2 \in V_{\sqrt{2}E_7}$ are of central
charge $1/2$ and $7/10$, respectively. The product and the inner
product in $\B$ are given by \eqref{eq:GC} and
\[
(\tilde{\omega}^1)_1 X^1 = \frac{1}{2}X^1,\quad
(\tilde{\omega}^2)_1 X^1 = \frac{3}{2}X^1,\quad (X^1)_1X^1 =
224\tom^1+ 480\tom^2, \quad \la X^1, X^1\ra =112.
\]

\begin{lemma}
Let $w=a\tom^1 +b\tom^2 +cX^1$ be a conformal vector in $\B$. Then
$(a,b,c)$ satisfies the following system of equations.
\begin{equation}\label{2a}
\begin{array}{ll}
a^2+112c^2-a=0, & \qquad
b^2+240c^2-b=0,\\
ac+3bc-2c=0. & \
\end{array}
\end{equation}
The central charge of $w$ is given by
$\frac{1}2a^2+\frac{7}{10}b^2+224c^2$.
\end{lemma}

The solutions $(a,b, c)$ of Equation \eqref{2a} are as follows.

\medskip
\noindent Central charge $1/2$: $ (1,0,0), \ (\frac{1}{8},
\frac{5}{8}, \frac{1}{32}), \ (\frac{1}{8}, \frac{5}{8},
-\frac{1}{32})$.

\medskip
\noindent Central charge $7/10$: $ (0,1,0), \ (\frac{7}{8},
\frac{3}{8}, \frac{1}{32}), \ (\frac{7}{8}, \frac{3}{8},
-\frac{1}{32})$.

\medskip
\noindent Central charge $6/5$: $(1,1,0)$.

\bigskip \textbf{$3A$ case.} In this case, $ L=L(2)\cong A_2\oplus
E_6$, $n_2=3$, and the conformal vectors $\tom^1 \in
V_{\sqrt{2}A_2}$ and $\tom^2 \in V_{\sqrt{2}E_6}$ are of central
charge $4/5$ and $6/7$, respectively. The product and the inner
product in $\B$ are given by \eqref{eq:GC} and
\begin{equation*}
\begin{split}
&(\tom^1)_1 X^1 =\frac{2}3 X^1,\quad  (\tom^1)_1 X^2 =\frac{2}3
X^2,\quad (\tom^2)_1 X^1 =\frac{4}3 X^1,\quad(\tom^2)_1 X^2
=\frac{4}3 X^2,\\
&(X^1)_1X^1=20 X^2, \quad   (X^2)_1X^2=20 X^1,\quad (X^1)_1X^2=135
\tom^1+ 252\tom^2,\\
& \la X^1, X^2\ra =81.
\end{split}
\end{equation*}

\begin{lemma}
Let $w=a\tom^1 +b\tom^2 +cX^1 +dX^2$ be a conformal vector in
$\B$. Then $(a,b,c,d)$ satisfies the following system of
equations.
\begin{equation}\label{3a}
\begin{array}{ll}
a^2+135cd-a=0,
& \qquad b^2+252cd-b=0,\\
15d^2 + 2ac + 4bc-3c=0, & \qquad 15c^2+ 2ad +4bd-3d=0.
\end{array}
\end{equation}
The central charge of $w$ is given by $\frac{4}5a^2+\frac{6}7
b^2+324cd$.
\end{lemma}

The solutions $(a,b,c,d)$ of Equation (\ref{3a}) are as follows,
where $\xi=e^{2\pi \sqrt{-1}/3}$ is a primitive cubic root of
unity.

\medskip
\noindent Central charge $1/2$  : $\big(\frac{5}{32},\frac{7}{16},
\frac{1}{32}\xi^j, \frac{1}{32}\xi^{2j}\big)$, $j=0,1,2$.

\medskip
\noindent Central charge $4/5$:  $(1,0,0,0)$, \quad
$\big(\frac{1}{16}, \frac{7}{8}, -\frac{1}{48}\xi^j,
-\frac{1}{48}\xi^{2j}\big)$, $j=0,1,2$.

\medskip
\noindent Central charge $6/7$: $(0,1,0,0)$,\quad
$\big(\frac{15}{16}, \frac{1}{8}, \frac{1}{48}\xi^j,
\frac{1}{48}\xi^{2j}\big)$, $j=0,1,2$.

\medskip
\noindent Central charge $81/70$: $\big(\frac{27}{32},
\frac{9}{16}, - \frac{1}{32}\xi^j, - \frac{1}{32}\xi^{2j}\big)$,
$j=0,1,2$.

\medskip
\noindent Central charge $58/35$: $(1,1,0,0)$.

\bigskip \textbf{$4A$ case.} In this case,  $L= L(3) \cong
A_3\oplus D_5$, $n_3=4$, and the conformal vectors $\tom^1\in
V_{\sqrt{2}A_3}$ and $\tom^2\in V_{\sqrt{2}D_5} $ are both of
central charge $1$. The product and the inner product in $\B$ are
given by \eqref{eq:GC} and
\begin{equation*}
\begin{array}{lll}
  (\tom^1)_{1} X^1 = \dfr{3}{4} X^1,
  &  (\tom^1)_{1} X^2 = X^2,\q
  &  (\tom^1)_{1} X^3 =\dfr{3}{4} X^3,
  \vsb\\
    (\tom^2)_{1} X^1 = \dfr{5}{4} X^1,
  & (\tom^2)_{1} X^2 = X^2,
  & (\tom^2)_{1} X^3 = \dfr{5}{4} X^3,
    \vspace{3mm}\\
    (X^1)_{1} X^1 = 16 X^2,
  & (X^1)_{1} X^2 = 15 X^3,
  & (X^1)_{1} X^3= 96 \tom^1+160 \tom^2
  \vsb\\
    (X^2)_{1} X^2 = 120 \tom^1 +120 \tom^2,\q
  & (X^2)_{1} X^3 = 15 X^1,
  &( X^3)_{1} X^3 = 16 X^2,
  \vsb\\
  \la X^1, X^3\ra = 64, & \la X^2,X^2\ra = 60. & \
\end{array}
\end{equation*}

\begin{lemma}
Let $w=a\tom^1 +b\tom^2 +cX^1 +dX^2 +eX^3$ be a conformal vector
in $\B$. Then $(a,b,c,d,e)$ satisfies the following system of
equations.
\begin{equation}\label{4a}
\begin{array}{ll}
a^2+60d^2+96ce-a=0,
& \qquad b^2+60 d^2+160ce-b=0,\\
3ac+5bc+60de-4c=0,
& \qquad 8c^2+8e^2+ad+bd-d=0,\\
3ae+5be+60cd-4e=0. & \
\end{array}
\end{equation}
The central charge of $w$ is given by $a^2+b^2+256c e+120 d^2$.
\end{lemma}

The solutions  $(a,b,c,d,e)$ of Equation \eqref{4a} are as
follows, where $\xi = e^{\pi \sqrt{-1}/4}$ is a primitive $8$-th
root of unity.

\medskip
\noindent  Central charge $1/2$: $\big( \frac{3}{16},
\frac{5}{16}, \frac{1}{32}(\sqrt{-1})^j, \frac{1}{32}(-1)^j,
    \frac{1}{32}(-\sqrt{-1})^j\big)$, \ $0 \le j \le 3$.

\medskip
\noindent Central charge $1$:
\begin{equation*}
\begin{array}{l}
  (1,0,0,0,0),\quad  (0,1,0,0,0), \vsb\\
  \l( (1+\sqrt{1-240 d^2})/2, (1-\sqrt{1-240 d^2})/2, 0, d,0
  \r) ,
  \vsb\\
  \l( (1-\sqrt{1-240 d^2})/2, (1+\sqrt{1-240 d^2})/2, 0, d,0
  \r) , \ \mathrm{where}\ d \in \C\setminus \{ 0 \} .
\end{array}
\end{equation*}

\medskip
\noindent  Central charge $3/2$: $\big( \frac{13}{16},
\frac{11}{16}, \frac{1}{32}(\sqrt{-1})^j, -\frac{1}{32}(-1)^j,
    \frac{1}{32}(-\sqrt{-1})^j\big)$, \ $0 \le j \le 3$.

\medskip
\noindent  Central Charge $6/7$: $\big( \frac{1}{7}, \frac{5}{7},
\frac{1}{28}\xi^{2j+1}, 0,
    \frac{1}{28}\xi^{-(2j+1)}\big)$, \ $0 \le j \le 3$.

\medskip
\noindent  Central Charge $8/7$: $\big( \frac{6}{7}, \frac{2}{7},
\frac{1}{28}\xi^{2j+1}, 0,
    \frac{1}{28}\xi^{-(2j+1)}\big)$, \ $0 \le j \le 3$.

\medskip
\noindent  Central Charge $2$: $(1,1,0,0,0)$.

\bigskip
\textbf{$5A$ case.}
In this case, $L=L(4) \cong A_4\oplus A_4$, $n_4=5$, and the conformal vectors
$\tom^1 \in V_{\sqrt{2}A_4}$ and $\tom^2 \in V_{\sqrt{2}A_4}$ are both of
central charge $8/7$. The product and the inner product in $\B$
are given by \eqref{eq:GC} and
$$
\begin{array}{l}
  \begin{array}{llll}
    (\tom^1)_1 X^1=\dfrac{4}5 X^1,
    & (\tom^2)_1 X^1=\dfrac{6}5 X^1,
    & (\tom^1)_1 X^2=\dfrac{6}5 X^2,
    & (\tom^2)_1 X^2=\dfrac{4}5 X^2,
    \vsb\\
    (\tom^1)_1 X^3=\dfrac{6}5 X^3,
    & (\tom^2)_1 X^3=\dfrac{4}5 X^3,
    & (\tom^1)_1 X^4=\dfrac{4}5 X^4,
    & (\tom^2)_1 X^4=\dfrac{6}5 X^4,
  \end{array}
  \vsb\\
  (X^i)_1 (X^j) =
  \begin{cases}
    12 X^{i+j} & \text{ if } i+j\not\equiv 0 \mod 5,
    \\
    70\tom^1+ 105\tom^2 & \text{ if } i=5-j=1\text { or } 4,
    \\
    105\tom^1+ 70\tom^2 & \text{ if } i=5-j=2\text { or } 3,
  \end{cases}
  \vsb\\
  \la X^1,X^4\ra = \la X^2,X^3\ra = 50.
\end{array}
$$

\begin{lemma}\label{cv}
  Let $w=a\tom^1 +b\tom^2 +c X^1 + d X^2 +e X^3+f X^4$ be a
  conformal vector in $\B$. Then $(a,b,c,d,e,f)$ satisfies the
  following system of equations.
  \begin{equation}\label{5ae}
  \begin{array}{ll}
    a^2+70cf+105de-a=0,
    & b^2+105cf+70de-b=0,
    \vsb\\
    4ac+6bc+30e^2+60df-5c=0,
    & 6ad+4bd+30c^2+60ef-5d=0,
    \vsb\\
    6ae+4be+30f^2+60cd-5e=0,
    & 4af+6bf+30d^2+60ce-5f=0.
  \end{array}
  \end{equation}
  The central charge of $w$ is given by
  \begin{equation}\label{5acc}
    \frac{8}7 (a^2+b^2)+200(cf+de).
  \end{equation}
\end{lemma}

In order to solve the above system of equations, we treat
$a,b,c,d,e,f$ as variables. Let $\C[a,b,c,d,e,f]$ be the
polynomial algebra with variables $a,b,c,d,e,f$ and $\I$ the ideal
generated by the six polynomials which appear on the left hand
side of (\ref{5ae}). We then compute the primary decomposition of
the ideal $\I=\cap \CP_i$ over the field $\Q$ of rational numbers
and solve the system corresponding to each prime ideal $\CP_i$. By
using the computer algebra system Risa/Asir, we found that $\I$ is
an intersection of $21$ prime ideals and there are $63$ different
nontrivial solutions of (\ref{5ae}). Thus we have the following
lemma.

\begin{lemma}
  There are exactly $63$ conformal vectors in $\B$.
\end{lemma}

The central charges of those conformal vectors are easily
calculated by \eqref{5acc}. We verified that there are only $43$
conformal vectors whose central charges are rational numbers.
Their central charges are shown in Table \ref{cclist1}.

\begin{table}[ht]
\begin{center}
\caption{Central charge (c.c.) and number of conformal vectors}
\label{cclist1}
\medskip
\begin{tabular}{c|c|c|c|c|c|cc}
\hline c.c. & $1/2$ & $25/28$ & $8/7$ & $16/7$ & $25/14$  &
$39/28$ \\
 \hline
number & $5$ &$10$ &$12$ &$1$ &$5$ &$10$\\
\hline
\end{tabular}
\end{center}
\end{table}

The following two lemmas were verified by computer.

\begin{lemma}\label{A-5A}
There are exactly five conformal vectors of central charge $1/2$
in $\B$, namely, $\sigma^j \hat{e}, 0 \le j \le 4$. Note that
\[
\sigma^j \hat{e}=\frac{7}{32} (\tom^1+\tom^2) + \frac{1}{32}(
\xi^j X^1+ \xi^{2j}X^2+ \xi^{3j}X^3+ \xi^{4j}X^4),
\]
where $\xi=e^{2 \pi \sqrt{-1}/5}$ is a primitive $5$-th root of
unity. Moreover, the inner product is $\la \sigma^i \hat{e},
\sigma^j \hat{e}\ra = 3/512$ for any $i\neq j$.
\end{lemma}

\begin{lemma}
There is a triple $(u,v,w)$ of mutually orthogonal conformal
vectors in $\B$ such that the central charges of $u, v, w$ are
$1/2, 25/28$, $25/28$, respectively and $u+v+w$ is equal to the
Virasoro element $\tilde{\omega}^1 + \tilde{\omega}^2$ of $U$. For
example,
\begin{align*}
u &=\hat{e}=\frac{7}{32} (\tom^1+\tom^2) + \frac{1}{32}( X^1+X^2+X^3+X^4),\\
v &= \frac{15}{64}\tom^1+\frac{35}{64}\tom^2 -
\frac{3}{64}(X^1+X^4) +\frac{1}{64}(X^2+X^3),\\
w &= \frac{35}{64}\tom^1+\frac{15}{64}\tom^2 +
\frac{1}{64}(X^1+X^4) -\frac{3}{64}(X^2+X^3).
\end{align*}
\end{lemma}

\medskip
\textbf{$6A$ case.}
In this case,  $L= L(5) \cong A_2 \oplus A_1 \oplus A_5$, $n_5 =6$, and
the conformal vectors $\tom^1\in V_{\sqrt{2}A_2}$,
$\tom^2\in V_{\sqrt{2}A_1}$, and $\tom^3\in V_{\sqrt{2}A_5}$ are of
central charge $4/5$, $1/2$, and $5/4$, respectively.
The product and the inner product in $\B$ are given by \eqref{eq:GC} and
\begin{equation}\label{6a-eq}
\begin{array}{lll}
  (\tom^1)_{1} X^1=\dfr{2}{3} X^1,
    & (\tom^2)_{1} X^1=\dfr{1}{2} X^1,
    & (\tom^3)_{1} X^1=\dfr{5}{6} X^1,
    \vsb\\
  (\tom^1)_{1} X^2=\dfr{2}{3} X^2,
    & (\tom^2)_{1} X^2=0,
    & (\tom^3)_{1} X^2 =\dfr{4}{3} X^2,
    \vsb\\
  (\tom^1)_{1} X^3 =0,
    & (\tom^2)_{1} X^3=\dfr{1}{2} X^3,
    & (\tom^3)_{1} X^3 = \dfr{3}{2} X^3,
    \vsb\\
  (\tom^1)_{1} X^4=\dfr{2}{3} X^4,
    & (\tom^2)_{1} X^4=0,
    & (\tom^3)_{1} X^4 =\dfr{4}{3} X^4,
    \vsb\\
  (\tom^1)_{1} X^5=\dfr{2}{3} X^5, \q
    & (\tom^2)_{1} X^5=\dfr{1}{2} X^5, \q
    & (\tom^3)_{1} X^5=\dfr{5}{6} X^5, \q
    \vsb\\
  (X^1)_{1} X^1 = 8 X^2,
    & (X^1)_{1} X^2= 9 X^3,
    & (X^1)_{1} X^3= 8 X^4,
    \vsb\\
  (X^1)_{1} X^4 = 10 X^5,
    & (X^1)_{1} X^5 = 60\tom^1 + 72\tom^2 + 48\tom^3,
    & (X^2)_{1} X^2 = 12 X^4,
    \vsb\\
  (X^2)_{1} X^3 = 10 X^5,
    & (X^2)_{1} X^4 = 75 \tom^1 + 96 \tom^3,
    & (X^2)_{1} X^5 = 10 X^1,
    \vsb\\
  (X^3)_{1} X^3 = 80 \tom^2 + 96 \tom^3,
    & (X^3)_{1} X^4 = 10 X^1,
    & (X^3)_{1} X^5 = 8 X^2,
    \vsb\\
  (X^4)_{1} X^4 = 12 X^2,
    & (X^4)_{1} X^5 = 9 X^3,
    & (X^5)_{1} X^5 = 8 X^4,
    \vsb\\
  \la X^1,X^5 \ra = 36,
    & \la X^2, X^4 \ra = 45,
    & \la X^3, X^3 \ra = 40.
\end{array}
\end{equation}

\begin{lemma}
Let $w = a\tom^1 + b\tom^2 + c\tom^3 + d X^1 + e X^2 + f X^e + g
X^4 + h X^5$ be a conformal vectors in $\B$. Then $
(a,b,c,d,e,f,g,h)$ satisfies the following system of equations.
\begin{equation}\label{6ae}
\begin{array}{l}
  a^2 + 60dh + 75eg - a = 0,
  \vsb\\
  b^2 + 40f^2 + 72dh - b = 0,
  \vsb\\
  c^2 + 48f^2 + 48dh + 96eg - c = 0,
  \vsb\\
  4ad + 3bd + 5cd + 60eh + 60fg - 6d = 0,
  \vsb\\
  12d^2 + 18g^2 + 2ae + 4ce + 24fh - 3e = 0,
  \vsb\\
  bf + 3cf + 18de + 18gh - 2f = 0,
  \vsb\\
  18e^2 + 12h^2 + 2ag + 4cg + 24df - 3g = 0,
  \vsb\\
  4ah + 3bh + 5ch + 60dg + 60ef - 6h = 0.
\end{array}
\end{equation}
The central charge of $w$ is given by
\begin{equation}
\frac{4}{5} a^2 + \frac{1}{2}b^2 +\fr{5}{4} c^2
+80f^2+144dh+180eg.
\end{equation}
\end{lemma}

\medskip

Again we treat $a,b,c,d,e,f,g,h$ as variables and let $\I$ be the
ideal in $\C[a,b,c,d,e,f,g,h]$ generated by the eight polynomials
which appear on the left hand side of (\ref{6ae}). By using the
computer algebra system Risa/Asir, we found that the ideal $\I$ is
an intersection of $112$ prime ideals over the field $\Q$ of
rational numbers and that there are totally $256$ conformal
vectors in $\B$. The number of conformal vectors whose central
charges are rational numbers less than $1$ are listed in Table 2.

\begin{table}[h]
\begin{center}
\caption{Central charge (c.c.) and number of conformal vectors}
\medskip
\label{cclist}
\begin{tabular}{c|c|c|c|c|c|c|c|c}
\hline c.c. & $1/2$ &  $7/10$ & $4/5$ & $6/7$ & $25/28$ & $11/12$
& $14/15$ & $21/22$ \\
\hline number &\ $7$ &\ $9$  &\ $7$ &\ $14 $ &\ $5 $ &\ $6$ &\ $6$ &\ $6$\\
\hline
\end{tabular}
\end{center}
\end{table}

The following three lemmas were verified by computer.

\begin{lemma}\label{A-6A}
There are exactly seven conformal vectors of central charge $1/2$
in $\B$, namely, $\tilde{\omega}^2$ and $\sigma^j\hat{e}, 0 \le j
\le 5$. Note that
\[
\sigma^j\hat{e} = \fr{5}{32}\tom^1 + \fr{1}{8} \tom^2 +
\fr{1}{4}\tom^3 + \fr{1}{32}( \xi^j X^1 + \xi^{2j} X^2 + \xi^{3j}
X^3 + \xi^{4j} X^4 + \xi^{5j} X^5),
\]
where $\xi = e^{\pi \sqrt{-1}/3}$ is a primitive $6$-th root of
unity. Moreover, the inner product among $\sigma^j\hat{e}, 0 \le j
\le 5$ are
\begin{equation*}
\la \sigma^i\hat{e}, \sigma^j\hat{e}\ra =
\begin{cases}
5/2^{10} &\text{ if } i-j\equiv \pm 1 \mod 6,\\
13/2^{10} &\text{ if } i-j\equiv \pm 2 \mod 6,\\
1/32 &\text{ if } i-j\equiv   3 \mod 6.
\end{cases}
\end{equation*}
\end{lemma}
\medskip

We are mainly interested in mutually orthogonal conformal vectors
whose sum is the Virasoro element $\omega' = \tom^1+\tom^2+\tom^3$
of $U$.

\begin{lemma}
There are exactly $10$ triples $(u,v,w)$ of mutually orthogonal
conformal vectors in $\B$ such that the central charge of $u,v,w$
are $4/5, 6/7, 25/28$, respectively and $u+v+w=\omega'$.
\end{lemma}

\begin{lemma}
There are exactly $6$ triples $(u,v,w)$ of mutually orthogonal
conformal vectors in $\B$ such that the central charge of $u,v,w$
are $7/10, 11/12, 14/15$, respectively and $u+v+w=\omega'$.
\end{lemma}

\medskip \textbf{$4B$ case.} In this case, $L=L(6)\cong A_1\oplus
A_7$, $n_6=4$, and the conformal vectors $\tom^1\in
V_{\sqrt{2}A_1}$ and $\tom^2\in V_{\sqrt{2}A_7}$ are of central
charge $1/2$ and $7/5$ respectively. The product and the inner
product in $\B$ are given by \eqref{eq:GC} and
\begin{equation}\label{4b-eq}
\begin{array}{lll}
  (\tom^1)_1 X^1 = \dfr{1}{2} X^1,&
  (\tom^1)_1 X^2 = 0,&
  (\tom^1)_1 X^3 = \dfr{1}{2} X^3,
  \vsb\\
  (\tom^2)_1 X^1 = \dfr{3}{2} X^1,&
  (\tom^2)_1 X^2 = 2 X^2,&
  (\tom^2)_1 X^3 = \dfr{3}{2} X^3,
  \vsb\\
  (X^1)_1 X^1 = 12 X^2,&
  (X^1)_1 X^2 = 15 X^3,&
  (X^1)_1 X^3 = 112 \tom^1 + 120 \tom^2,
  \vsb\\
  (X^2)_1 X^2 = 200 \tom^2,\quad &
  (X^2)_1 X^3 = 15 X^1,\quad &
  (X^3)_1 X^3 = 12 X^2,
  \vsb\\
  \la X^1,X^3 \ra = 56,&
  \la X^2,X^2 \ra = 70.& \
\end{array}
\end{equation}

\begin{lemma}
Let $w=a\tom^1+b\tom^2 +cX^1 +dX^2+eX^3$ be a conformal vector in
$\B$. Then $(a,b,c,d,e)$ satisfies the following system of
equations.
\begin{equation}\label{4b}
\begin{array}{ll}
a^{2} + 112ce-a=0, \qquad &
b^{2} + 100d^{2} + 120ce-b=0, \\
ac + 3bc + 30de-2c=0, \qquad &
6c^{2} + 6e^{2} + 2bd-d=0, \\
ae + 3be + 30cd-2e=0. & \
\end{array}
\end{equation}
The central charge of $w$ is given by ${\frac {1}{2}} \,a^{2} + {
\frac {7 }{5}} \,b^{2} + 224\,c\,e + 140\,d^{2}$.
\end{lemma}

The solutions $(a,b,c,d,e)$ of Equation \eqref{4b} are as follows,
where $\xi=e^{\pi \sqrt{-1}/4}$ is a primitive $8$-th root of
unity.

\medskip
\noindent Central charge 1/2: $(1,0,0,0,0)$, \ $\big(\frac{1}{8},
\frac{5}{16}, \frac{1}{32}(\sqrt{-1})^j, \frac{1}{32}(-1)^j,
\frac{1}{32}(-\sqrt{-1})^j\big)$,\ $0 \le j \le 3$.

\medskip
\noindent Central charge 7/10: $(0,\frac{1}{2},0,\pm
\frac{1}{20},0)$,\ $\big(\frac{7}{8}, \frac{3}{16},
\frac{1}{32}(\sqrt{-1})^j, \frac{3}{160}(-1)^j,
\frac{1}{32}(-\sqrt{-1})^j\big)$,\ $0 \le j \le 3$.

\medskip
\noindent Central charge $21/22$: $\big(\frac{7}{11},
\frac{5}{11}, \frac{1}{22}\xi^{2j+1}, 0,
\frac{1}{22}\xi^{-(2j+1)}\big)$, \ $0 \le j \le 3$.

\medskip
\noindent Central charge $52/55$: $\big(\frac{4}{11},
\frac{6}{11}, \frac{1}{22}\xi^{2j+1}, 0,
\frac{1}{22}\xi^{-(2j+1)}\big)$, \ $0 \le j \le 3$.

\medskip
\noindent Central charge $6/5$: $(1,\frac{1}{2},0,\pm
\frac{1}{20},0)$,\ $\big(\frac{1}{8}, \frac{13}{16},
\frac{1}{32}(\sqrt{-1})^j, -\frac{3}{160}(-1)^j,
\frac{1}{32}(-\sqrt{-1})^j\big)$,\ $0 \le j \le 3$.

\medskip
\noindent Central charge $7/5$: $(0,1,0,0,0)$,\ $\big(\frac{7}{8},
\frac{11}{16}, \frac{1}{32}(\sqrt{-1})^j, -\frac{1}{32}(-1)^j,
\frac{1}{32}(-\sqrt{-1})^j\big)$,\ $0 \le j \le 3$.

\medskip
\noindent Central charge $19/10$:  $(1,1,0,0,0)$.

\bigskip
\textbf{$2B$ case.} In this case, $L=L(7)\cong D_8$, $n_7=2$, and
the conformal vector $\tom^1 \in V_{\sqrt{2}D_8}$ is of central
charge $1$. The product and the inner product in $\B$ are given by
\eqref{eq:GC} and
\[
(\tom^1)_1 X^1=2X^1, \quad (X^1)_1X^1= 512 \tom^1,\quad \la
X^1,X^1\ra =128.
\]

\begin{lemma}
Let $w=a\tom^1+bX^1$ be a conformal vector in $\B$. Then we have
$a=1/2$ and $b=\pm 1/32$, and the central charge of $w$ is $1/2$.
In other words, there are exactly two conformal vectors in $\B$
and both of them are of central charge $1/2$.
\end{lemma}

\medskip \textbf{$3C$ case.} In this case, $L=L(8)\cong A_8$,
$n_8=3$, and the conformal vector $\tilde{\omega}^1 \in
V_{\sqrt{2}A_8}$ is of central charge $16/11$. The product and the
inner product in $\B$ are given by \eqref{eq:GC} and
\begin{equation}\label{3c-eq}
\begin{array}{lll}
(\tilde{\omega}^1)_1 X^1 = 2X^1, & \quad (\tilde{\omega}^1)_1 X^2
= 2X^2, & \
\vsb\\
(X^1)_1X^1= 20 X^2, & \quad (X^2)_1X^2= 20 X^1, & \quad
(X^1)_1X^2=231\tilde{\omega}^1,
\vsb\\
\la X^1, X^2\ra=84. & \  & \
\end{array}
\end{equation}

\begin{lemma}\label{A-3C}
Let $w=a\tom^1 +bX^1 +cX^2$ be a conformal vector in $\B$. Then
$(a,b,c)$ satisfies the following system of equations.
\begin{equation}\label{3c}
\begin{array}{ll}
a^2+924bc-a=0, & \qquad
2ab+10c^2-b=0,\\
2ac + 10b^2-c=0. & \
\end{array}
\end{equation}
The central charge of $w$ is given by $\frac{16}{11}a^2+336 bc$.
\end{lemma}

The solutions $(a,b,c)$ of Equation \eqref{3c} are as follows,
where $\xi=e^{2\pi \sqrt{-1}/3}$ is a primitive cubic root of
unity.

\medskip
\noindent  Central charge $1/2$: $\big(\frac{11}{32},
\frac{1}{32}\xi^j, \frac{1}{32}\xi^{2j}\big)$,\ $j=0,1,2$.

\medskip
\noindent Central charge $21/22$: $\big(\frac{21}{32},
-\frac{1}{32}\xi^j, -\frac{1}{32}\xi^{2j}\big)$,\ $j=0,1,2$.

\medskip
\noindent Central charge $16/11$: $(1,0,0,0)$.

\section{Characters of $W$-algebras and the structure of $U$ for
the cases of  $3C$ and $5A$}

In this appendix, we shall determine the structure of the coset
subalgebra $U$ defined by \eqref{eq:UDEF} for the cases of $3C$
and $5A$. The main tool is the character of parafermion algebras.
First, let us recall a construction of parafermion algebras from
\cite{dl}.

Let $\Lambda_0$ and $\Lambda_1$ be the fundamental weights of the
affine Lie algebra $\hat{sl}_2(\C)$. For any positive integer
$\ell$ and $0\leq j\leq \ell$, let $\L(\ell,j)$ be the irreducible
highest weight module of $\hat{sl}_2(\C)$ with the highest weight
$(\ell-j)\Lambda_0 +j\Lambda_1$. Note that $\L(\ell,0)$ has a
natural VOA structure and $\{\L(\ell,j)|\, 0\leq j\leq \ell\}$ is
the set of all inequivalent irreducible modules of $\L(\ell,0)$
(cf. \cite{fz}).

Now let ${A_1}^{\ell} = \Z \epsilon_1\oplus \cdots \oplus \Z
\epsilon_\ell$ be an even lattice with $\la
\epsilon_i,\epsilon_j\ra =2\delta_{i,j}$ and $V_{{A_1}^{\ell}}$
the lattice VOA associated with ${A_1}^{\ell}$. Then
$V_{{A_1}^{\ell}}\cong (V_{A_1})^{\tensor \ell}\cong
\L(1,0)^{\tensor \ell}$.

Set $H^{(\ell)}=\epsilon_1(-1)\cdot 1+\cdots +
\epsilon_\ell(-1)\cdot 1$, $E^{(\ell)}=e^{\epsilon_1}+\cdots
+e^{\epsilon_\ell}$, and $F^{(\ell)}=e^{-\epsilon_1}+\cdots
+e^{-\epsilon_\ell}$. Then $\C H^{(\ell)}+\C E^{(\ell)}+\C
F^{(\ell)}$ forms a simple Lie algebra $sl_2(\C)$ inside the
weight one subspace of $V_{{A_1}^{\ell}}$. Moreover, the subVOA
generated by $\{ H^{(\ell)}, E^{(\ell)}, F^{(\ell)} \}$ is
isomorphic to $\L(\ell,0)$ (cf. \cite{dl}).

Let $\gamma= \epsilon_1+\cdots+\epsilon_\ell$. Then
$\gamma(-1)\cdot 1= H^{(\ell)}$ and it is easy to verify that
\[
e^{\gamma}=\frac{1}{\ell!}((E^{(\ell)})_{-1})^{\ell-1}E^{(\ell)}.
\]
Thus $\L(\ell,0)$ contains a subalgebra isomorphic to the lattice
VOA $V_{\Z \gamma}$.

Let $A= \la e^\gamma \ra$ be the subgroup of $\C\{{A_1}^{\ell}\}$
generated by $e^\gamma$ and denote by $\Omega_{\L(\ell,j)}^A$ the
set of all highest weight vectors for $V_{\Z\gamma}$ in
$\L(\ell,j)$, $0\leq j\leq \ell$.  It is shown in \cite{dl} that
$\Omega_{\L(\ell,0)}^A= \oplus_{k=0}^{\ell-1}
(\Omega_{\L(\ell,0)}^A)^k$ is a generalized VOA and
$\Omega_{\L(\ell,j)}^A$ are irreducible
$\Omega_{\L(\ell,0)}^A$-modules. Note that
$(\Omega_{\L(\ell,0)}^A)^0$ itself is a VOA, which we shall denote
by $W_\ell$.

Now let
\[
  \L(\ell,j)= \bigoplus_{k=0}^{2\ell -1} V_{(k/2\ell)\gamma +
  \Z\gamma}\tensor W_\ell(j,k)
\]
be the decomposition of $\L(\ell,j)$ as a $V_{\Z\gamma}\otimes
W_\ell$-module, where $W_{\ell}(j,k)$ is the multiplicity of
$V_{k\gamma/2\ell + \Z\gamma}$ in $\L(\ell,j)$. By \cite{dl},
$W_\ell(j,k)=0$ if $j+k \equiv 1 \mod 2$ and so
\begin{equation}\label{eq:LDEC}
  \L(\ell,j)=
\begin{cases}
\displaystyle
  \bigoplus_{k=0}^{\ell -1} V_{(k/\ell)\gamma + \Z\gamma} \tensor
  W_\ell(j,2k)& \text{ if }j \text{ is
  even},\\
\displaystyle \bigoplus_{k=0}^{\ell -1} V_{((2k+1)/2\ell)\gamma +
\Z\gamma}\tensor W_\ell(j,2k+1)& \text{ if }j \text{ is
  odd}.
\end{cases}
\end{equation}

\begin{proposition}[cf. Dong and Lepowsky \cite{dl}]
  All $W_\ell(j,k)$, $0\leq j\leq \ell$, $0\leq k \leq
  2\ell-1$, $j\equiv k \mod 2$, are irreducible
  $W_\ell$-modules.
\end{proposition}

By \eqref{eq:LDEC}, we can actually compute the character of
$W_\ell(j,k)$ by using the character of $\L(\ell,j)$ and
$V_{(k/2\ell)\gamma + \Z\gamma}$. Recall that the character $\ch
_M(q)$ of a module $M$ of a VOA $V$ is defined by $\ch_M(q) =
\tr_M q^{L(0)}$, where $L(0) = \omega_1$ with $\omega$ being the
Virasoro element of $V$. In Kac and Raina \cite{kr}, the following
formula is proved.
\begin{equation}\label{eq:tr1}
  \tr_{\L(\ell,j)} \sigma (-\tfr{1}{2}z \gamma)
  q^{L(0)-\ell/8(\ell+2)}
  = \dfr{\theta_{j+1,\ell+2}(\tau,z)
    -\theta_{-j-1,\ell+2}(\tau,z)}
    {\theta_{1,2}(\tau,z)-\theta_{-1,2}(\tau,z)},
\end{equation}
where $q = e^{2\pi\sqrt{-1}\tau}$, $\sigma (z \gamma) =
e^{2\pi\sqrt{-1} z \gamma(0)}$ for $z\in \Q$, and
\begin{equation*}
  \theta_{n,m}(\tau,z) = \dsum_{j\in n/2m +\Z}
  e^{2\pi\sqrt{-1} mzj} q^{mj^2}.
\end{equation*}

By using the $q$-integers $[n]_a = (a^n-a^{-n})/(a-a^{-1})
=a^{n-1}+a^{n-3}+\cds +a^{-n+3}+a^{-n+1}$ for $n\in \Z$, we can
rewrite the formula \eqref{eq:tr1} in the following form.
\begin{equation}\label{eq:tr2}
\tr_{\L(\ell,0)}\sigma(-\tfr{1}{2}z \gamma) q^{L(0)}
=\dfr{1+\dsum_{m \in \Z \setminus \{ 0\}}
      [2(\ell+2)m+1]_{e^{\pi\sqrt{-1} z}} \cdot
      q^{(\ell+2)m^2+m}}{1+\dsum_{n \in \Z
      \setminus \{ 0\}}
      [4n+1]_{e^{\pi\sqrt{-1} z}} \cdot q^{2n^2+n}} .
\end{equation}

The character of $V_{(k/2\ell)\gamma +\Z\gamma}$ is also well
known (cf. \cite{flm}). It is given by
\begin{equation}\label{eq:CHV}
\ch_{V_{(k/2\ell)\gamma +\Z\gamma}} (q)= \tr_{V_{(k/2\ell)\gamma
+\Z\gamma}}q^{L(0)}= \frac{1}{\prod_{n=1}^{\infty} (1-q^n)}
\theta_{(k/2\ell)\gamma +\Z\gamma}(q),
\end{equation}
where $\theta_{(k/2\ell)\gamma +\Z\gamma}(q)=\sum_{\al\in
(k/2\ell)\gamma +\Z\gamma} q^{\la \al,\al\ra /2}$.

\medskip
Next, we consider another construction of the VOA
$W_\ell=W_\ell(0,0)$ given in \cite{ly3} by using the lattice VOA
$V_{\sqrt{2}A_{\ell-1}}$, namely,
\[
W_\ell\cong \{ u\in V_{\sqrt{2}A_{\ell-1}}|\ s(A_{\ell -1})_1 u=0
\},
\]
where
\[
s(A_{\ell -1})=\frac{1}{2(\ell+2)}\sum_{\al\in \Phi^+(A_{\ell
-1})}\left( \al(-1)^2\cdot 1 -2(e^{\sqrt{2}\al}+
e^{-\sqrt{2}\al})\right)
\]
and $\Phi^+(A_{\ell -1})$ is the set of all positive roots of the
lattice of type $A_{\ell -1}$ (cf. \eqref{cv1}).
\vsb

Let $N=\mathrm{span}_\Z \{ \epsilon_1-\epsilon_2,
\epsilon_2-\epsilon_3, \ldots, \epsilon_{\ell-1}-\epsilon_\ell\}
\subset {A_1}^{\ell}$. Then $N \cong \sqrt{2}A_{\ell-1}$.
Moreover,
\[
{A_1}^{\ell} = \bigcup_{k=0}^{\ell-1} (k \epsilon_\ell + N +
\Z\gamma) = \bigcup_{k=0}^{\ell-1} \Big( (k \eta +N)+
(\frac{k}{\ell}\gamma +\Z\gamma)\Big),
\]
where $\gamma= \epsilon_1+\cdots +\epsilon_\ell$ and $\eta=
\frac{1}{\ell} (-\epsilon_1-\cdots -\epsilon_{\ell-1}
+(\ell-1)\epsilon_\ell)$. Indeed, $|{A_1}^\ell/(N + \Z\gamma)| =
\ell$ and $\eta + \gamma/\ell = \epsilon_\ell$. Note that $\la
\eta,\gamma \ra = 0$ and $\la N,\gamma\ra=0$.
\vsb

By the above argument, we obtain the following proposition.

\begin{proposition}\label{propB-2}
For any positive integer $\ell$ and $k=0, \dots,\ell-1$,
\[
W_\ell(0, 2k) \cong \{ u\in V_{k\eta +N}\,|\, s(A_{\ell-1})_1 u=0
\}.
\]
\end{proposition}

\subsection{The case for $3C$}

We consider an embedding of a $\sqrt{2}E_8$ lattice into $\R^9$ by
using $N$ and $\eta$ for the case $\ell = 9$. Indeed,
$\sqrt{2}\alpha_j$, $0 \le j \le 8$, where $\alpha_0,\alpha_1,
\ldots,\alpha_8$ are the nine nodes in the extended $E_8$ diagram
\eqref{ext}, can be realized as $\sqrt{2}\alpha_j =
\epsilon_{j+1}-\epsilon_{j+2}$, $0 \le j \le 7$, and
$\sqrt{2}\alpha_8 = 3\eta + (\epsilon_7 - \epsilon_8) +
2(\epsilon_8 - \epsilon_9) \in 3\eta + N$. Then $N=\sqrt{2}L(8)$,
where $L(8)$ is the sublattice spanned by
$\alpha_0,\alpha_1,\ldots,\alpha_7$ (cf. Section 3). Hence
$\sqrt{2}E_8 = N \cup (3\eta+N) \cup (6\eta+N)$. This implies that
\[
  V_{\sqrt{2}E_8}= V_{\sqrt{2}A_8} \oplus V_{3\eta + \sqrt{2}A_8}
  \oplus V_{6\eta +\sqrt{2}A_8}.
\]
Furthermore, by the definition \eqref{eq:UDEF} of the coset
subalgebra, $U = \{ u \in V_{\sqrt{2}E_8}\,|\,s(A_8)_1 u=0\}$.
So by Proposition \ref{propB-2}, we have the following result.

\begin{proposition}
  As a module of $W_9$,
  \begin{equation}\label{eq:3CW}
    U \cong W_9(0,0)\oplus W_9(0,6)\oplus W_9(0,12).
  \end{equation}
\end{proposition}

\begin{remark}
  It is easy to see that $W_9(0,6)$ and $W_9(0,12)$ have the same
  character by their construction.
  In fact, $W_9(0,6)$ is the dual module of $W_9(0,12)$.
  Thus by \eqref{eq:3CW}, $\mathrm{ch}_{U}(q) \equiv \mathrm{ch}_{W_9}(q) \mod 2$.
\end{remark}

As discussed in Subsection \ref{sub5-3c}, the coset subalgebra $U$
can be decomposed into a direct sum of irreducible modules of
$L(1/2,0) \otimes L(21/22,0)$ in the following form.
\begin{equation}\label{eq:3CUbis}
\begin{split}
   U
   & \cong  m_1 L(\frac{1}2,0)\tensor L(\frac{21}{22},0)
   \oplus m_2 L(\frac{1}{2},0)\tensor L(\frac{21}{22},8)
   \\
   & \quad \oplus m_3 L(\frac{1}2,\frac{1}2)\tensor  L(\frac{21}{22},\frac{7}{2})
     \oplus m_4 L(\frac{1}2,\frac{1}2)\tensor
     L(\frac{21}{22},\frac{45}{2})
   \\
   & \quad \oplus m_5 L(\frac{1}2,\frac1{16})\tensor L(\frac{21}{22},\frac{31}{16})
      \oplus m_6 L(\frac{1}2,\frac{1}{16})\tensor L(\frac{21}{22},\frac{175}{16}),
\end{split}
\end{equation}
where $m_j \in \Z$ denotes the multiplicity of each summand.
We want to show that $m_j=1$ for $j=1,\ldots,6$.

By \eqref{eq:LDEC},
$\L(9,0) = \oplus_{k=0}^8 V_{(k/9)\gamma + \Z\gamma} \otimes W_9(0,2k)$.
Let $z=j/9$ in \eqref{eq:tr1}.
Then $\sigma (-\frac{j}{18}\gamma)$ acts on $V_{(k/9)\gamma +\Z\gamma}$
as a scalar $e^{-2\pi\sqrt{-1}kj/9}$ and acts on $W_9(0,2k)$ as the identity.
Hence
\begin{equation*}
  \tr_{\L(9,0)} \sigma (-\frac{j}{18}\gamma) q^{L(0)} = \sum_{k=0}^8
  e^{-2\pi\sqrt{-1}kj/9} \ch_{V_{(k/9)\gamma +\Z\gamma}}(q)
  \ch_{W_9(0,2k)}(q)
\end{equation*}
and thus
\begin{equation*}
  \ch_{V_{(k/9)\gamma +\Z\gamma}}(q)
  \ch_{W_9(0,2k)}(q) = \frac{1}{9} \sum_{j=0}^8 \tr_{\L(9,0)}
  \sigma(-\frac{j}{18}\gamma) q^{L(0)} e^{2\pi\sqrt{-1}kj/9}.
\end{equation*}

Now, consider the case for $k=0$, $3$, and $6$. By \eqref{eq:tr1}
and \eqref{eq:CHV}, we can calculate that
$$
\begin{array}{l}
  \begin{array}{ll}
    \ch_{W_9(0,0)}(q)
    & = 1+ q^2+ 2q^3+ 4q^4 + 6q^5 +11q^6+16q^7+27q^8
    \vsb\\
    & \quad +40q^9 +62q^{10}+90q^{11}+137q^{12} +194q^{13}+284q^{14}
    \vsb\\
    & \quad + 400q^{15} + 569q^{16}+788q^{17} + 1102q^{18} +1504q^{19}
    \vsb\\
    & \quad +2066q^{20} + 2792q^{21} +3776q^{22}+5046q^{23}+\cdots ,
  \end{array}
  \vsb\\
  \begin{array}{ll}
    \ch_{W_9(0,6)}(q)=\ch_{W_9(0,12)}(q)
    & = q^2+ q^3+3q^4+5q^5+9q^6+ 14q^7
    \vsb\\
    & \quad + 25q^8 + 36q^9+58q^{10} + 86q^{11}+\cdots .
  \end{array}
\end{array}
$$
Let $[h_1,h_2]_n$ be the weight $n$ subspace of $L(1/2,h_1)
\otimes L(21/22,h_2)$. Its dimension for small $n$ is as follows.

\begin{align*}
  & \dim [0,0]_2=2,
  && \dim [0,0]_4=5,
  && \dim [\frac{1}{16},\frac{31}{16}]_4=4,
  \\
  & \dim [0,0]_8=27,
  && \dim [\frac{1}{16},\frac{31}{16}]_8=36,
  && \dim [\frac{1}{2},\frac{7}{2}]_8=13,
  \\
  & \dim [0,0]_{11}=75,
  && \dim [\frac{1}{16},\frac{31}{16}]_{11}=130,
  && \dim [\frac{1}{2},\frac{7}{2}]_{11}=51,
  \\
  & \dim [0,8]_{11}=5,
  && \dim [0,0]_{23}=3073,
  && \dim [\frac{1}{16},\frac{31}{16}]_{23}=7040,
  \\
  &\dim [\frac{1}{2},\frac{7}{2}]_{23}=3510,
  && \dim [0,8]_{23}=946,
  && \dim [\frac{1}{16},\frac{175}{16}]_{23}=490.
\end{align*}

Comparing the characters of \eqref{eq:3CW} and \eqref{eq:3CUbis},
we have $m_1=m_2=m_3=m_5=m_6=1$ and $m_4 \equiv 1\mod 2$.
Moreover, $m_4\leq 1$ since $L(1/2,1/2) \otimes L(21/22,45/2)$ is
a simple current module. Hence we obtain the following theorem.

\begin{theorem}\label{B-3C}
As a module of $L(1/2,0) \otimes L(21/22,0)$,
\begin{equation*}
\begin{split}
 U & \cong  L(\frac{1}2,0)\tensor L(\frac{21}{22},0)
 \oplus L(\frac{1}{2},0)\tensor L(\frac{21}{22},8)\\
 & \quad \oplus  L(\frac{1}2,\frac{1}2)\tensor L(\frac{21}{22},\frac{7}{2})
 \oplus L(\frac{1}2,\frac{1}2)\tensor
 L(\frac{21}{22},\frac{45}{2})\\
  & \quad  \oplus L(\frac{1}{2},\frac{1}{16})\tensor L(\frac{21}{22},\frac{31}{16})
  \oplus L(\frac{1}{2},\frac{1}{16})\tensor
  L(\frac{21}{22},\frac{175}{16}).
\end{split}
\end{equation*}
\end{theorem}

\subsection{The case for $5A$}

In this case we consider an embedding of a $\sqrt{2}E_8$ lattice
into $\R^{10}$ by using two $N$'s and two $\eta$'s for the case
$\ell = 5$. Let $\epsilon_i, \epsilon'_i \in \R^{10}$, $1 \le i
\le 5$, be such that $\la \epsilon_i,\epsilon_j \ra = \la
\epsilon'_i, \epsilon'_j \ra = 2\delta_{ij}$ and $\la
\epsilon_i,\epsilon'_j \ra = 0$. Set
\begin{equation*}
N = \mathrm{span}_\Z \{ \epsilon_i - \epsilon_{i+1};\ 1 \le i \le
4\},
\qquad N' = \mathrm{span}_\Z \{ \epsilon'_i -
\epsilon'_{i+1};\ 1 \le i \le 4\},
\end{equation*}
\begin{equation*}
\eta = \frac{1}{5}(-\epsilon_1 - \epsilon_2 - \epsilon_3 -
\epsilon_4 + 4\epsilon_5),
\qquad \eta' = \frac{1}{5}(-\epsilon'_1
- \epsilon'_2 - \epsilon'_3 - \epsilon'_4 + 4\epsilon'_5).
\end{equation*}

Let $\sqrt{2}\alpha_j = \epsilon_{j+1} - \epsilon_{j+2}$, $0 \le j
\le 3$, $\sqrt{2}\alpha_7 = \epsilon'_1 - \epsilon'_2$,
$\sqrt{2}\alpha_6 = \epsilon'_2 - \epsilon'_3$, $\sqrt{2}\alpha_5
= \epsilon'_3 - \epsilon'_4$, $\sqrt{2}\alpha_8 = \epsilon'_4 -
\epsilon'_5$, $\sqrt{2}\alpha_4 = \eta + 2\eta' + (\epsilon'_4 -
\epsilon'_5) \in \eta + 2\eta' + N'$. Then
$\alpha_0,\alpha_1,\ldots,\alpha_8$ are the nine nodes in the
extended $E_8$ diagram \eqref{ext} and $N \oplus N' =
\sqrt{2}L(4)$, where $L(4)$ is the sublattice spanned by
$\alpha_j$,$0 \le j \le 8$, $j \ne 4$ (cf. Section 3). Hence
$|\sqrt{2}E_8/N\oplus N'| = 5$ and
\begin{align*}
\sqrt{2}E_8 &= N\oplus N' \cup (\eta+N)\oplus(2\eta'+N') \cup
(2\eta+N)\oplus(4\eta'+N')\\
& \qquad \cup (3\eta+N)\oplus(\eta'+N') \cup
(4\eta+N)\oplus(3\eta'+N').
\end{align*}

This implies that
\begin{align*}
  V_{\sqrt{2}E_8} &= (V_{\sqrt{2}A_4}\otimes V_{\sqrt{2}A_4})\oplus
  (V_{\eta+\sqrt{2}A_4}\otimes V_{2\eta+\sqrt{2}A_4}) \oplus
  (V_{2\eta+\sqrt{2}A_4}\otimes V_{4\eta+\sqrt{2}A_4})
  \\
  &\qquad \oplus (V_{3\eta+\sqrt{2}A_4}\otimes
  V_{\eta+\sqrt{2}A_4})\oplus (V_{4\eta+\sqrt{2}A_4}\otimes
  V_{3\eta+\sqrt{2}A_4}).
\end{align*}

Furthermore, the following proposition holds by the definition
\eqref{eq:UDEF} of the coset subalgebra $U$ and Proposition
\ref{propB-2}.
\begin{proposition}
  As a module of $W_5 \otimes W_5$,
  \begin{equation}\label{eq:5AW}
  \begin{split}
    U
    & \cong (W_5(0,0)\otimes W_5(0,0))\oplus (W_5(0,2)\otimes
      W_5(0,4)) \oplus (W_5(0,4)\otimes W_5(0,8))
      \\
    & \qquad  \oplus ( W_5(0,6)\otimes W_5(0,2))
      \oplus (W_5(0,8)\otimes W_5(0,6)).
  \end{split}
  \end{equation}
\end{proposition}

Note that $\mathrm{ch}_{V_{\eta+\sqrt{2}A_4}}(q)=
\mathrm{ch}_{V_{4 \eta+\sqrt{2}A_4}}(q)$ and
$\mathrm{ch}_{V_{2\eta+\sqrt{2}A_4}}(q)= \mathrm{ch}_{V_{3
\eta+\sqrt{2}A_4}}(q)$. Thus by the construction, it is also easy
to see that
\[
  \mathrm{ch}_{W_5(0,2)\tensor W_5(0,4)}(q) =
  \mathrm{ch}_{W_5(0,4)\tensor W_5(0,8)}(q) =
  \mathrm{ch}_{ W_5(0,6)\tensor W_5(0,2)}(q) =
  \mathrm{ch}_{W_5(0,8)\tensor W_5(0,6)}(q).
\]
Hence
\begin{equation}\label{eq:5AWbis}
  \ch_U(q)=\ch_{W_5(0,0)}(q)^2
  +4 \ch_{W_5(0,2)}(q) \ch_{W_5(0,4)}(q).
\end{equation}

By \eqref{eq:LDEC}, $\L(5,0) = \oplus_{k=0}^4
V_{(k/5)\gamma+\Z\gamma} \otimes W_5(0,2k)$. Let $z=j/5$ in
\eqref{eq:tr1}. Then $\sigma (-\frac{j}{10}\gamma)$ acts on
$V_{(k/5)\gamma +\Z\gamma}$ as a scalar $e^{-2\pi\sqrt{-1}kj/5}$
and acts on $W_5(0,2k)$ as the identity. Then arguing as in the
case for $3C$, we get
\begin{equation*}
  \ch_{V_{(k/5)\gamma+\Z\gamma}}(q) \ch_{W_5(0,2k)}(q)
= \frac{1}{5} \sum_{j=0}^4 e^{2\pi\sqrt{-1}kj/5} \tr_{\L(5,0)}
    \sigma(-\frac{j}{10} \gamma) q^{L(0)}.
\end{equation*}

Now using \eqref{eq:tr1} and \eqref{eq:CHV}, we can calculate that
\begin{equation*}
\begin{array}{ll}
  \ch_{W_5(0,0)}(\tau)=
  & 1+q^2+2q^3+4q^4+6q^5+10q^6+14q^7+23q^8+32q^9+48q^{10}
  \vsb\\
  & +66q^{11}+96q^{12}+130q^{13}+183q^{14}+246q^{15}+\cds,
  \vsb\\
  \ch_{W_5(0,2)}(\tau)=
  & q^{4/5}(1+q+2q^2+3q^3+6q^4+8q^5+14q^6+20q^7+31q^8+43q^9
  \vsb\\
  & +64q^{10}+87q^{11}+125q^{12}+169q^{13}+234q^{14}+313q^{15}
    +\cds ),
  \vsb\\
  \ch_{W_5(0,4)}(\tau)=
  & q^{6/5}(1+q+3q^2+4q^3+7q^4+10q^5+17q^6+23q^7+36q^8+50q^9
  \vsb\\
  & +73q^{10}+100q^{11}+142q^{12}+191q^{13}+265q^{14}+
    353q^{15}+\cds ).
\end{array}
\end{equation*}

In Subsection \ref{sub5-5A} we have shown that $U$ contains a
subalgebra isomorphic to $L(1/2,0) \otimes L(25/28,0) \otimes
L(25/28,0)$. We also know all irreducible modules $L(1/2,h_1)
\otimes L(25/28,h_2) \otimes L(25/28,h_3)$ with integral weights.
Note that $U$ is a direct sum of these irreducible modules.
Comparing the characters of those irreducible modules and
\eqref{eq:5AWbis}, we can verify the following theorem.

\begin{theorem}\label{B-5A}
As a module of $L(1/2,0) \otimes L(25/28,0) \otimes L(25/28,0)$,
\[
  \begin{array}{l}
    U \cong L(\frac{1}{2},0)\tensor L(\frac{25}{28},0)\tensor L(\frac{25}{28},0)
      \oplus L(\frac{1}{2},\frac{1}{16})\tensor
        L(\frac{25}{28},\frac{5}{32})\tensor
        L(\frac{25}{28},\frac{57}{32})
    \vsb\\
      \qquad \oplus L(\frac{1}{2},\frac{1}{16})\tensor
        L(\frac{25}{28},\frac{57}{32})\tensor
        L(\frac{25}{28},\frac{5}{32})
      \oplus L(\frac{1}{2},\frac{1}{2})\tensor
        L(\frac{25}{28},\frac{3}{4})\tensor L(\frac{25}{28},\frac{3}{4})
    \vsb\\
      \qquad \oplus L(\frac{1}{2},0)\tensor L(\frac{25}{28},\frac{3}{4})\tensor
        L(\frac{25}{28},\frac{13}{4})
      \oplus L(\frac{1}{2},0)\tensor L(\frac{25}{28},\frac{13}{4})\tensor
        L(\frac{25}{28},\frac{3}{4})
    \vsb\\
      \qquad \oplus L(\frac{1}{2},\frac{1}{16})\tensor
        L(\frac{25}{28},\frac{57}{32})\tensor
        L(\frac{25}{28},\frac{165}{32})
      \oplus L(\frac{1}{2},\frac{1}{16})\tensor
        L(\frac{25}{28},\frac{165}{32})\tensor
        L(\frac{25}{28},\frac{57}{32})
    \vsb\\
      \qquad \oplus L(\frac{1}{2},\frac{1}{2})\tensor
        L(\frac{25}{28},\frac{13}{4})\tensor L(\frac{25}{28},\frac{13}{4})
      \oplus L(\frac{1}{2},\frac{1}{2})\tensor L(\frac{25}{28},0)\tensor
        L(\frac{25}{28},\frac{15}{2})
    \vsb\\
      \qquad \oplus L(\frac{1}{2},\frac{1}{2})\tensor
        L(\frac{25}{28},\frac{15}{2})\tensor L(\frac{25}{28},0)
      \oplus L(\frac{1}{2},0)\tensor L(\frac{25}{28},\frac{15}{2})\tensor
        L(\frac{25}{28},\frac{15}{2}).
  \end{array}
\]
\end{theorem}

\subsection{Character of the 6A-algebra}

Let $U_{\mathrm{6A}}$ be the 6A-algebra, the coset subalgebra of
$V_{\sqrt{2} E_8}$ constructed from the coset decomposition
$E_8/(A_5\oplus A_2\oplus A_1)$. We show here that
\begin{equation}\label{1}
\begin{array}{l}
  U_{\mathrm{6A}}
  \simeq U_{\mathrm{3A}}(0)\tensor L(25/28,0)
  \bigoplus U_{\mathrm{3A}}(5/7)\tensor L(25/28,9/7)
  \vsb\\
  \hspace{2cm}\bigoplus U_{\mathrm{3A}}(1/7)\tensor L(25/28,34/7),
\end{array}
\end{equation}
where $U_{\mathrm{3A}}(0)$ denotes the 3A-algebra and
$U_{\mathrm{3A}}(h)$, $h=5/7, 1/7$, denote the irreducible
$U_{\mathrm{3A}}(0)$-modules whose top weights are equal to $h$. We
show the above isomorphism by computing the $q$-character of
$U_{\mathrm{6A}}$. Recall the notion of $W_\ell$-algebras in
Appendix B. By considering an isometric embedding $\sqrt{2}
A_5\oplus \sqrt{2} A_2 \oplus \sqrt{2} A_1 \hookrightarrow
A_1^6\oplus A_1^3\oplus \sqrt{2}A_1$, one can easily verify that
$$
\begin{array}{l}
  U_{\mathrm{6A}}
  \simeq W_6(0,0)\tensor W_3(0,0)\tensor L(1/2,0)
  \bigoplus W_6(0,2)\tensor W_3(0,4)\tensor L(1/2,1/2)
  \vsb\\
  \bigoplus W_6(0,4)\tensor W_3(0,2)\tensor L(1/2,0)
  \bigoplus W_6(0,6)\tensor W_3(0,0)\tensor L(1/2,1/2)
  \vsb\\
  \bigoplus W_6(0,8)\tensor W_3(0,4)\tensor L(1/2,0)
  \bigoplus W_6(0,10)\tensor W_3(0,2)\tensor L(1/2,1/2),
\end{array}
$$
where we have used the fact that the lattice VOA $V_{\sqrt{2} A_1}$
is isomorphic to the code VOA associated to a binary code $\{
(00),(11)\}$. Since $W_3(0,0)\simeq L(4/5,0)\oplus L(4/5,3)$ and
$W_3(0,2)\simeq W_3(0,4)\simeq L(4/5,2/3)$ as $L(4/5,0)$-modules, we
only need to compute the characters of $W_6(0,2s)$, $0\leq s\leq 5$.
However, by noticing the dual module relations, we know that
$\ch_{W_6(0,2)}(q)=\ch_{W_6(0,10)}(q)$ and
$\ch_{W_6(0,4)}(q)=\ch_{W_6(0,8)}(q)$. Therefore, we should compute
characters of $W_6(0,2s)$ for $s=0,1,2,3$. By a method in the
preceding, we can obtain the following results.
$$
\begin{array}{l}
  \ch_{W_6(0,0)}=1 + q^2 + 2 q^3 + 4 q^4 + 6 q^5 + 11 q^6+\cds ,
  \vsb\\
  \ch_{W_6(0,2)}=q^{5/6} + q^{11/6} + 2 q^{17/6} + 3 q^{23/6} +
    6 q^{29/6} + 9 q^{35/6}+\cds,
  \vsb\\
  \ch_{W_6(0,4)}=q^{4/3} + q^{7/3} + 3 q^{10/3} + 4 q^{13/3} +
    8 q^{16/3}+\cds ,
  \vsb\\
  \ch_{W_6(0,6)}=q^{3/2} + q^{5/2} + 3 q^{7/2} + 5 q^{9/2}
    + 8 q^{11/2}+\cds .
\end{array}
$$
Then by comparing characters, we can establish the desired
isomorphism \eqref{1}.

\subsubsection{Highest weight vector of weight $(0,1/7,34/7)$}

In this section, we prove that the 6A-algebra is generated by its
weight two subspace as a vertex operator algebra. We use the
notation as in section 5.6 of our preprint. Let
$w^1=\tilde{\omega}^1=\tilde{\omega}(A_2)$,
$w^2=\tilde{\omega}(E_6)$ and
$w^3=\tilde{\omega}(A_1)+\tilde{\omega}(A_5)-w^2$. Then
$w^1,w^2,w^3$ are mutually orthogonal conformal vector of central
charges 4/5, 6/7 and 25/28, respectively, and the sum $w^1+w^2+w^3$
is the Virasoro vector of $U_{\mathrm{6A}}$ and the sum $w^1+w^2$ is
the Virasoro vector of $U_{\mathrm{3A}}(0)$ in the decomposition
\eqref{1}. We can write down $w^i$ explicitly as follows.
$$
  w^1=\tilde{\omega}(A_2),\quad
  w^2=\fr{2}{7}\tilde{\omega}(A_1)+\fr{4}{7}\tilde{\omega}(A_5)
      +\fr{1}{14}X^3,\quad
  w^3=\fr{5}{7}\tilde{\omega}(A_1)+\fr{3}{7}\tilde{\omega}(A_5)
      -\fr{1}{14}X^3.
$$
%By the isomorphism \eqref{1}, we know that any irreducible
%$\vir (w^1)\tensor \vir (w^2)\tensor \vir (w^3)$-module appears in
%$U_{\mathrm{6A}}$ with multiplicity one.
Let $W$ be the subalgebra of $U_{\mathrm{6A}}$ generated by its
weight two subspace. It is clear that $W$ contains both
$U_{\mathrm{3A}}(0)\tensor L(25/28,0)$ and
$U_{\mathrm{3A}}(5/7)\tensor L(25/28,9/7)$. So we only have to show
that $W$ contains an irreducible $\vir (w^1)\tensor \vir
(w^2)\tensor \vir (w^3)$-module of highest weight $(0,1/7,34/7)$.
Let $v$ be a highest weight vector for $\vir (w^1)\tensor \vir (w^2)
\tensor \vir (w^3)$ with highest weight $(0,5/7,9/7)$. Since such a
vector is unique up to linearity, we may take $v$ as follows.
$$
  v=400 \tilde{\omega}(A_1)-16\tilde{\omega}(A_5)+X^3.
$$
Then by the fusion rules for the unitary Virasoro VOAs, we know that
$v_{-2}v$ is contained in the weight 5 subspace of
$$
  L(4/5,0)\tensor L(6/7,5/7)\tensor L(25/28,9/7)
  \oplus
  L(4/5,0)\tensor L(6/7,1/7)\tensor L(25/28,34/7).
$$
Assume that $W$ is a proper subalgebra. Then $W$ is isomorphic to
$U_{\mathrm{3A}}(0)\tensor L(25/28,0)\oplus
U_{\mathrm{3A}}(5/7)\tensor L(25/28,9/7)$ and $v_{-2}v$ is in the
weight 5 subspace of
$$
  L(4/5,0)\tensor L(6/7,5/7)\tensor L(25/28,9/7).
$$
In this case, by computing inner products, we must have the
following equality:
$$
\begin{array}{ll}
  v_{-2}v
  =& \displaystyle
    - \fr{132}{7} w^2_{-1} w^2_0 v
    + \fr{1870}{49} w^2_{-2}v
    - \fr{80}{21} w^3_{-1} w^3_0 v
    + \fr{1128}{49} w^3_{-2} v
    + \fr{288}{23} w^2_0 w^3_{-1} v
  \vsb\\
  & \displaystyle
    - \fr{40}{23} w^2_0 w^3_0 w^3_0 v
    + \fr{440}{19} w^2_{-1} w^3_0 v
    - \fr{198}{19} w^2_0 w^2_0 w^3_0 v.
\end{array}
$$
However, the squared length of the right hand side is
1913313600/437, whereas that of the left hand side is 10209600. This
means that there is a highest weight vector of weight $(0,1/7,34/7)$
in $v_{-2}v$, and we can understand that the difference
$$
  10209600-1913313600/437=2548281600/437
$$
is the squared length of the highest weight vector. Thus the
6A-algebra $U_{\mathrm{6A}}$ is generated by its weight two
subspace.

\section{ Classification for the irreducible modules of the $3C$  and $5A$ algebras}

In this appendix, we shall give the classification of all the
irreducible modules for the $3C$ and $5A$ algebras.
For convenience, we introduce the following notation.
Let $V$ be a VOA and $M$ its module.
For subsets $A\subset V$ and $B\subset M$, we set
$$
  A\cd B=\mathrm{span}\{ a_n v \mid a\in A, v\in B, n\in \Z\} .
$$
It is shown in Lemma 3.12 of \cite{Li15} that for $a,b\in V$, $v \in M$ and
$p,q\in \Z$, there are $m,n\geq 0$ such that
\begin{equation}\label{eq:c1}
  a_p b_q v
  = \sum_{i=0}^m \sum_{j=0}^n \binom{p-n}{i}\binom{n}{j}(a_{p-n-i+j}b)_{q+n+i-j} v.
\end{equation}
In particular, $V\cd (V\cd v)= V\cd v$ for any $v\in M$
so that $V\cd v$ is a submodule of $M$.

Recall the notion of the fusion products (cf.\ \cite{hl5} \cite{li1}).
We denote by $M^1\boxtimes_V M^2$ the fusion product of $V$-modules
$M^1$ and $M^2$.
The basic result (loc.\ cit.) is that the fusion product exists if
$V$ is rational.

We shall study an extension of a rational VOA by an irreducible
module which is not a simple current module.
Our settings are as follows.

Let $V$ be a simple rational VOA and $W$ an irreducible $V$-module such that
\begin{equation}\label{eq:W}
  \dim I_V \binom{M^2}{W\quad M^1} \leq 1
\end{equation}
for any irreducible $V$-modules $M^1$ and $M^2$, where
$I_V\binom{M^2}{W\ M^1}$ denotes the space of $V$-intertwining operators
of type $\binom{M^2}{W\ M^1}$.
Assume that the space $\tilde{V}=V\oplus W$ has a simple VOA structure
$(\tilde{V},\tilde{Y}(\cd,z))$ which is an extension of $V$ such that
for any $u,v\in W$,
\[
  \tilde{Y}(u,z)v= (\mathcal{I}(u,z)+ \mathcal{J}(u,z))v,
\]
where
\[
  \mathcal{I}(\cd,z)\in I_V \binom{V}{W\quad W}  \quad \text{ and }\quad
  \mathcal{J}(\cd,z)\in  I_V \binom{W} {W\quad W}
\]
are non-zero intertwining operators.
Note that the simplicity of $\tilde{V}$ implies that $V$ and $W$ are inequivalent
$V$-modules.
For, if $V$ and $W$ are isomorphic, then the isomorphic image of the vacuum vector
of $V$ in $W$ is a vacuum-like vector (cf.\ \cite{Li}).
Since the vertex operator of a vacuum-like vector commutes with all the vertex
operators on $\tilde{V}$ and $\tilde{V}$ is simple, every vacuum-like vector is
a scalar multiple of the vacuum vector.
Thus $V$ and $W$ are inequivalent.
Then it follows from \eqref{eq:c1} that $\tilde{V}=W\cd W$.
By fixing one VOA structure on $\tilde{V}$, we shall show that
the module structure of certain types of irreducible $\tilde{V}$-modules
are uniquely determined by their $V$-module structures.

\begin{lemma}\label{C1}
  Let $M$ be an irreducible $\tilde{V}$-module.
  Assume that $M$ contains an irreducible $V$-submodule $M^0$ such that
  $M$ is a direct sum $nM^0$ of $n$ copies of $M^0$ as a $V$-module.
  Then $M=M^0$, i.e., $n=1$.
\end{lemma}

\begin{proof}
By \eqref{eq:c1}, we know that $W\cd M^0$ is a $V$-submodule of $M$ which is a
direct sum of some copies of $M^0$.
Since $\tilde{V}$ is simple, it is clear that $W\cd M^0\ne 0$
(cf.\ Proposition 11.9 of \cite{dl}).
By the universal property of the fusion product (cf.\ \cite{hl5} \cite{li1}),
there exists a $V$-epimorphism from $W\boxtimes_V M^0$ onto $W\cd M^0$.
Since $\dim I_V\binom{M^0}{W\ M^0}=1$ by \eqref{eq:W}, $W\boxtimes_V M^0$ contains
a $V$-submodule isomorphic to $M^0$ with multiplicity one.
Therefore, $W\cd M^0$ is an irreducible $V$-submodule isomorphic to $M^0$.
Similarly, $W\cd (W\cd M^0)$ is also an irreducible $V$-submodule of $M$.
If $W\cd M^0=M^0$, then $\tilde{V}\cd M^0= (V\cd M^0)+(W\cd M^0)=M^0$ as a $V$-module
so that $M=M^0$ and we are done.
Assume that $W\cd M^0\ne M^0$.
In this case $M^0+(W\cd M^0)=M^0\oplus (W\cd M^0)$ since both $M^0$ and $W\cd M^0$
are irreducible $V$-submodules.
Then by the irreducibility we have $M=\tilde{V}\cd M^0=M^0\oplus (W\cd M^0)$.
Consider $(W\cd W)\cd M^0$.
As we have seen, $W\cd W=V\oplus W=\tilde{V}$ so that $(W\cd W)\cd M^0=M$.
By the associativity formula
\begin{equation}\label{asso}
  (a_m b)_n=\sum_{i=0}^\infty (-1)^i \binom{m}{i}\{ a_{m-i}b_{n+i}-(-1)^m b_{m+n-i}a_i\},
\end{equation}
we see that $(W\cd W)\cd M^0 \subset W\cd (W\cd M^0)$.
Since $W\cd (W\cd M^0)$ is an irreducible $V$-submodule, so is
$(W\cd W)\cd M^0=\tilde{V}\cd M^0$.
But this is a contradiction.
Hence $W\cd M^0=M^0$ and $M=M^0$.
\end{proof}

\begin{lemma}\label{C2}
  Let $M$ be an irreducible $\tilde{V}$-module which is also irreducible as a $V$-module.
  Then there is exactly one irreducible $\tilde{V}$-module structure on $M$ up to isomorphism.
\end{lemma}

\begin{proof}
Suppose that $(M, Y_M(\cd,z))$ and $(M, Y'_M(\cd,z))$ are two irreducible $\tilde{V}$-module
structures on an irreducible $V$-module $M$.
Without loss, we may assume that $Y_M(u,z)=Y'_M(u,z)$ for all $u\in V$.
Since $M$ is irreducible and $\dim I_V\binom{M}{W\ M}\leq 1$ by \eqref{eq:W},
there exists $\lambda\neq 0$ such that
$$
  Y'_M(a,z)v=\lambda Y_M(a,z)v\quad \text{ for any } a\in W,\ v\in M.
$$
By the associativity, for any $a,b\in W$ and $v\in M$, there exists $k>0$ such
that
\begin{equation}\label{c.1}
\begin{split}
  (z_0+z_2)^k Y_M(a,z_0+z_2)Y_M(b,z_2)v
  =&(z_0+z_2)^k Y_M(Y_{\tilde{V}}(a,z_0)b, z_2)v
  \\
  =& (z_0+z_2)^k Y_M(\mathcal{I}(a,z_0)b+\mathcal{J}(a,z_0)b, z_2)v.
\end{split}
\end{equation}
Similarly, we have
\[
\begin{split}
  (z_0+z_2)^kY'_M(a,z_0+z_2) Y'_M(b,z_2)v
  =& (z_0+z_2)^k Y'_M(Y_{\tilde{V}}(a,z_0)b, z_2)v
  \\
  =& (z_0+z_2)^k Y'_M(\mathcal{I}(a,z_0)b+\mathcal{J}(a,z_0)b, z_2)v
\end{split}
\]
and thus
\begin{equation}\label{c.2}
\begin{split}
  &\lambda^2(z_0+z_2)^kY_M(a,z_0+z_2)Y_M(b,z_2)v
  \\
  =& (z_0+z_2)^k ( Y_M(\mathcal{I}(a,z_0)b,z_2)v
     +\lambda Y_M(\mathcal{J}(a,z_0)b, z_2)v).
\end{split}
\end{equation}
By \eqref{c.1} and \eqref{c.2}, we have
\begin{equation}\label{m.1}
  (z_0+z_2)^k Y_M((\lambda^2-1)\mathcal{I}(a,z_0)b
  +(\lambda^2-\lambda)\mathcal{J}(a,z_0)b, z_2)v=0.
\end{equation}
We note that
$$
  Y_M((\lambda^2-1)\mathcal{I}(a,z_0)b +(\lambda^2-\lambda)\mathcal{J}(a,z_0)b, z_2)v
  \in M((z_0))[[z_2,z_2^{-1}]] .
$$
Thus \eqref{m.1} implies that
$$
  (\lambda^2-1)\mathcal{I}(a,z_0)b=0 \quad \text{ and } \quad
  (\lambda^2-\lambda)\mathcal{J}(a,z_0)b=0.
$$
Therefore, we have $\lambda^2=\lambda=1$.
Hence there is only one irreducible $\tilde{V}$-module structure on $M$.
\end{proof}

\begin{lemma}\label{C3}
  Let $M$ be an irreducible $\tilde{V}$-module.
  Assume that there are two inequivalent irreducible $V$-submodules $M^0$ and $M^1$
  of $M$ such that $M$ is isomorphic to $m M^0\oplus n M^1$ with $m,n>0$ as a $V$-module.
  If $\dim I_V\binom{M^0}{W\ M^0}=0$, then $m=n=1$ and there is exactly one
  $\tilde{V}$-module structure on $M$ up to isomorphism.
\end{lemma}

\begin{proof}
First, we shall show that $W\cd M^0=M^1$ and $M=W\cd M^1=M^0\oplus M^1$.
It is clear from \eqref{eq:c1} that both $W\cd M^0$ and $W\cd M^1$ are $V$-submodules
of $M$ which are non-zero by Proposition 11.9 of \cite{dl}.
By the universal property of the fusion product, we have a $V$-epimorphism from
$W\boxtimes_V M^0$ onto $W\cd M^0$.
By our setting \eqref{eq:W} and the assumption $\dim I_V\binom{M^0}{W\ M^0}=0$,
$W\cd M^0$ is an irreducible $V$-submodule isomorphic to $M^1$.
Thus $M^0+(W\cd M^0)$ is a direct sum.
Since $M$ is an irreducible $\tilde{V}$-module, $M=\tilde{V}\cd M^0
=(V\cd M^0)+ (W\cd M^0)=M^0\oplus (W\cd M^0)$.
This proves $W\cd M^0=M^1$ and $M=M^0\oplus M^1$.
Then $W\cd M^1=W\cd (W\cd M^0)=(W\cd W)\cd M^0=\tilde{V}\cd M^0=M$
by \eqref{eq:c1} and \eqref{asso}.

Let us suppose that $(M, Y_M(\cd,z))$ and $(M, Y'_M(\cd,z))$ are two
irreducible $\tilde{V}$-module structures on the $V$-module $M=M^0\oplus M^1$.
Without loss, we may assume that $Y_M(u,z) =Y'_M(u,z)$ for all $u\in V$.
Let $a\in W$, $v^0\in M^0$ and $v^1\in M^1$ be arbitrary.
As we have shown, $W\cd M^0=M^1$ under both structures.
Thus there is a non-zero $V$-intertwining operator $T(\cd,z)$ of type $\binom{M^1}{W\ M^0}$
and a scalar $\gamma \ne 0$ such that
$$
  Y_M(a,z)v^0=T(a,z)v^0\quad \text{and}\quad Y'_M(a,z)v^0=\gamma T(a,z) v^0.
$$
We have also shown that $W\cd M^1=M^0\oplus M^1$ under both structures.
Thus there exist non-zero $V$-intertwining operators
$I(\cd,z)\in I_V\binom{M^0}{W\ M^1}$, $J(\cd,z)\in I_V\binom{M^1}{W\ M^1}$
and scalars $\mu,\lambda \ne 0$ such that
$$
  Y_M(a,z)v^1= I(a,z)v^1+ J(a,z)v^1 \quad \text{and}\quad
  Y'_M(a,z)v^1= \mu {I}(a,z)v^1 +\lambda {J}(a,z)v^1 .
$$
Then
$$
  Y_M(u+a,z)(v^0+v^1)
  = Y_M(u,z)v^0 +I(a,z)v^1+T(a,z)v^0+(Y_M(u,z)+J(a,z))v^1,
$$
where $u\in V$.
Note that $Y_M(u,z)v^0\in M^0((z))$, $I(a,z)v^1\in M^0((z))$,
$T(a,z)v^0\in M^1((z))$ and $(Y_M(u,z)+J(a,z))v^1\in M^1((z))$.
Hence the vertex operator $Y_M(\cd,z)$  can be written as the following
matrix form.
$$
  Y_M\left(
  \begin{bmatrix}
    u\\ a
  \end{bmatrix},z\right)
  \begin{bmatrix}
    v^0\\ v^1
  \end{bmatrix}
  =
  \begin{bmatrix}
    Y_M(u,z) &  {I}(a,z)\\
    T(a,z) & Y_M(u,z)+ {J}(a,z)
  \end{bmatrix}
  \begin{bmatrix}
    v^0 \\ v^1
  \end{bmatrix} .
$$
Likewise,
$$
  Y'_M\left(
  \begin{bmatrix}
    u\\ a
  \end{bmatrix},z\right)
  \begin{bmatrix}
    v^0\\ v^1
  \end{bmatrix}
  =
  \begin{bmatrix}
    Y_M(u,z) & \mu {I}(a,z)\\
    \gamma T(a,z) & Y_M(u,z)+\lambda {J}(a,z)
  \end{bmatrix}
  \begin{bmatrix}
    v^0\\ v^1
  \end{bmatrix}.
$$

For simplicity, we use the following notation.
For $f(z_0,z_2), g(z_0,z_2)\in M[[z_0,z_0^{-1},z_2,z_2^{-1}]]$,
$f(z_0,z_2)\sim g(z_0,z_2)$ means that there exists $k>0$ such that
$$
  (z_0+z_2)^k f(z_0,z_2)= (z_0+z_2)^k g(z_0,z_2).
$$
Then, by the associativity on $Y'_M(\cd,z)$, for any $a,b \in W$ we have
the following system of equations:
\[
\begin{split}
  &
  \begin{bmatrix}
    Y_M(\mathcal{I}(a,z_0) b,z_2)
    & \mu {I}(\mathcal{J}(a, z_0)b,z_2)
    \\
    \gamma T(\mathcal{J}(a,z_0)b,z_2)
    & Y_M(\mathcal{I}(a,z_0)b,z_2)+\lambda {J}(\mathcal{J}(a,z_0)b,z_2)
  \end{bmatrix}
  \\
  \sim
  &
  \begin{bmatrix}
    \mu\gamma{I}(a,z_0+z_2) T(b,z_2)
    & \mu \lambda{I}(a, z_0+z_2)J(b,z_2)
    \\
    \lambda \gamma {J}(a,z_0+z_2)T(b,z_2)
    & \mu\gamma T(a,z_0+z_2)I(b,z_2) +\lambda^2 {J}(a,z_0+z_2){J}(b,z_2)
\end{bmatrix}
\end{split}
\]
Similarly, by the associativity on $Y_M(\cd,z)$, we have
\[
\begin{split}
  &
  \begin{bmatrix}
    Y_M(\mathcal{I}(a,z_0)b,z_2)
    &  {I}(\mathcal{J}(a, z_0)b,z_2)
    \\
    T(\mathcal{J}(a,z_0)b,z_2)
    & Y_M(\mathcal{I}(a,z_0)b,z_2)+{J}(\mathcal{J}(a,z_0)b,z_2)
  \end{bmatrix}
  \\
  \sim
  &
  \begin{bmatrix}
    {I}(a,z_0+z_2)T(b,z_2)
    & {I}(a, z_0+z_2)J(b,z_2)
    \\
    {J}(a,z_0+z_2)T(b,z_2)
    & T(a,z_0+z_2)I(b,z_2)+{J}(a,z_0+z_2){J}(b,z_2)
  \end{bmatrix}.
\end{split}
\]

Then by a similar argument as in the proof of Lemma \ref{C2},
we have $\mu\gamma=1$, $\mu=\mu\lambda$ and $\lambda\gamma=\gamma$,
from which we conclude that $\lambda=1$ and $\mu=1/\gamma$.
Thus we obtain the following relation:
$$
  Y_M'\left(
  \begin{bmatrix}
    u\\ a
  \end{bmatrix},z\right)
  \begin{bmatrix}
    v^0\\ v^1
  \end{bmatrix}
  =
  \begin{bmatrix}
    Y_M(u,z) &  \fr{1}{\gamma}{I}(a,z)\\
    \gamma T(a,z) & Y_M(u,z)+ {J}(a,z)
  \end{bmatrix}
  \begin{bmatrix}
    v^0 \\ v^1
  \end{bmatrix}.
$$
Now define $\psi: M^0\oplus M^1 \to M^0\oplus M^1$ by
$$
  \psi(v^0)=\frac{1}\gamma v^0\quad \text{ for } v^0\in M^0,\q
  \psi(v^1)= v^1 \quad \text{for } v^1\in M^1.
$$
Then
\[
\begin{split}
  \psi\left( Y_M\left(
  \begin{bmatrix}
    u \\ a
  \end{bmatrix}
  ,z\right)
  \begin{bmatrix}
    v^0 \\ v^1
  \end{bmatrix}
  \right ) =
  &
  \begin{bmatrix}
    \frac{1}\gamma (Y_M(u,z)v^0+{I}(a,z)v^1)
    \\
    T(a,z)v^0 + Y_M(u,z)v^1+ {J}(a,z)v^1
  \end{bmatrix}
  \\
  =&
  \begin{bmatrix}
    Y_M(u,z)\psi(v^0)+\mu{I}(a,z)\psi(v^1)
    \\
    \gamma T(a,z)\psi(v^0) + Y_M(u,z)\psi(v^1)+ {J}(a,z)\psi(v^1)
  \end{bmatrix}
  \\
  =&
   Y'_M\left(
  \begin{bmatrix}
    u \\ a
  \end{bmatrix}
  ,z\right)
  \begin{bmatrix}
    \psi(v^0) \\ \psi(v^1)
  \end{bmatrix}.
\end{split}
\]
Hence $\psi$ induces a $\tilde{V}$-isomorphism from $(M,Y_M(\cd,z))$ to
$(M,Y_M'(\cd,z))$.
\end{proof}

\subsection{Irreducible modules for the $3C$ algebra}

First we shall classify all irreducible modules for the
$3C$-algebra $U_{3C}$.
Recall that
\[
\begin{split}
  U_{3C}
  & \cong  L(\frac{1}2,0)\tensor L(\frac{21}{22},0)
    \oplus L(\frac{1}{2},0)\tensor L(\frac{21}{22},8)
    \\
  & \quad \oplus  L(\frac{1}2,\frac{1}2)\tensor L(\frac{21}{22},\frac{7}{2})
    \oplus L(\frac{1}2,\frac{1}2)\tensor  L(\frac{21}{22},\frac{45}{2})
    \\
  & \quad  \oplus L(\frac{1}{2},\frac{1}{16})\tensor L(\frac{21}{22},\frac{31}{16})
    \oplus L(\frac{1}{2},\frac{1}{16})\tensor L(\frac{21}{22},\frac{175}{16})
\end{split}
\]
as a module of $ L(\frac{1}2,0)\tensor L(\frac{21}{22},0)$.
For simplicity, we shall use $[h_1,h_2]$ to denote the module
$L(\frac{1}2,h_1)\tensor L(\frac{21}{22},h_2)$.

Let $e\in U_{3C}$ be the Virasoro element of the VOA $L(\fr{1}{2},0)$ and
$\tau_e$ the corresponding Miyamoto involution.
Then the corresponding fixed point subalgebra of $\tau_e$ is
\[
  U^{\tau_e}_{3C}= [0,0] \oplus [0,8] \oplus [\frac{1}{2},\frac{7}{2}]
 \oplus [\frac{1}{2},\frac{45}{2}] .
\]

Set $V=[0,0]\oplus [\fr{1}{2},\fr{45}{2}]$ and $W=[0,8]\oplus[\fr{1}{2},\fr{7}{2}]$.
Then $V$ is a subalgebra of $U^{\tau_e}_{3C}$ and $W$ is an irreducible $V$-submodule
of $U^{\tau_e}_{3C}$.
It is shown in the proof of Theorem \ref{thm:3.36} that $W\cd W=V\oplus W=U^{\tau_e}_{3C}$.
Therefore, we can use the representation theory for $V$ to classify irreducible
$U^{\tau_e}_{3C}$-modules.
Note that $V$ is a $\Z_2$-graded simple current extension of $[0,0]$ so that
$V$ is rational.
Moreover, irreducible $V$-modules and their fusion rules are easily determined
(cf.\ \cite{lam,lly,Y}).

\begin{lemma}[cf. \cite{lam,lly,Y}] \label{lem:C4}
  The irreducible modules for $V=[0,0]\oplus [\fr{1}{2},\fr{45}{2}]$ are given as follows.
  $$
  \begin{array}{lll}
    [0, h_{i,j}]\oplus [\frac{1}{2}, h_{i, 12-j}],\
    & [\frac{1}{16}, h_{i,k}]\oplus [\frac{1}{16}, h_{i, 12-k}], \
    & [\frac{1}{16}, h_{i,6}]^+\ \text{and}\quad  [\frac{1}{16}, h_{i, 6}]^-,
  \end{array}
  $$
  where $1\leq i\leq 5$, $j=1,3,5,7,9,11$  and $k=2,4$.
\end{lemma}

Note that all irreducible modules are integrally graded, i.e.,
if $M$ is an irreducible module, then
\[
  M=\bigoplus_{n\in \Z} M_{\lambda+n}\quad \text{ with }\quad
  L(0)|_{M_{\lambda+n}}=\lambda+n
\]
for some $\lambda\in \C$.
Thus, by using the fusion rules of $V$-modules and the integrally graded condition above,
we can list up possible irreducible $U^{\tau_e}_{3C}$-modules as follows.

\begin{lemma}\label{C5}
  Let $M$ be an irreducible $U^{\tau_e}_{3C}$-module.
  Assume that $M$ does not contain $[\fr{1}{16},h_{i,6}]$, $1\leq i\leq 5$, as
  $L(\fr{1}{2},0)\tensor L(\fr{21}{22},0)$-submodules.
  Then as a module of $L(\frac{1}2,0)\tensor L(\frac{21}{22},0)$, $M$ is isomorphic
  to one of the following.
  $$
  \begin{array}{l}
    [0, h_{i,1}]
    \oplus [\frac{1}{2}, h_{i, 5}]
    \oplus [0, h_{i,7}]
    \oplus [\frac{1}{2}, h_{i, 11}],
    \vsb\\
    {} [\frac{1}{2}, h_{i,1}]
      \oplus [0, h_{i, 5}]
      \oplus [\frac{1}{2}, h_{i,7}]
      \oplus [0, h_{i, 11}],
    \vsb\\
    {} [\frac{1}{16}, h_{i,4}]
      \oplus [\frac{1}{16}, h_{i, 8}],
    \end{array}
  $$
  where $1\leq i\leq 5$.
  Moreover, the $U^{\tau_e}_{3C}$-module structure on $M$ is uniquely determined
  by its $V$-module structure.
\end{lemma}

\begin{proof}
Let $M$ be an irreducible $U^{\tau_e}_{3C}$-module.
Since $U^{\tau_e}_{3C}=V\oplus W$ and $V$ is rational,
$M$ is a direct sum of irreducible $V$-modules.
By the list of irreducible $V$-modules shown in Lemma \ref{lem:C4} and
the fusion rules of $\Vir(e)$-modules, we have the following two cases:
As a $\Vir(e)$-module, $M$ is a sum of $L(\fr{1}{2},\fr{1}{16})$
or $M$ does not contain $L(\fr{1}{2},\fr{1}{16})$.
In the former case, by the fusion rules of $V$-modules and the integrally graded
condition, we see that all irreducible $V$-submodules of $M$ are mutually isomorphic
and they are isomorphic to one of $[\fr{1}{16},h_{i,4}]\oplus [\fr{1}{16},h_{i,8}]$
with $1\leq i\leq 5$.
Then by Lemmas \ref{C1} and \ref{C2}, $M$ is in fact irreducible as a $V$-module and
its $U^{\tau_e}_{3C}$-module structure is uniquely determined by its $V$-module
structure.
Hence $M$ is as in the assertion.

If $M$ does not contain $L(\fr{1}{2},\fr{1}{16})$ as $\Vir(e)$-submodules,
then by the list of irreducible $V$-modules shown in Lemma \ref{lem:C4} and
the fusion rules of $V$-modules, we are in the situation as in Lemma \ref{C3}
and the integrally graded condition leads to that $M$ is in the list of the
assertion.
In this case the uniqueness of $U^{\tau_e}_{3C}$-module structure is already shown in
Lemma \ref{C3}.
\end{proof}

Finally, we have the following classification theorem.

\begin{theorem}
  There are exactly five irreducible $U_{3C}$-modules.
  As $L(\fr{1}{2},0)\tensor L(\fr{21}{22},0)$-modules, they are of the
  following form:
  \begin{align*}
    U(0)
    & \cong  [0,0] \oplus [0,8] \oplus [\frac{1}{2},\frac{7}{2}]
      \oplus [\frac{1}{2},\frac{45}{2}] \oplus [\frac{1}{16},\frac{31}{16}] \oplus
      [\frac{1}{16},\frac{175}{16}] \, (=U_{3C}),
    \\
    U(2)
    & \cong [0,\frac{13}{11}] \oplus [0,\frac{35}{11}] \oplus [\frac{1}{2},\frac{15}{22}]
      \oplus [\frac{1}{2},\frac{301}{22}] \oplus [\frac{1}{16},\frac{21}{176}]
      \oplus [\frac{1}{16}, \frac{901}{176}],
    \\
    U(4)
    & \cong [0,\frac{6}{11}]  \oplus [0,\frac{50}{11}] \oplus [\frac{1}{2},\frac{1}{22}]
      \oplus [\frac{1}{2},\frac{155}{22}] \oplus [\frac{1}{16},\frac{85}{176}] \oplus
      [\frac{1}{16}, \frac{261}{176}],
    \\
    U(6)
    & \cong  [0,\frac{1}{11}] \oplus [0,\frac{111}{11}] \oplus [\frac{1}{2},\frac{35}{22}]
      \oplus [\frac{1}{2},\frac{57}{22}] \oplus [\frac{1}{16}, \frac{5}{176}]
      \oplus [\frac{1}{16},\frac{533}{176}],
    \\
    U(8)
    & \cong [0,\frac{20}{11}] \oplus [0,\frac{196}{11}] \oplus [\frac{1}{2},\frac{7}{22}]
      \oplus [\frac{1}{2},\frac{117}{22}] \oplus [\frac{1}{16},\frac{133}{176}] \oplus
      [\frac{1}{16},\frac{1365}{176}] .
  \end{align*}
\end{theorem}

\begin{proof}
Set $U^-_{3C}=[\fr{1}{16},\fr{31}{16}]\oplus [\fr{1}{16},\fr{175}{16}]$.
Then $\tau_e\in \aut(U_{3C})$ acts on $U^-_{3C}$ as $-1$ and
$U_{3C}=U^{\tau_e}_{3C}\oplus U_{3C}^-$ is a $\Z_2$-graded extension of $U^{\tau_e}_{3C}$.
Let $M$ be an irreducible $U_{3C}$-module.
Denote by $M^0$ the sum of irreducible $\vir(e)$-submodules of $M$ isomorphic to
$L(\fr{1}{2},0)$
%%!!!%%%
or
%%!!!%%%
$L(\fr{1}{2},\fr{1}{2})$ and by $M^1$ the sum of
irreducible $\vir(e)$-submodules of $M$ isomorphic to $L(\fr{1}{2},\fr{1}{16})$.
Then $M=M^0\oplus M^1$.
By the fusion rules of $L(\fr{1}{2},0)$-modules, $M^0$ and $M^1$ are inequivalent
irreducible $U_{3C}^{\tau_e}$-submodules and $M$ carries the following $\Z_2$-grading:
$U^-_{3C}\cd M^0=M^1$ and $U^-_{3C}\cd M^1=M^0$.
In Lemma \ref{C5} we have classified irreducible $U^{\tau_e}_{3C}$-modules having
no $\Vir(e)$-submodules isomorphic to $L(\fr{1}{2},\fr{1}{16})$.
Then by the integrally graded condition, we see that $M$ cannot contain
$[\fr{1}{16},h_{i,6}]$, $1\leq i\leq 5$ as
$L(\fr{1}{2},0)\tensor L(\fr{21}{22},0)$-submodules.
Now thanks to Lemma \ref{C5} we can classify the possible pairs $(M^0,M^1)$ of
irreducible $U_{3C}^{\tau_e}$-modules such that $M^0\oplus M^1$ are integrally graded.
As a result, we see that $M$ is isomorphic to one of $U(2i)$, $i=0,1,\dots,4$,
as an $L(\fr{1}{2},0)\tensor L(\fr{21}{22},0)$-module.
We already know that all $U(2i)$, $i=0,1,\dots,4$, appear as $U_{3C}$-submodules of
$V_{\sqrt{2}E_8}$.
Thus there exist irreducible $U_{3C}$-modules of the form $U(2i)$.
It remains to show that there is only one irreducible $U_{3C}$-module structure
on each $U(2i)$.
Let $U(2i)=M^0\oplus M^1$ be the decomposition above.
It is clear that both $M^0$ and $M^1$ are self-dual $U^{\tau_e}_{3C}$-modules.
By case by case verifications, we can deduce the fusion rule
$\dim I_{U_{3C}^{\tau_e}}\binom{M^1}{U^-_{3C}\ M^0}=1$
from the fusion rules of $L(\fr{1}{2},0)\tensor L(\fr{21}{22},0)$-modules.
Since both $M^0$ and $M^1$ are self-dual, we also have
$\dim I_{U_{3C}^{\tau_e}} \binom{M^0}{U^-_{3C}\ M^1}=1$.
Then by a standard argument as in \cite{lam,lly,Y} we can easily show that there exists
only one irreducible $U_{3C}$-module structure on $U(2i)=M^0\oplus M^1$, since
$U_{3C}=U_{3C}^{\tau_e}\oplus U_{3C}^-$ is a $\Z_2$-graded extension of $U^{\tau_e}_{3C}$.
\end{proof}

\subsection{Irreducible modules for the $5A$ algebra}

Next we shall classify the irreducible modules for the $5A$-algebra $U_{5A}$.
The argument used here is exactly the same as in the case of $U_{3C}$ so that
we omit details in this case.
Recall that
\[
  \begin{array}{l}
    U_{5A} \cong L(\frac{1}{2},0)\tensor L(\frac{25}{28},0)\tensor L(\frac{25}{28},0)
      \oplus L(\frac{1}{2},\frac{1}{16})\tensor
        L(\frac{25}{28},\frac{5}{32})\tensor
        L(\frac{25}{28},\frac{57}{32})
    \vsb\\
      \qquad \oplus L(\frac{1}{2},\frac{1}{16})\tensor
        L(\frac{25}{28},\frac{57}{32})\tensor
        L(\frac{25}{28},\frac{5}{32})
      \oplus L(\frac{1}{2},\frac{1}{2})\tensor
        L(\frac{25}{28},\frac{3}{4})\tensor L(\frac{25}{28},\frac{3}{4})
    \vsb\\
      \qquad \oplus L(\frac{1}{2},0)\tensor L(\frac{25}{28},\frac{3}{4})\tensor
        L(\frac{25}{28},\frac{13}{4})
      \oplus L(\frac{1}{2},0)\tensor L(\frac{25}{28},\frac{13}{4})\tensor
        L(\frac{25}{28},\frac{3}{4})
    \vsb\\
      \qquad \oplus L(\frac{1}{2},\frac{1}{16})\tensor
        L(\frac{25}{28},\frac{57}{32})\tensor
        L(\frac{25}{28},\frac{165}{32})
      \oplus L(\frac{1}{2},\frac{1}{16})\tensor
        L(\frac{25}{28},\frac{165}{32})\tensor
        L(\frac{25}{28},\frac{57}{32})
    \vsb\\
      \qquad \oplus L(\frac{1}{2},\frac{1}{2})\tensor
        L(\frac{25}{28},\frac{13}{4})\tensor L(\frac{25}{28},\frac{13}{4})
      \oplus L(\frac{1}{2},\frac{1}{2})\tensor L(\frac{25}{28},0)\tensor
        L(\frac{25}{28},\frac{15}{2})
    \vsb\\
      \qquad \oplus L(\frac{1}{2},\frac{1}{2})\tensor
        L(\frac{25}{28},\frac{15}{2})\tensor L(\frac{25}{28},0)
      \oplus L(\frac{1}{2},0)\tensor L(\frac{25}{28},\frac{15}{2})\tensor
        L(\frac{25}{28},\frac{15}{2}).
  \end{array}
\]

Again, we shall use $e\in U_{5A}$ to denote the Virasoro element of the
VOA $L(1/2,0)$ and $\tau_e$ the corresponding Miyamoto involution.
We shall also use $[h_1,h_2,h_3]$ to denote the module of the form
$L(\frac{1}{2},h_1)\tensor L(\frac{25}{28},h_2)\tensor  L(\frac{25}{28},h_3)$.

Note that the fixed point subalgebra of $\tau_e$ is as follows.
\[
\begin{split}
  U^{\tau_e}_{5A} \cong
  & [0,0,0]\oplus [0,\fr{15}{2},\fr{15}{2}]\oplus [\fr{1}{2},\fr{15}{2},0]
    \oplus [\fr{1}{2},0,\fr{15}{2}]
    \\
  & \oplus  [0,\fr{3}{4},\fr{13}{4}]\oplus [0,\fr{13}{4},\fr{3}{4}]\oplus
    [\fr{1}{2},\fr{13}{4},\fr{13}{4}] \oplus [\fr{1}{2},\fr{3}{4},\fr{3}{4}] .
\end{split}
\]
It contains a subalgebra
\[
  V=[0,0,0]\oplus  [0,\fr{15}{2},\fr{15}{2}]\oplus
  [\fr{1}{2},\fr{15}{2},0]\oplus [\fr{1}{2},0,\fr{15}{2}],
\]
which is a $(\Z_2\times\Z_2)$-graded simple current extension of $[0,0,0]$.
Set
$$
   W =[0,\fr{3}{4}, \fr{13}{4}]\oplus [0,\fr{13}{4},\fr{3}{4}]\oplus
   [\fr{1}{2},\fr{13}{4},\fr{13}{4}] \oplus [\fr{1}{2},\fr{3}{4},\fr{3}{4}].
$$
Then $W$ is an irreducible $V$-submodule of $U_{5A}^{\tau_e}$ and we have
a decomposition $U_{5A}^{\tau_e}=V\oplus W$.
It is shown in Lemma \ref{lem:3.16} that $W\cd W=V\oplus W=U_{5A}^{\tau_e}$.
Therefore, we can use the representation theory of $V$ to classify irreducible
$U_{5A}$-modules.

\begin{lemma}[cf. \cite{lam,lly,Y}]
  The irreducible modules for $V$ are given as follows.
  $$
  \begin{array}{l}
  [0,h_{i,m}, h_{k,n}]
    \oplus [\frac{1}{2}, h_{i, 8-m}, h_{k,n}]
    \oplus [\frac{1}2,h_{i,m}, h_{k,8-n}]
    \oplus [0, h_{i, 8-m}, h_{k,8-n}],
  \vsb\\
  {}[\frac{1}{16}, h_{i,j}, h_{k,\ell}]
    \oplus [\frac{1}{16}, h_{i,8-j}, h_{k,\ell}]
    \oplus [\frac{1}{16}, h_{i,j}, h_{k,8-\ell}]
    \oplus [\frac{1}{16}, h_{i,8-j}, h_{k,8-\ell}],
  \vsb\\
  {}([\frac{1}{16}, h_{i,4}, h_{k,\ell}]
    \oplus [\frac{1}{16}, h_{i,4},h_{k,8-\ell}])^\pm
    \quad \text{or} \quad
    ([\frac{1}{16}, h_{i,j}, h_{k,4}]
    \oplus [\frac{1}{16}, h_{i,8-j}, h_{k,4}])^\pm,
  \vsb\\
  {}[\frac{1}{16}, h_{i,4},h_{k,4}]\tensor Q,
  \end{array}
  $$
  where $1\leq i,k\leq 3$, $m,n=1,3,5,7$, $j,\ell=2,6$ and
  $Q$ is the unique 2-dimensional irreducible module of the quaternion group of order 8.
\end{lemma}

Now by using the fusion rules and the fact that the weights of
irreducible modules are integrally graded, we have

\begin{lemma}\label{lem:C.8}
  Let $M$ be an irreducible module for $U_{5A}^{\tau_e}$.
  Assume that $M$ does not contain $[\fr{1}{16},h_{p,4},h_{q,4}]$,   $1\leq p,q\leq 3$,
  as $L(\fr{1}{2},0)\tensor L(\fr{25}{28},0)\tensor L(\fr{25}{28},0)$-submodules.
  Then, as an module of $L(\frac{1}{2},0)\tensor L(\frac{25}{28},0)\tensor
  L(\frac{25}{28},0)$,  $M$ is isomorphic to one of the following,
  where $i, k=1,2$ or $3$.
  $$
  \begin{array}{l}
    [0,h_{i,1}, h_{k,1}]
      \oplus [\frac{1}{2}, h_{i, 7}, h_{k,1}]
      \oplus [\frac{1}2,h_{i,1}, h_{k,7}]
      \oplus [0, h_{i, 7}, h_{k,7}]
    \vsb\\
      \quad
      \oplus [0,h_{i,3}, h_{k,5}]
      \oplus [\frac{1}{2}, h_{i, 3},  h_{k,3}]
      \oplus [\frac{1}2,h_{i,5}, h_{k,5}]
      \oplus [0, h_{i, 5},h_{k,3}],
    \vsv\\
    {}[0,h_{i,1}, h_{k,7}]
      \oplus [\frac{1}{2}, h_{i, 1}, h_{k,1}]
      \oplus [\frac{1}2,h_{i,7}, h_{k,7}]
      \oplus [0, h_{i, 7}, h_{k,1}]
    \vsb\\
      \quad
      \oplus [0,h_{i,3}, h_{k,3}]
      \oplus [\frac{1}{2}, h_{i, 3},h_{k,5}]
      \oplus [\frac{1}2,h_{i,5}, h_{k,3}]
      \oplus [0, h_{i, 5},h_{k,5}],
    \vsv\\
    {}[\frac{1}{16}, h_{i,2}, h_{k,6}]
      \oplus [\frac{1}{16}, h_{i,6}, h_{k,6}]
      \oplus [\frac{1}{16}, h_{i,6}, h_{k,2}]
      \oplus [\frac{1}{16},h_{i,2}, h_{k,2}],
    \vsv\\
    {}[\frac{1}{16}, h_{i,4}, h_{k,2}]
      \oplus [\frac{1}{16}, h_{i,4},h_{k,6}]
      \oplus \ [\frac{1}{16}, h_{i,6}, h_{k,4}]
      \oplus [\frac{1}{16}, h_{i,2}, h_{k,4}].
      \vsb
  \end{array}
  $$
  Moreover, the $U_{5A}^{\tau_e}$-module structure on $M$ is uniquely determined by
  its $V$-module structure.
\end{lemma}

Finally by using the same method for the $3C$ case, we have the
following theorem.

\begin{theorem}
  There are exactly nine irreducible modules $U(i,j),i,j=1,3,5$ for   $U_{5A}$.
  As $L(\fr{1}{2},0) \otimes L(\fr{25}{28},0)\otimes L(\fr{25}{28},0)$-modules,
  they are of the following form.
  \begin{equation*}
  \begin{split}
    U(i,j)
    & \cong [0, h_{i,1}, h_{j,1}]\oplus [0, h_{i,3}, h_{j,5}]
    \oplus [0, h_{i,5}, h_{j,3}]\oplus [0, h_{i,7}, h_{j,7}]
    \\
    & \quad \oplus [\frac{1}2, h_{i,1}, h_{j,7}]\oplus [\frac{1}2,h_{i,3}, h_{j,3}]
    \oplus [\frac{1}2, h_{i,5}, h_{j,5}]\oplus [\frac{1}2, h_{i,7}, h_{j,1}]
    \\
    & \quad \oplus [\frac{1}{16}, h_{i,2}, h_{j,4}]\oplus [\frac{1}{16}, h_{i,4}, h_{j,2}]
      \oplus [\frac{1}{16}, h_{i,6},h_{j,4}]\oplus [\frac{1}{16}, h_{i,4}, h_{j,6}].
  \end{split}
\end{equation*}
\end{theorem}

\begin{proof}
Let $M$ be an irreducible $U_{5A}$-module.
By the integrally graded condition, we know that $M$ cannot contain
$[\fr{1}{16},h_{i,4},h_{k,4}]$, $1\leq i,k\leq 3$, as
$L(\fr{1}{2},0)\tensor L(\fr{25}{28},0)\tensor L(\fr{25}{28},0)$-submodules.
Then by Lemma \ref{lem:C.8} and the integrally graded condition, $M$ is isomorphic to
one of $U(i,j)$, $i,j=1,3,5$, as an
$L(\fr{1}{2},0)\tensor L(\fr{25}{28},0)\tensor L(\fr{25}{28},0)$-module.
We already know that all $U(i,j)$, $i,j=1,3,5$, appear as $U_{5A}$-submodules of
$V_{\sqrt{2}E_8}$.
Thus there exist irreducible $U_{5A}$-modules of the form $U(i,j)$.
It remains to show the uniqueness of the irreducible $U_{5A}$-module structure on
each $U(i,j)$.
By a similar argument as in the case of $U_{3C}$-modules, we can establish the
uniqueness.
\end{proof}

\bibliographystyle{amsplain}

\begin{thebibliography}{99}

\bibitem{abe} T. Abe, Fusion rules for the charge conjugate orbifold, \emph{J. Algebra}
\textbf{242} (2001), 624--655.

\bibitem{c1}  J. H. Conway, A simple construction for the Fisher-Griess
Monster group, \emph{Invent. Math}. \textbf{79} (1985), 513-540.

\bibitem{cs} J. H.~Conway and N.J.A.~Sloane, \emph{Sphere Packings, Lattices and
Groups}, Springer-Verlag, 1988.

\bibitem{d1}  C. Dong, Vertex algebras associated with even lattices,\emph{ J.
Algebra} \textbf{161} (1993), 245-265.

%\bibitem{dong}  C. Dong, The representation of Moonshine module vertex
%operator algebra, \emph{Contemporary Math.}, \textbf{175} (1994), 27-36.

\bibitem{dg} C. Dong and  R. L. Griess Jr.,
Rank One Lattice Type Vertex Operator Algebras and Their
Automorphism Groups,  \emph{J. Algebra} \textbf{208},1998,
262--275.

%\bibitem{dgh}  C. Dong, R. L. Griess Jr. and G. H\"{o}hn, Framed vertex
%operator algebras, codes and the moonshine module, \emph{Comm.
%Math. Phys.} 193 No.2(1998), pp. 407-448.


\bibitem{dl}  C. Dong and J. Lepowsky, \emph{Generalized vertex algebras and
relative vertex operators}, Progress in Math. Vol. 112,
Birkh\"{a}user, Boston 1993.

%\bibitem{dlm1}  C. Dong, H. Li and G. Mason, Regularity of rational vertex
%operator algebras, \emph{Adv. Math.} 132 (1997), 148-166.

%\bibitem{dlm2}  C. Dong, H. Li and G. Mason, Simple currents and extensions
%of vertex operator algebras, \emph{Comm. Math. Phys.}, 180(1996),
%671-707.

%\bibitem{dlm3} C. Dong, H. Li and G. Mason, Some twisted
%sectors for the Moonshine module, \emph{Moonshine, the Monster,
%and related topics (South Hadley, MA, 1994)}, Contemp. Math., 193,
%Amer. Math. Soc., Providence, RI, 1996, 25 -- 43.

%\bibitem{dlin}  C. Dong and Z. Lin, Induced modules for vertex operator
%algebras, \emph{Comm. Math. Phys.}, 179 (1996), 154-184.

\bibitem{dlmn}  C. Dong, H. Li, G. Mason and S.P. Norton, Associative
subalgebras of Griess algebra and related topics, \emph{Proc. of
the Conference on the Monster and Lie algebra at the Ohio State
University}, May 1996, ed. by J. Ferrar and K. Harada, Walter de
Gruyter, Berlin - New York, 1998, 27--42.

\bibitem{dly} C. Dong, C. Lam and H. Yamada, Decomposition of the vertex
operator algebra $V_{\sqrt{2}A_3}$, \emph{J. Algebra} 222 (1999),
500-510.

\bibitem{dly2} C. Dong, C. Lam, and H. Yamada, Decomposition of the
vertex operator algebra $V\sb {\sqrt{2}D\sb l}$.  \emph{Commun.
Contemp. Math.}  3  (2001),  no. 1, 137--151.

\bibitem{dm}  C. Dong and G. Mason, On Quantum Galois Theory, \emph{Duke Math.
J.} 86 (1997), no. 2, 305-321.

\bibitem{dm-z3}  C. Dong and G. Mason, The construction of the moonshine
module as a $Z\sb p$-orbifold. Mathematical aspects of conformal
and topological field theories and quantum groups (South Hadley,
MA, 1992), Contemp. Math., 175, Amer. Math. Soc., Providence, RI,
1994, 37--52.

\bibitem{dmz}  C. Dong, G. Mason and Y. Zhu, Discrete series of the Virasoro
algebra and the moonshine module, \emph{Pro. Symp. Pure. Math.},
American Math. Soc. 56 II (1994), 295-316.

\bibitem{ff}  B. L. Feigin and D. B. Fuchs, Verma modules over the Virasoro
algebra, \emph{Topology} (Leningrad, 1982), Lecture Notes in
Math., {\bf 1060}, 230--245, Springer, Berlin-New York, 1984.

\bibitem{FFR}
A. J. Feingold, I. B. Frenkel and J. F.X. Ries, Spinor
construction of vertex operator algebras, triality, and
$E_8^{(1)}$, \emph{Contemp. Math.} {\bf 121} (1991).

\bibitem{FRW}
A. J. Feingold, J. F.X. Ries and M. Weiner, Spinor construction of
the $c=\frac12$ minimal model. Moonshine, the Monster, and related
topics (South Hadley, MA, 1994), \emph{Contemp. Math.} {\bf 193},
Amer. Math. Soc., Providence, RI, 1996, 45--92.

\bibitem{fhl}  I. B.\ Frenkel, Y. Huang and J. Lepowsky, \emph{On axiomatic
approaches to vertex operator algebras and modules}, Mem. Amer.
Math. Soc. 104, 1993.

\bibitem{flm}  I. B.\ Frenkel, J. Lepowsky, and A. Meurman,\emph{Vertex Operator
Algebras and the Monster}, Pure and Applied Math., Vol. 134,
Academic Press, 1988.

\bibitem{fz}  I. B. Frenkel and Y. Zhu, Vertex operator algebras associated
to representations of affine and Virasoro algebras, \emph{Duke
Math. J.} 66(1992), 123-168.

\bibitem{GN} G. Glauberman and S. P. Norton, On McKay's connection between
the affine $E_8$ diagram and the Monster, \emph{CRM Proceedings
and Lecture Notes}, Vol. \textbf{30}, Amer. Math. Soc.,
Providence, 2001, 37--42.

\bibitem{griess} R. Griess, The Friendly Giant, \emph{ Invent. Math.}
{\bf 69} (1982), 1-102.

\bibitem{gko} P. Goddard, A. Kent and D. Olive, Unitary
representations of the Virasoro and super-Virasoro algebras, \emph
{Comm. Math. Phys.} \textbf{103} (1986), 105-119.

%\bibitem{hk} M. Harada and M. Kitazume, $\Z_4$-Code Constructions
%for the Niemeier Lattices and their Embeddings in the Leech
%Lattice, \emph{Europ. J. Combinatorics } (2000) 21, 473--485.



%\bibitem{ho} G. H\"ohn, Selbstduale Vertexoperatorsuperalgebren und das
%Babymonster. Ph.D. thesis, University of Bonn, 1995.

%\bibitem{h1}
%Y.-Z. Huang, A theory of tensor products for module categories for
%a vertex operator algebra, IV, \emph{ J. Pure Appl. Alg.} {\bf
%100} (1995), 173-216.

%\bibitem{h2-5}
%Y.-Z. Huang, A nonmeromorphic extension of the moonshine module
%vertex operator algebra, in: \emph{ Moonshine, the Monster and
%related topics, Proc. Joint Summer Research Conference, Mount
%Holyoke, 1994,} ed. C. Dong and G. Mason, Contemporary Math., Vol.
%193, Amer. Math. Soc., Providence, 1996, 123--148.

%\bibitem{h2}
%Y.-Z. Huang, Virasoro vertex operator algebras, (nonmeromorphic)
%operator product expansion and the tensor product theory, \emph{
%J. Algebra} {\bf 182} (1996), 201--234.


%\bibitem{hl1} Y.-Z. Huang and J. Lepowsky, Toward a theory of
%tensor products for representations of a vertex operator algebras,
%\emph{ Proc. 20th International Conference on Differential
%Geometric Methods in Theoretical Physics, New York, 1991}, ed. S.
%Catto and A. Rocha, World Scientific, Singapore, 1992, 344--354.

%\bibitem{hl2}
%Y.-Z. Huang and J. Lepowsky, A theory of tensor products for
%module categories for a vertex operator algebra, I, \emph{ Selecta
%Mathematica, New Series} {\bf 1} (1995), 699-756.

%\bibitem{hl3}
%Y.-Z. Huang and J. Lepowsky, A theory of tensor products for
%module categories for a vertex operator algebra, II, \emph{
%Selecta Mathematica, New Series} {\bf 1} (1995), 757--786.

%\bibitem{hl4}
%Y.-Z. Huang and J. Lepowsky, Tensor products of modules for a
%vertex operator algebras and vertex tensor categories, in:
%     \emph{ Lie Theory and Geometry,
%in honor of Bertram Kostant,} ed. R. Brylinski, J.-L. Brylinski,
%V. Guillemin, V. Kac, Birkh\"{a}user, Boston, 1994, 349--383.

\bibitem{hl5}
Y.-Z. Huang and J. Lepowsky, A theory of tensor products for
module categories for a vertex operator algebra, III, \emph{ J.
Pure Appl. Alg.} {\bf 100} (1995),  141-171.

%\bibitem{kw}  V. Kac and W. Wang, Vertex operator superalgebras and
%representations, \emph{Contemporary Math.} 175 (1994), 161-191.

\bibitem{kr} V. Kac and A. K. Raina, \emph{Bombary Lectures on Highest weight
representations of infinite dimensional Lie algebra},  Adv. Ser.
Math. Phys. Vol. 2, World Scientific, 1987.

\bibitem{kly} M. Kitazume, C. Lam and H. Yamada, Decomposition of the
Moonshine vertex operator algebra as Virasoro modules, \emph{J.
Algebra}, \textbf{226} (2000), 893-919.

\bibitem{kly3}  M. Kitazume, C.H. Lam, and H. Yamada, $3$-state Potts model,
Moonshine vertex operator algebra and 3A-elements of the Monster
group, \emph{IMRN}, 2003, No. 23, 1269 - 1303.


\bibitem{kmy}  M. Kitazume, M. Miyamoto and H. Yamada, Ternary codes and
vertex operator algebras, \emph{J. Algebra}, \textbf{ 223} (2000),
379-395.


\bibitem{lam}  C. Lam, Induced modules for orbifold vertex operator
algebras, \emph{J. Math. Soc. of Japan}, 53 (2001), no. 3,
541-557.

\bibitem{lam-3c} C. Lam, Lattice vertex operator
algebra $V_{\sqrt{2}E_8}$ and an algebra of Miyamoto of central
charge $\frac{1}2+\frac{21}{22}$, RIMS Kokyuroku 1327, Kyoto
University,  Japan, June 2003, 159--169.

\bibitem{lly}  C.H. Lam, N. Lam, and H. Yamauchi, Extension of unitary
Virasoro vertex operator algebra by a simple module, \emph{IMRN}
2003, No. 11, 577 - 611.

\bibitem{ly}  C. Lam and H. Yamada, $\mathbb{Z}_{2}\times
\mathbb{Z}_{2}$ codes and vertex operator algebras, \emph{J.
Algebra} \textbf{224} (2000), 268-291.


\bibitem{ly2}  C. Lam and H. Yamada, Tricrtical $3$-state Potts model and vertex operator
algebras constructed from ternary codes, Comm. Algebra,
\textbf{32} (2004), 4197-4220.

\bibitem{ly3}  C. Lam and H. Yamada, Decomposition of the lattice
vertex operator algebra $V_{\sqrt{2}A_l}$, \emph{J. Algebra},
\textbf{272} (2004),  614-624.

\bibitem{lyy}  C. Lam, H. Yamada and H. Yamauchi,
Vertex operator algebras, extended $E_8$ diagram, and McKay's
observation on the Monster simple group, to appear in Trans. Amer.
Math. Soc.

%\bibitem{lm}  J. Lepowsky and A. Meurman, An $E_8$-approach to
%the Leech lattice and the Conway group, \emph{J. Algebra} 77
%(1982),  484--504.

%\bibitem{kmy}  M. Kitazume, M. Miyamoto and H. Yamada, Ternary codes and
%vertex operator algebras, preprint.

%\bibitem{lam}  C. Lam, Twisted representations of code vertex operator
%algebras, \emph{J. Algebra} 217 (1999), 275-299.

%\bibitem{li}  H. Li, Local systems of vertex operators, vertex superalgebras
%and modules, \emph{J. Pure Appl. Algebra} 109 (1996), no. 2,
%143--195.

\bibitem{Li}  H. Li, Symmetric invariant bilinear forms on vertex
operator algebras, \emph{J. Pure and Appl. Algbera} \textbf{96}
(1994), 279--297.

\bibitem{li1}
  H. Li, An analogue of the Hom functor and a generalized nuclear democracy
  theorem, {\it Duke Math. J.} {\bf 93} (1998), 73--114.


%\bibitem{li2}  H. Li, Extension of vertex operator algebras by a self-dual
%simple module, \emph{J. Algebra} 187 (1997), 236-267.

\bibitem{Li15}
  H. Li, The regular representation, Zhu's $A(V)$-theory and induced modules,
  {\it J. Algebra} {\bf 238} (2001), 159--193.

\bibitem{mat} A. Matsuo, Norton's trace formulae for the Griess algebra
 of a vertex operator agebra with larger symmetry,\emph{ Comm. Math. Phys.}
\textbf{ 224} (2001), 565--591.

\bibitem{McK} J. McKay, Graphs, singularities, and finite groups,
\emph{Proc. Symp. Pure Math.}, Vol. \textbf{37}, Amer. Math. Soc.,
Providence, 1980, 183--186.


\bibitem{m1}  M. Miyamoto, Griess algebras and conformal vectors in vertex
operator algebras, \emph{J. Algebra} 179 (1996), 523-548.

\bibitem{m2}  M. Miyamoto, Binary codes and vertex operator (super)algebras,
\emph{J. Algebra} 181 (1996), 207-222.

\bibitem{m5}  M. Miyamoto, $3$-state Potts model and automorphism of vertex
operator algebra of order $3$, \emph{J. Algebra} \textbf{239}
(2001), 56--76.

\bibitem{m6}  M. Miyamoto, VOAs generated by two conformal vectors whose
$\tau$-involutions generate $S_3$,  \emph{J. Algebra},
\textbf{268} (2003), no. 2, 653--671. .

\bibitem{m4}  M. Miyamoto, A new construction of the Moonshine
vertex operator algebras over the real number field, \emph{Ann. of
Math.}, \textbf{159} (2004), no. 2, 535--596.

%\bibitem{mo}  G. Mossberg, Axiomatic vertex algebras and the Jacobi
%identity, \emph{J. Algebra} 170 (1994), 956-1010.

%\bibitem{mse} G. Moore and N. Seiberg, Classical and quatumn
%conformal field theory, \emph{Comm. Math. Phys.} \textbf{123}
%(1989), 177-254.

\bibitem{sy} S. Sakuma and H. Yamauchi, Vertex operator algebra with automorphism
group $S_3$, \emph{J. Algebra} \textbf{267} (2003), 272--297.

\bibitem{Shi} H. Shimakura, Decompositions of the moonshine module with
respect to subVOAs associated to codes over $\Z_{2k}$, \emph{J.
Algebra} \textbf{251} (2002), 308--322.

\bibitem{TY}
  K. Tanabe and H. Yamada, The fixed point subalgebra of a lattice
  vertex operator algebra by an automorphism of order three, preprint.

\bibitem{Y}
  H. Yamauchi, Module category of simple current extensions of vertex
  operator algebras, {\it J. Pure Appl. Algebra} {\bf 189} (2004) 315--328.

\bibitem{W}
W. Wang, Rationality of Virasoro vertex operator algebras, {\it
Duke Math. J. IMRN}, {\bf Vol. 71}, No. 1 (1993), 197-211.

%\bibitem{zhu}  Y. Zhu, Vertex operator algebras, elliptic functions and
%modular form, Ph.D. dissertation, Yale University, 1990.

\bibitem{zf1} A.B. Zamolodchikov and V.A. Fateev, Nonlocal
(parafermion) currents in two dimensional conformal quantum field
theory and self-dual critical points in $\Z_N$-symmetric
statistical systems, \emph{Sov. Phys. JETP} \textbf{62}(1985),
215-225.


\end{thebibliography}

\end{document}